%% file: main.tex
\documentclass[12pt]{amsart}
\usepackage{amssymb}
\input amssym.def
\usepackage{pstricks}
\newtheorem{theorem}{Theorem}[section]
\newtheorem{lemma}[theorem]{Lemma}
\newtheorem{prop}[theorem]{Proposition}
\newtheorem{cor}[theorem]{Corollary}
\theoremstyle{definition}
\newtheorem{rem}[theorem]{Remark}
\newcommand{\foorp}{\hfill$\Box$}
\newcommand{\M}{\mathcal{M}}
\newcommand{\PM}{\mathcal{PM}}
\newcommand{\Z}{\mathbb{Z}}
\newcommand{\R}{\mathbb{R}}
\newcommand{\T}{\mathcal{T}}
\newcommand{\Dif}{\textrm{Diff}}
\newcommand{\bdr}[1]{\partial\! #1}
\newcommand{\C}{\mathcal{C}}
\newcommand{\Xo}{\mathcal{C}^{\mathrm{ord}}(F)}

\newcommand{\p}[1]{p(\left<#1\right>)}
\newcommand{\lr}[1]{\left<#1\right>}
\numberwithin{equation}{section}
\author{B\l{}a\.zej Szepietowski}
\title[Presentation for the mapping class group]
{A presentation for the mapping class group of a non-orientable surface from the action on the complex of curves}
\address[]{Institute of Mathematics, Gda\'nsk University, Wita Stwosza 57,
80-952 Gda\'nsk, Poland}
\email{blaszep@math.univ.gda.pl}
\thanks{Supported by KBN 1 P03A 024 26}
\begin{document}
\begin{abstract}
We study the action of the mapping class group $\M(F)$ on the
complex of curves of a non-orientable surface $F$. Following the
outline of \cite{B} we obtain, using the result of \cite{Br}, a
presentation for $\M(F)$ defined in terms of the mapping class
groups of the complementary surfaces of collections of curves,
provided that $F$ is not sporadic, i.e. the complex of curves of
$F$ is simply connected. We also compute a  finite presentation
for the mapping class group of each sporadic surface.
\end{abstract}
%
%\noindent{\small Accepted for publication in\\
%Osaka Journal of Mathematics.}
\maketitle
\input{intro}
\input{preli}
\input{stab}
\input{orbits}
\input{present}
\input{sporadic}

%\bigskip

%\noindent{\bf Acknowledgment.} Author wishes to thank the referee
%for his/her helpful suggestions.
%
%
%%%%%%%%%%%%%%%%%%%%%%%%%%%%%%%%%%%%%%%%%%%%%%%%%%%%%%%
%%%  Bibliography  %%%
%%%%%%%%%%%%%%%%%%%%%%%

%

\end{document}

%% file: intro.tex
\section{\label{intro}Introduction}

Presentations for the mapping class group $\M(F_g^n)$ of a compact
orientable surface of genus $g$ with $n$ boundary components have
been found by various authors. Hatcher and Thurston \cite{HT}
derived a presentation for $\M(F_g^1)$ from its action on a simply
connected 2-dimensional complex, the {\it cut system complex}. This
complex was simplified by Harer \cite{Har1} and using this
simplified complex, Wajnryb \cite{W} obtained a simple presentation
for $\M(F_g^1)$ and $\M(F_g^0)$. Starting from Wajnryb's result,
Gervais \cite{G} found a simple presentation for $\M(F_g^n)$ for any
$n$ and $g\ge 1$. Benvenuti \cite{B} and Hirose \cite{Hir} showed
independently how the Gervais presentation can be recovered using
two different modifications of the classical {complex of curves}
introduced by Harvey \cite{H}. Benvenuti used the {\it ordered
complex of curves} and obtained a presentation for $\M(F_g^n)$ in
terms of the mapping class groups of the complementary surfaces of
collections of curves.

If $F_g^n$ is a non-orientable surface of genus $g$ with $n$
boundary components (i.e. $F_g^n$ is homeomorphic to the connected
sum of $g$ projective planes, from which $n$ open discs have been
removed), then presentations for $\M(F_g^n)$ are known only for
$g\le 3$ and small $n$. The complex of curves of $F_g^n$ has been
studied by various authors. Ivanov \cite{I} determined its homotopy
type used it to compute the virtual cohomological dimension of the
mapping class group $\M(F_g^n)$.

In this paper we study the action of the mapping class group $\M(F)$
on the complex of curves of a non-orientable surface $F=F_g^n$. Our
main result says that $\M(F)$ can be presented in terms of the
isotropy subgroups of the collections of curves, provided that $F$
is not sporadic, i.e. the complex of curves of $F$ is simply
connected. On the other hand we show that a presentation for the
isotropy subgroup of a collection of curves $A$ can be obtained from
a presentation for the mapping class group of the surface obtained
by cutting $F$ along $A$. Thus our result recursively produces a
presentation for $\M(F)$, provided that we know presentations for
the mapping class groups of all sporadic subsurfaces. In this paper
we compute an explicit finite presentation for the mapping class
group of each sporadic surface.

The paper is organized as follows, In the next two sections we
present basic definitions and preliminary results about simple
closed curves. In Section \ref{stab} we determine the structure of
the stabilizer of a simplex of the complex of curves, and in Section
\ref{orbits} we determine $\M(F)$-orbits of simplices. In Section
\ref{pre} we use the ordered complex of curves to obtain, by a
result of Brown \cite{Br}, a presentation for the mapping class
group. Then we show how this presentation can be simplified.
Finally, in Section \ref{sporadic} we compute presentations for
mapping class groups of sporadic surfaces.

%% file: preli.tex
\section{\label{preli}Basic definitions.}
Let $F$ denote a  smooth, compact, connected surface, orientable
or not, possibly with boundary. Define $\Dif(F)$ to be the group
of all (orientation preserving if $F$ is orientable) diffeomorphisms
$h\colon F\to F$ such that $h$ is the identity on the boundary of
$F$. The {\it mapping class group} $\M(F)$ is the group of isotopy
classes in $\Dif(F)$. By abuse of notation we will use the same
symbol to denote a diffeomorphism and its isotopy class.
If $g$ and $h$ are two diffeomorphisms, then the composition $gh$
means that $h$ is applied first.

\medskip

By a {\it simple closed curve} in $F$ we mean an embedding
$a\colon S^1\to F$. Note that $a$ has an orientation; the curve
with opposite orientation but same image will be denoted by
$a^{-1}$. By abuse of notation, we also use $a$ for the image of
$a$. If $a_1$ and $a_2$ are isotopic, we write $a_1\simeq a_2$.

We say that $a\colon S^1\to F$ is {\it non-separating} if $F\backslash a$ is connected
and {\it separating} otherwise.
According to whether a regular neighborhood of $a$ is an annulus
or a M\"obius strip, we call $a$ respectively {\it two-} or {\it one-sided}.
If $a$ is one-sided, then we denote by $a^2$ its double,
i.e. the curve $a^2(z)=a(z^2)$ for $z\in S^1\subset\mathbb{C}$.
Note that although $a^2$ is not simple, it is freely homotopic to a
two-sided simple closed curve.

We say that $a$ is {\it essential} if it neither
bounds a disk nor is isotopic to a boundary curve. We say that $a$
is {\it generic} if it is essential and does not bound a M\"obius
strip. Note that every one-sided curve is generic.

Define a {\it generic r-family of disjoint curves} to be a $r$-tuple
$(a_1,\dots,a_r)$ of generic simple closed curves satisfying:
\begin{itemize}
\item $a_i\cap a_j=\emptyset$, for $i\ne j$;
\item $a_i$ is neither isotopic to $a_j$ nor to $a_j^{-1}$, for $i\ne j$.
\end{itemize}

We say that two generic $r$-families of disjoint curves
$(a_1,\dots,a_r)$ and $(b_1,\dots,b_r)$ are {\it equivalent} if
there exists a permutation $\sigma\in\Sigma_r$ such that
$a_i\simeq b_{\sigma(i)}^{\pm 1}$for each $1\le i\le r$. We write
$[a_1,\dots,a_r]$ for the equivalence class of a generic
$r$-family of disjoint curves.

\medskip

The {\it complex of curves} of $F$ is the simplicial complex
$\C(F)$ whose $r$-simplices are the equivalence classes of generic
$(r+1)$-families of disjoint curves in $F$. Vertices of $\C(F)$
are the isotopy classes of unoriented generic curves. The mapping
class group $\M(F)$ acts simplicially on $\C(F)$ by
$h[a_1,\dots,a_r]=[h\circ a_1,\dots,h\circ a_r]$.

\section{A few results about simple closed curves}

A {\it bigon} cobounded by two transversal simple closed curves
$a$ and $b$ is a region in $F$, whose interior is an open disc and
whose boundary is the union of an arc of $a$ and an arc of $b$.
Moreover, we assume that except for the endpoints, these arcs are
disjoint from $a\cap b$, and that the endpoints do not coincide.
If the endpoints coincide (i.e. the arcs are closed curves), then
we say that $a$ and $b$ cobound a {\it degenerate bigon}.

\begin{lemma}[Epstein \cite{E}]\label{epstein}
Let $a,b$ be two two-sided essential curves in $F$, and suppose $a$ is isotopic to $b$.

i) If $a\cap b=\emptyset$, then there exists an annulus in $F$ whose boundary components
are $a$ and $b$.

ii) If $a\cap b\ne\emptyset$, and they intersect transversely, then $a$ and $b$ cobound
a bigon.
\foorp\end{lemma}

\begin{lemma}\label{epstein2}
Let $a,b$ be two one-sided simple closed curves and suppose $a$ is isotopic to $b$.
Then $a\cap b\ne\emptyset$. If they intersect transversely, then:

i) if $|a\cap b|=1$, then $a$ and $b$ cobound a degenerate bigon,

ii) if $|a\cap b|>1$, then $a$ and $b$ cobound a bigon.
\end{lemma}
\proof We choose a regular neighborhood $N_a$ of $a$,
diffeomorphic to the M\"obius strip, and denote by $a'$ its
boundary curve which is homotopic to $a^2$. Similarly we define
$N_b$ and $b'$ homotopic to $b^2$. Now $a'$ and $b'$ are simple
closed curves and $a'\simeq b'$, since $a\simeq b$.

If $F$ is the projective plane or the M\"obius strip, then the proof is trivial.
In the other case $a'$ and $b'$ are essential and we can apply Lemma \ref{epstein}.

Assume $a\cap b=\emptyset$. Then we can choose $N_a$ and $N_b$ disjoint. By Lemma
\ref{epstein}, $a'$ and $b'$ cobound an annulus $A$.
But then $F=A\cup N_a\cup N_b$ is diffeomorphic to the
Klein bottle and $a$ and $b$ are clearly not isotopic. Thus we have proved that $a$ and
$b$ intersect.

Assume that $a$ and $b$ intersect transversely. Then we can choose
$N_a$ and $N_b$ in such a way that $a'$ and $b'$ also intersect
transversely and $|a'\cap b'|=4|a\cap b|$. By Lemma \ref{epstein}
$a'$ and $b'$ cobound a bigon $D$. If $|a\cap b|=1$ then
$M=N_a\cup N_b\cup D$ is a M\"obius strip which contains $a$ and
$b$. In this case $a$ and $b$ cobound a degenerate bigon in $M$.
Assume that $|a\cap b|\ge 2$. Then there exist an arc $c$ of $a$,
an arc $d$ of $b$ and closed subsets $N_c\subset N_a$ and
$N_d\subset N_b$ such that: $|c\cap d|=2$ and the interior of
$N_c\cup N_d\cup D$ is homeomorphic to an open disc. Now $c$ and
$d$ cobound a bigon in $N_c\cup N_d\cup D$. \foorp

\medskip

The next two propositions are proved in \cite{PR} (Propositions
3.5 and 3.10) for orientable surfaces. Their proofs are based on
Lemma \ref{epstein} and can by applied also in the non-orientable
case if the involved curves are two-sided. Therefore, in the
proofs we restrain ourselves to the case of one-sided curves,
where we use Lemma \ref{epstein2} instead of Lemma \ref{epstein}.

\medskip

By a {\it subsurface} $N$ of $F$ we mean a closed subset which is also a surface.
We say furthermore that $N$ is {\it essential} if no boundary curve of $N$ bounds a
disk in $F$.

\begin{prop}\label{subsurface}
Let $N$ be an essential subsurface of $F$, and let $a,b\colon S^1\to N$ be two essential
simple closed curves. (In particular $a$ is not isotopic to a boundary curve of $N$.) Then
$a$ is isotopic to $b$ in $F$ if and only if $a$ is isotopic to $b$ in $N$.
\end{prop}

\proof The nontrivial thing to show is that if $a$ and $b$ are
isotopic in $F$, then they are also isotopic in $N$. We assume
that $a$ and $b$ are one-sided. By Lemma \ref{epstein2} they
intersect. We may assume that they intersect transversally and
argue by induction on $|a\cap b|$.

If $|a\cap b|=1$, then by Lemma \ref{epstein2}, $a$ and $b$ cobound a degenerate bigon $D$ in $F$.
Since $N$ is essential, $D\cap\bdr{N}=\emptyset$ and hence $D\subset N$.
Now we can use $D$ to define an isotopy in $N$ from $a$ to $b^{\pm 1}$. If $a\simeq b^{-1}$
in $N$, then $b\simeq b^{-1}$ in $F$, which can only happen if $F$ is the projective plane
(c.f. \cite{E}, Theorem 1.7). But the projective plane does not contain any essential subsurface.
Thus $a\simeq b$ in $N$.

If $|a\cap b|>1$, then by Lemma \ref{epstein2} $a$ and $b$ cobound a bigon $D\subset F$. As above,
$D\subset N$ and we can use $D$ to define an isotopy in $N$ from $b$ to a curve $b'$ with
$|a\cap b'|=|a\cap b|-2$. By the inductive hypothesis, $b'$ is isotopic to $a$ in $N$, hence so is $b$.
\foorp

\begin{prop}\label{families}
Let $(a_1,\dots,a_r)$, $(b_1,\dots,b_r)$ be two generic $r$-families of disjoint curves
such that $a_i\simeq b_i$ for all $1\le i\le r$. Then there exists an isotopy
$h_t\colon F\to F$, $t\in[0,1]$, such that $h_0=\textrm{identity}$ and
$h_1\circ a_i=b_i$ for all $1\le i\le r$.
\end{prop}

\proof We use induction on $r$. The proposition is obvious for
$r=1$ and we assume that it is true for $(r-1)$-families.
Replacing each $a_i$ by $h_1\circ a_i$, we may assume that
$a_i=b_i$ for $1\le i\le r-1$. Then $a_r$ and $b_r$ are disjoint
from $a_i=b_i$ for $i<r$ and $a_r\simeq b_r$. Now it suffices to
show that there is an isotopy of $F$ which takes $a_r$ to $b_r$
and does not move the curves $a_i=b_i$ for $i<r$. We assume that
$a_r$ and $b_r$ are one-sided and intersect transversally. We
argue by induction on $|a_r\cap b_r|$.

If $|a_r\cap b_r|=1$, then by Lemma \ref{epstein2}, $a_r$ and
$b_r$ cobound a degenerate bigon $D$ in $F$. Since the curves
$a_i=b_i$ for $i<r$ are generic, they are all disjoint from $D$.
Now it is easy to construct an isotopy of $F$, which takes $a_r$
to $b_r$ across $D$ and is equal to the identity outside a
neighborhood of $D$, so the other curves do not move.

If $|a_r\cap b_r|>1$, then by Lemma \ref{epstein2}, $a_r$ and
$b_r$ cobound a bigon $D$ in $F$. As above, the curves $a_i=b_i$
for $i<r$ are disjoint from $D$, and there is an isotopy of $F$,
fixed outside a neighborhood of $D$, which takes $a_r$ across $D$
and reduces the number $|a_r\cap b_r|$ without moving the other
curves. By the inductive hypothesis there is a final isotopy
taking $a_r$ to $b_r$. \foorp

\medskip

Given a two-sided simple closed curve $a$ we can define a Dehn
twist $t_a$ about $a$. Since we are dealing with non-orientable
surfaces, it is impossible to distinguish between right and left
twists. The direction of a twist $t_a$ has to be specified for
each curve $a$. Equivalently we may choose an orientation of a
tubular neighborhood of $a$. Then $t_a$ denotes the right Dehn
twist with respect to the chosen orientation. Unless we specify
which of the two twists we mean, $t_a$ denotes (the isotopy class
of) any of the two possible twists.

The next Proposition is proved in \cite{PR} for orientable surfaces and in
\cite{S} for
non-orientable surfaces.

\begin{prop}\label{twists}
Suppose that $F$ is not homeomorphic to the Klein bottle.
Consider $r$ two-sided simple closed curves $a_1,\dots,a_r$ satisfying:

i) $a_i$ is either generic or isotopic to a boundary curve;

ii)  $a_i\cap a_j=\emptyset$, for $i\ne j$;

iii) $a_i$ is neither isotopic to $a_j$ nor to $a_j^{-1}$, for $i\ne j$.

Then the subgroup of $\M(F)$ generated by Dehn twists $t_{a_1},\dots, t_{a_r}$ is
a free abelian group of rank $r$.
\foorp\end{prop}

Note that if $F$ is homeomorphic to the Klein bottle, then up to isotopy there is only
one generic two-sided curve $a$, and $t_a$ has order 2.

%% file: stab.tex
\section{\label{stab}The structure of the stabilizer}

In this section we follow the outline of Paragraph 6 of \cite{P} to
expresses the stabilizer of a simplex of $\C(F)$ by means of the
mapping class group of the complementary surface. Our Proposition
\ref{sequence} is a generalization to the case of a non-orientable
surface of Proposition 6.3 of \cite{P}.

\medskip

Let $A=(a_1,\dots,a_r)$ be a generic $r$-family of disjoint curves. Denote by $F_A$ the
compact surface obtained by cutting $F$ along $A$, i.e. the natural compactification of
$F\backslash(\bigcup_{i=1}^{r}a_i)$. Note that $F_A$ is in general not connected.
Denote by $N_1,\dots,N_k$ the connected components of $F_A$. Then we write
$$\M(F_A)=\M(N_1)\times\dots\times\M(N_k).$$
Denote by $\rho_A\colon F_A\to F$ the continuous map induced by
the inclusion of $F\backslash(\bigcup_{i=1}^{r}a_i)$ in $F$. The
map $\rho_A$ induces a homomorphism
$\rho_\ast\colon\M(F_A)\to\M(F)$.

\medskip

A {\it pair of pants} is a compact surface homeomorphic to a sphere with $3$ holes.
We say that the family $A$ determines a {\it pants decomposition} if each component
of $F_A$ is a pair of pants. Such a family exists if and only if the Euler
characteristic of $F$ is negative. In such case, a generic family
$A$ determines a pants decomposition if and only if $A$ represents
a maximal simplex in $\C(F)$. Given a generic family $A=(a_1,\dots,a_r)$ we can always
complete it to a pants decomposition, i.e. there exist generic curves
$(a_{r+1},\dots,a_s)$ such that $(a_1,\dots,a_s)$ determines a pants decomposition.
Recall that if $N$ is a pair of pants then $\M(N)$ is the free abelian group of rank $3$
generated by Dehn twists along the boundary curves.

\begin{lemma}\label{kernel}
Assume that $F$ has negative Euler characteristic.
Let $A=(a_1,\dots,a_r)$ be a generic
family of disjoint curves in $F$ such that $a_1,\dots,a_p$ are two-sided and
$a_{p+1},\dots,a_r$ are one-sided.
For each $i\in\{1,\dots,p\}$ let $a'_i$ and $a''_i$ denote the
boundary curves of $F_A$ such that $\rho_A\circ a'_i=\rho_A\circ
a''_i=a_i$, and choose $t_{a'_i}$ and $t_{a''_i}$  so that
$\rho_\ast(t_{a'_i})=\rho_\ast(t_{a''_i})$.
For each $j\in\{p+1,\dots,r\}$ let $a'_j$ denote the boundary curve
of $F_A$ such that $\rho_A\circ a'_j=a^2_j$.
Then $\ker\rho_\ast$ is generated by
$\{t_{a'_1}t_{a''_1}^{-1},\dots,t_{a'_p}t_{a''_p}^{-1},t_{a'_{p+1}},\dots,t_{a'_r}\}$ and is
a free abelian group of rank $r$.
\end{lemma}
\proof
Let $G$ denote the subgroup of $\M(F_A)$ generated by
$$\{t_{a'_1}t_{a''_1}^{-1},\dots,t_{a'_p}t_{a''_p}^{-1},t_{a'_{p+1}},\dots,t_{a'_r}\}.$$
Clearly $G\subseteq\ker\rho_\ast$ and it follows from Proposition \ref{twists}
that $G$ is a free abelian group of rank $r$. It remains to show that
$\ker\rho_\ast\subseteq G$.

Let $c_1,\dots,c_n$ denote the boundary curves of $F$ and
$c'_1,\dots,c'_n$ the corresponding boundary curves of $F_A$ (i.e.
$\rho_A\circ c'_i=c_i$). Complete $A$ to a pants decomposition
$A'=(a_1,\dots,a_r,a_{r+1},\dots,a_q,\dots,a_s)$, where
$a_{r+1},\dots,a_q$ are two-sided and $a_{q+1},\dots,a_s$ one-sided.
Let $a'_{r+1},\dots,a'_s$ denote the generic curves in $F_A$ such that
$\rho_A\circ a'_j=a_j$ for $r+1\le j\le s$.

Let $h$ be an element of $\ker\rho_\ast$ and $j\in\{r+1,\dots,s\}$. We have
$\rho_A\circ h\circ a'_j\simeq\rho_A\circ a'_j$ and it follows by Proposition
\ref{subsurface} that $h\circ a'_j\simeq a'_j$. Hence, by Proposition
\ref{families} we may assume that $h\circ a'_j=a'_j$. Now
$h$ induces a diffeomorphism of $F_{A'}$, and hence by the
structure of the mapping class group of the pair of pants we can write:
$$h=t^{u_1}_{a'_1}t^{v_1}_{a''_1}\dots t^{u_p}_{a'_p}t^{v_p}_{a''_p}
t^{u_{p+1}}_{a'_{p+1}}\dots t^{u_q}_{a'_q} t^{w_1}_{c'_1}\dots
t^{w_n}_{c'_n},$$ where $u_1,\dots,w_n\in\Z$. The equality
$$1=\rho_\ast(h)=t^{u_1+v_1}_{a_1}\dots t^{u_p+v_p}_{a_p}
t^{u_{r+1}}_{a_{r+1}}\dots t^{u_q}_{a_q}
t^{w_1}_{c_1}\dots t^{w_n}_{c_n}$$
implies by Proposition \ref{twists}:
$$u_1+v_1=\dots=u_p+v_p=u_{r+1}=\dots=u_q=w_1=\dots=w_n=0,$$
and we have
$h=(t_{a'_1}t^{-1}_{a''_1})^{u_1}\dots(t_{a'_p}t^{-1}_{a''_p})^{u_p}
t^{u_{p+1}}_{a'_{p+1}}\dots t^{u_r}_{a'_r}$.\foorp

\medskip

Denote by $[A]$ the simplex in $\C(F)$ represented by the family
$A=(a_1,\dots,a_r)$, and by $\textrm{Stab}([A])$ the stabilizer of
$[A]$ in $\M(F)$.

Define the {\it cubic group} $\textrm{Cub}_r$ to be the group of
linear transformations $\phi\in GL(\R^r)$ such that $\phi(e_i)=\pm
e_j$ for all $1\le i\le r$, where $\{e_1,\dots,e_r\}$ denotes the
canonical basis of $\R^r$. There is a natural homomorphism
$\Phi_A:\textrm{Stab}([A])\to \textrm{Cub}_r$ defined as follows:
$$\Phi_A(h)(e_i)=\left\{\begin{array}{ll}
e_j & \textrm{if $h\circ a_i\simeq a_j$},\\
-e_j & \textrm{if $h\circ a_i\simeq a_j^{-1}$}.
\end{array}\right.$$
Denote by $\textrm{Stab}^+([A])$ the kernel of $\Phi_A$. By
Proposition \ref{families}, each element of $\textrm{Stab}^+([A])$
is represented by a diffeomorphism $h\in\Dif(F)$, such that
$h\circ a_i=a_i$ for all $1\le i\le r$. Consider the subgroup $H$
of $\textrm{Stab}^+([A])$ consisting of the isotopy classes of
diffeomorphisms preserving each curve of $A$ with its orientation
and preserving orientation of a tubular neighborhood of each two-sided
curve of $A$. If $A$ contains $p$ two-sided curves, then there is
an obvious homomorphism $\textrm{Stab}^+([A])\to(\Z_2)^p$ with
kernel $H$. Finally observe that $H$ is equal to
$\textrm{Im$\rho_\ast$}$.

Now we can summarize the considerations of this Section in the
following proposition.

\begin{prop}\label{sequence}
Assume that $F$ is a surface of negative Euler characteristic.
Let $A$ be a generic
$r$-family of disjoint curves containing $p$ two-sided curves ($0\le p\le r$). Then we have the following exact sequences:
\begin{eqnarray*}
1\to\Z^r\to\M(F_A)\stackrel{\rho_\ast}{\to}{\rm{Stab}}^+([A])\to(\Z_2)^p,\\
1\to{\rm{Stab}}^+([A])\to{\rm{Stab}}([A])\stackrel{\Phi_A}{\to}{\rm{Cub}}_r.
\end{eqnarray*}
\end{prop}

\begin{rem}
The homomorphisms $\textrm{Stab}^+([A])\to(\Z_2)^p$ and $\Phi_A$ are in general not surjective. By an easy analysis case by case
it is  possible to describe their images exactly.
\end{rem}

%% file: orbits.tex
\section{\label{orbits}The orbits}

For the rest of this paper we assume that $F=F_g^n$ is a {\it
non-orientable} surface of genus $g$ with $n$ boundary components
($n\ge 0$). Recall that this means that $F$ is diffeomorphic to
the connected sum of $g$ projective planes, from which $n$
disjoint open discs have been removed. We also assume that $F$ has
negative Euler characteristic, i.e. $g+n>2$. In this section we
determine the $\M(F)$-orbits of simplices of the complex of curves
$\C(F)$. We say that two simplices $[A]$ and $[B]$ of $\C(F)$ are
$\M(F)$-{\it equivalent} if they are in the same $\M(F)$-orbit. If
$A=(a_1,\dots,a_r)$, $B=(b_1,\dots,b_r)$, then $[A]$ and $[B]$ are
$\M(F)$-equivalent if and only if there exist $h\in\Dif(F)$ and
permutation $\sigma\in\Sigma_r$, such that $h\circ a_i\simeq
b^{\pm 1}_{\sigma(i)}$. By Proposition \ref{families} that is
equivalent to existence of $h\in\Dif(F)$, such that $h\circ a_i=
b^{\pm 1}_{\sigma(i)}$.

\medskip

Let $A=(a_1,\dots,a_r)$ be a generic family of disjoint curves.
Let us fix boundary curves $c_1,\dots, c_n$ of $F$.
By abuse of notation we also denote by $c_i$ the boundary curve
$c_i\colon S^1\to\bdr{N}$ such that $\rho_A\circ c_i=c_i$, where $N$ is a connected component of
$F_A$. We say that $c_i$ is an {\it exterior boundary curve} of $N$.

Let $a_i\colon S^1\to F$ be a two-sided curve in the family $A$. There exist
two connected components $N'$ and $N''$ of $F_A$, and two distinct
curves $a'_i\colon S^1\to\bdr{N'}$ and $a''_i\colon
S^1\to\bdr{N''}$ such that $\rho_A\circ a'_i=\rho_A\circ
a''_i=a_i$. We say that $a_i$ is a {\it separating limit curve} of
$N'$ (and $N''$) if $N'\ne N''$, and $a_i$ is a {\it
non-separating two-sided limit curve} of $N'$ if $N'=N''$.

Let $a_i\colon S^1\to F$ be a one-sided curve in $A$. There exists a component
$N$ of $F_A$ and a curve $a'_i\colon S^1\to\bdr{N}$ such that
$\rho_A\circ a'_i=a^2_i$. We say that $a_i$ is
a {\it one-sided limit curve} of $N$.

\begin{lemma}\label{boundary slide}
Suppose that $N$ is a non-orientable connected surface and
$c\colon S^1\to\bdr{N}$ is a boundary curve in $N$. There exists a diffeomorphism
$h\colon N\to N$ such that $h$ is the identity on $\bdr{N}\backslash c$ and
$h\circ c=c^{-1}$.
\end{lemma}

\proof Let $N'$ be the surface obtained by gluing a disc $D$ to
$N$ along $c$. Let $p$ be the center of $D$, and $\alpha\colon
(S^1,1)\to(N'\backslash\bdr{N'},p)$ any one-sided simple loop
based at $p$. There exists an isotopy $h_t\colon N'\to N'$, $0\le
t\le 1$, such that: $h_0=\textrm{identity}$,
$h_t(p)=\alpha(e^{2\pi t})$, $h_t$ is the identity on $\bdr{N'}$
for all $t$, and $h_1\circ c=c^{-1}$. We define $h\colon N\to
N$ to be the restriction of $h_1$ to $N$. Such diffeomorphism is
called {\it the boundary slide} (c.f. \cite{K}). \foorp

\begin{prop}\label{simplices}
Let $A=(a_1,\dots,a_r)$ and $B=(b_1,\dots,b_r)$ be two generic $r$-families of
disjoint curves. The simplices $[A]$ and $[B]$ are
$\M(F)$-equivalent if and only if there exists a permutation $\sigma\in\Sigma_r$, such that
for all subfamilies $A'\subseteq A$ and $B'\subseteq B$, such that $a_i\in A' \iff b_{\sigma(i)}\in B'$,
there exists a one to one
correspondence between the connected components of $F_{A'}$ and those
of $F_{B'}$, such that  for
every pair $(N,N')$ where $N$ is any component of $F_{A'}$ and $N'$
is the corresponding component of $F_{B'}$, we have:
\begin{itemize}
\item
$N$ and $N'$ are either both orientable or both non-orientable, of the same genus;
\item
if $c_i$ is an exterior boundary curve of $N$, then it is also an exterior boundary curve of $N'$;
\item
if $N$ is orientable and $c_i$ and $c_j$ induce the same orientation of $N$, then
they also induce the same orientation of $N'$;
\item
if $a_i$ is a separating limit curve of $N$, then $b_{\sigma(i)}$ is a separating limit curve
of $N'$;
\item
if $a_i$ is a non-separating two-sided limit curve of $N$, then $b_{\sigma(i)}$ is a non-separating two-sided limit curve
of $N'$;
\item
if $a_i$ is a one-sided limit curve of $N$, then $b_{\sigma(i)}$ is a one-sided limit curve
of $N'$.
\end{itemize}
\end{prop}

\proof Suppose $[A]$ and $[B]$ are $\M(F)$-equivalent. Then there
exist $h\in\Dif(F)$ and $\sigma\in\Sigma_r$, such that $h\circ
a_i=b_{\sigma(i)}^{\pm 1}$ for $1\le i\le r$. For each subfamily
$A'\subseteq A$, $h$ induces a diffeomorphism $h'\colon F_{A'}\to
F_{B'}$, such that $h\circ\rho_{A'}=\rho_{B'}\circ h'$. We define
a correspondence between the connected components of $F_{A'}$ and
those of $F_{B'}$ as follows. If $N$ is any component of $F_{A'}$
then $N'=h'(N)$ is the corresponding component of $F_{B'}$. Note
that we have $h'\circ c_i=c_i$ and hence $c_i$ is an exterior
boundary curve of $N$ if and only if it is an exterior boundary
curve of $N'$. Furthermore, if $N$ is orientable and $c_i$, $c_j$
induce the same orientation of $N$, then they also induce the same
orientation of $N'$. Suppose that $a_i\in A'$ is a two-sided limit
curve of $N$. Then $a_i=\rho_{A'}\circ a'_i$ for $a'_i\colon
S^1\to\bdr{N}$ and $b_{\sigma(i)}=h\circ a^{\pm
1}_i=h\circ\rho_{A'}\circ (a'_i)^{\pm 1}=\rho_{B'}\circ h'\circ
(a'_i)^{\pm 1}$. Hence $b_{\sigma(i)}$ is a two-sided limit curve
of $N'$. Clearly if $a_i$ is separating then so is
$b_{\sigma(i)}$. Similarly, if $a_i$ is a one-sided limit curve of
$N$ and $a_i^2=\rho_{A'}\circ a'_i$, then
$b^2_{\sigma(i)}=\rho_{B'}\circ h'\circ (a'_i)^{\pm 1}$ and
$b_{\sigma(i)}$ is a one-sided limit curve of $N'$.

Assume now that there exists a permutation $\sigma\in\Sigma_r$, such that for each
subfamily $A'\subseteq A$ the conditions of the proposition are satisfied. Let us
assume, for simplicity, that $\sigma$ is the trivial permutation
$\sigma(i)=i$ for $1\le i\le r$. Denote by
$N_1,\dots,N_k$ the connected components of $F_A$, and by $N'_1,\dots,N'_k$
the corresponding components of $F_B$. By the classification of compact
surfaces there exist diffeomorphisms $h_i\colon N_i\to N'_i$,
$1\le i\le k$, such that for each exterior boundary curve
$c_l\colon S^1\to\bdr N_i$ we have $h_i\circ c_l=c_l^{\pm 1}$,
and if $a_j$ is a limit curve of $N_i$, then $\rho_B\circ
h_i\circ a'_j=b^{\pm 1}_j$ if $a_j=\rho_A\circ a'_j$, and
$\rho_B\circ h_i\circ a'_j=(b^2_j)^{\pm 1}$ if $a_j^2=\rho_A\circ
a'_j$. We will show that we can choose $h_i$ so that
for each boundary curve
\begin{equation}\label{w1} h_i\circ c_l=c_l,\end{equation}
and for each two-sided limit curve $a_j$ of $N_i$ and $N_m$, if
$a_j=\rho_A\circ a'_j=\rho_A\circ a''_j$ then
\begin{equation}\label{w2} \rho_B\circ h_i\circ a'_j=b_j\iff\rho_B\circ h_m\circ a''_j=b_j.\end{equation}
Notice that if $h_i$ satisfy (\ref{w1}) and (\ref{w2}),
then they induce $h\in\Dif(F)$ such that $h\circ a_j= b^{\pm 1}_j$ for $1\le j\le r$, which proves
Proposition \ref{simplices}.

If all $N_i$ are non-orientable, then by Lemma \ref{boundary
slide} we can compose $h_i$ with suitable boundary slides, so that
(\ref{w1}) and (\ref{w2}) are satisfied. Suppose that $N_1,\dots,
N_s$, $1\le s\le k$  are all orientable components of $F_A$. We
define $A'\subseteq A$ to be any maximal subfamily consisting of
separating limit curves of $N_1,\dots, N_s$ such that: the surface
$M$ obtained by gluing $\coprod_{i=1}^s N_i$ along $A'$ is
orientable; each $a_i\in A'$ separates $M$, i.e. $M\backslash a_i$
has more connected components than $M$. Notice that $A'$ may be
empty. The surface $M$ is in general disconnected and it is the
sum of all orientable components of $F_{A\backslash A'}$. Let $M'$
denote the surface obtained by gluing $\coprod_{i=1}^s N'_i$ along
$B'$, where $b_i\in B'\Leftrightarrow a_i\in A'$. Notice that $M'$
is the sum of all orientable components of $F_{B\backslash B'}$.
We claim that we can choose $h_i$ for $i\le s$, so that (\ref{w2})
holds for each $a_j\in A'$. First notice, that after re-numbering
the orientable components of $F_A$ if necessary, we may assume
that for each $m\le s$ there is at most one $a_j\in A'$ such that
$a_j$ is a separating limit curve of $N_m$ and $N_i$ for $i<m$.
Now we define $h_i$ inductively. We choose any $h_1$. Suppose that
we have chosen $h_i$ for all $i<m\le s$. If there is $a_j\in A'$
such that $a_j$ is a separating limit curve of $N_m$ and $N_i$ for
$i<m$, then we choose $h_m$ so that (\ref{w2}) is satisfied. If
there is no such curve, then we choose any $h_m$.
Such chosen $h_i$ induce
$\widetilde{h}\colon M\to M'$, so that $\widetilde{h}\circ
c_l=c^{\pm 1}_l$ for each exterior boundary curve of $M$. Let
$c_i$, $c_j$ be two exterior boundary curves of one component of
$M$. Since $A\backslash A'$ and $B\backslash B'$ satisfy the
conditions of the proposition, $c_i$ and $c_j$ induce the same
orientation of the component of $M$ if and only if they induce the
same orientation of the corresponding component of $M'$, hence
$\widetilde{h}\circ c_i=c_i\Leftrightarrow \widetilde{h}\circ
c_j=c_j$. Now it is clear that composing if necessary some $h_i$
with orientation reversing diffeomorphism, we can assume
$\widetilde{h}\circ c_l=c_l$ for each exterior boundary curve of
$M$. Thus $h_i$ also satisfy (\ref{w1}).

Suppose that $a_j\in A\backslash A'$ is a two-sided limit curve of
$N_i$ and $N_m$, $i\le m\le s$. Since $A'$ is maximal, $a_j$ is a
non-separating limit curve of some component $M_j$ of $M$, i.e.
$a_j=\rho_{A\backslash A'}\circ a'_j=\rho_{A\backslash A'}\circ
a''_j$ for $a'_j,a''_j\colon S^1\to M_j$. Then
$b_j=\rho_{B\backslash B'}\circ b'_j=\rho_{B\backslash B'}\circ
b''_j$ for $b'_j=\widetilde{h}\circ(a'_j)^{\pm 1}$,
$b''_j=\widetilde{h}\circ(a''_j)^{\pm 1}$. Note that the surface
obtained from $M_j$ by gluing along $a_j$ is orientable if and
only if $a'_j$ and $a''_j$ induce opposite orientations of $M_j$.
Since $A\backslash(A'\cup\{a_j\})$ and
$B\backslash(B'\cup\{b_j\})$ satisfy the conditions of the
proposition, the surface obtained from $M$ by gluing along $a_j$
is diffeomorphic to the surface obtained by gluing $M'$ along
$b_j$. In particular, one of these surfaces is orientable if and
only if the other one is. Hence $a'_j$ and $a''_j$ induce the same
orientation of $M_j$ if and only if $b'_j$ and $b''_j$ induce the
same orientation of $\widetilde{h}(M_j)$. Thus $\widetilde{h}\circ
a'_j=b'_j\Leftrightarrow \widetilde{h}\circ a''_j=b''_j$ and so
(\ref{w2}) holds for $a_j$.

Once we have chosen $h_i$ for $i\le s$, it is easy to construct, using Lemma \ref{boundary
slide}, diffeomorphisms $h_i$ for $i>s$ satisfying
(\ref{w1}) and (\ref{w2}) for all curves.
\foorp

\begin{cor}\label{finite}
There are only finitely many $\M(F)$-orbits in $\C(F)$.
\end{cor}

\proof Let $N$ be a disjoint union of $g+n-2$ pairs of pants. Choose
boundary curves of $N$
\begin{equation}\label{curves}
c_1,\dots,c_n, a'_1,\dots,a'_s, a''_1,\dots,a''_r,
\end{equation}
where $r\le s$, $n+r+s=3(g+n-2)$. Consider the surface $M=N/\sim$,
where $\sim$ identifies pairs of boundary points as follows:
$a'_i(z)=a''_i(z)$ for $i\le r$, $a'_i(z)=a'_i(z^2)$ for $i>r$.
Let $\rho\colon N\to M$ denote the canonical projection. Generic
family of disjoint curves $(a_1,\dots,a_s)$, where $a_i=\rho\circ
a'_i$ for $i\le r$, $a_i^2=\rho\circ a'_i$ for $i>r$, determines a
pants decomposition of $M$. Notice that for some choices of curves
(\ref{curves}) we have $M=F_g^n$, i.e. $M$ is a connected,
non-orientable surface of genus $g$. Furthermore, every pants
decomposition of $F_g^n$ can be obtained in this way, and thus, by
Proposition \ref{simplices}, there is at most as many
$\M(F)$-orbits of pants decompositions, as the number of different
(i.e. not isotopic) choices of curves (\ref{curves}). Since that
number is finite and every generic family of disjoint curves can
be completed to a pants decomposition, there are only finitely
many $\M(F)$-orbits in $\C(F)$. \foorp

\begin{figure}%[!htbp]
\begin{center}
\input{fig1_ai}
\caption{\label{ai} Non-separating curves.} \end{center}
\end{figure}
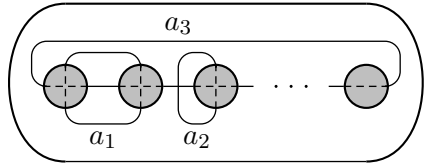

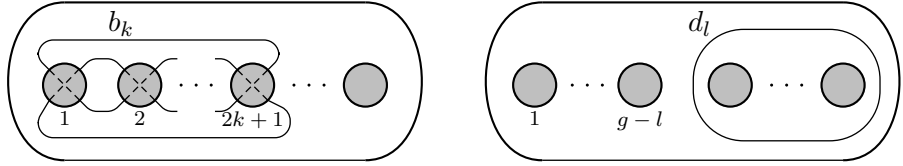
\begin{figure}%[!htbp]
\begin{center}
\begin{tabular}{cc}
\input{fig2_bk} & \input{fig3_dl}
\end{tabular}
\caption{\label{bkdl} Separating curves.} \end{center}
\end{figure}

\medskip

Let us list all $M(F)$-orbits of the vertices of $\C(F)$. We call
a vertex $[a]$ one- or two-sided,  and separating or
non-separating if $a$ has the appropriate property.

Suppose that $F$ is closed and has genus $g\ge 3$.
Consider the three non-separating curves $a_1$, $a_3$, $a_3$
in Figure \ref{ai}. In this figure, and also in other figures in this paper, the shaded discs represent
crosscaps; this means that their interiors should be removed and then the antipodal points in each
boundary component should be identified. We have:
\begin{itemize}
\item $a_1$ is two-sided, $F_{a_1}$ is non-orientable;
\item $a_2$ is one-sided, $F_{a_2}$ is non-orientable;
\item $F_{a_3}$ is orientable, $a_3$ is one-sided if $g$ is odd, and
two-sided if $g$ is even.
\end{itemize}
For each integer $k$, such that $1\le k\le\frac{g}{2}-1$ and for each $l$ such that
$2\le l\le\frac{g}{2}$ we define separating generic curves $b_k$ and $d_l$ represented
in Figure \ref{bkdl}. We have:
\begin{itemize}
\item one component of $F_{b_k}$ is orientable and has genus $k$, the other
component is non-orientable and has genus $g-2k$;
\item both components of $F_{d_l}$ are non-orientable and have genera $l$ and $g-l$.
\end{itemize}
By Proposition \ref{simplices}, every vertex of $\C(F)$ is
$\M(F)$-equivalent to one of the vertices $[a_1]$, $[a_2]$,
$[a_3]$, $[b_k]$, $[d_l]$. Thus we have $3$ orbits of
non-separating vertices and $2([\frac{g}{2}]-1)$ orbits of
separating vertices, where $[\frac{g}{2}]$ denotes the integer
part of $\frac{g}{2}$.

Now suppose that $F$ has boundary, that is $n\ge 1$, and
$g$ is arbitrary such that $\chi(F)=2-g-n<0$.
For each pair $\{I,I'\}$ of sets such that $I\cup I'=\{1,\dots,n\}$, $I\cap I'=\emptyset$
there is one $\M(F)$-orbit consisting of all non-separating vertices $[a]$ such that
\begin{itemize}
\item $F_a$ is orientable, and $c_i$, $c_j$ induce the same
orientation of $F_a$ if and only if $\{i,j\}\subseteq I$ or
$\{i,j\}\subseteq I'$.
\end{itemize}
There are $2^{n-1}$ such orbits. The remaining non-separating
vertices have form $[a]$, where $F_a$ is non-orientable. If $g=1$
then there are no such vertices. If $g=2$ then they are all
one-sided and form one $\M(F)$-orbit. If $g\ge 3$ then they form
$2$ orbits, one contains all one-sided vertices, the other one
contains all two-sided vertices.

The orbits of separating vertices are of two types, like for closed $F$.
For every integer $k$ such that $0\le k\le\frac{g-1}{2}$, and pair $\{I,J\}$ of disjoint subsets of $\{1,\dots,n\}$
such that $g+n-2\ge 2k+\#(I\cup J)\ge 2$ there is one $\M(F)$-orbit consisting of
all separating vertices $[b]$ such that
\begin{itemize}
\item $F_b$ has one orientable component $N_o$ of genus $k$ and
one non-orientable component $N_n$ of genus $g-2k$; \item
$c_i\subset N_o\Leftrightarrow i\in(I\cup J)$; $c_i$, $c_j$ induce
the same orientation of $N_o$ if and only if $\{i,j\}\subseteq I$
or $\{i,j\}\subseteq J$.
\end{itemize}
For every integer $l$ such that $1\le l\le\frac{g}{2}$ and every
$I\subseteq\{1,\dots,n\}$ such that $l+\#I\ge 2$ there is one
$\M(F)$-orbit consisting of all separating vertices $[d]$ such
that
\begin{itemize}
\item $F_d$ has two non-orientable components $N_1$ and $N_2$ of genera $l$ and $g-l$ respectively;
$c_i\subset N_1\Leftrightarrow i\in I$.
\end{itemize}

%% file: fig1_ai.tex
\pspicture*(6,2.5)
\psline[linewidth=.5pt](1,.25)(5,.25)
\psline(1,2.35)(5,2.35)
\psbezier(1,.25)(0,.25)(0,2.35)(1,2.35)
\psbezier(5,.25)(6,.25)(6,2.35)(5,2.35)
\pscircle*[linecolor=lightgray](1,1.25){.3}
\pscircle(1,1.25){.3}
\pscircle*[linecolor=lightgray](2,1.25){.3}
\pscircle(2,1.25){.3}
\pscircle*[linecolor=lightgray](3,1.25){.3}
\pscircle(3,1.25){.3}
\pscircle*[linecolor=lightgray](5,1.25){.3}
\pscircle(5,1.25){.3}
\rput[b](2.5,1.95){\small$a_3$}
\psline[linewidth=.5pt,linearc=.2](.7,1.25)(.55,1.25)(.55,1.85)(5.45,1.85)(5.45,1.25)(5.3,1.25)
\psline[linewidth=.5pt,linestyle=dashed,dash=3pt 2pt](.7,1.25)(1.3,1.25)
\psline[linewidth=.5pt](1.3,1.25)(1.7,1.25)
\psline[linewidth=.5pt,linestyle=dashed,dash=3pt 2pt](1.7,1.25)(2.3,1.25)
\psline[linewidth=.5pt](2.3,1.25)(2.7,1.25)
\psline[linewidth=.5pt,linestyle=dashed,dash=3pt 2pt](2.7,1.25)(3.3,1.25)
\psline[linewidth=.5pt](3.3,1.25)(3.5,1.25)
\pscircle*(3.8,1.25){.02}
\pscircle*(4,1.25){.02}
\pscircle*(4.2,1.25){.02}
\psline[linewidth=.5pt](4.5,1.25)(4.7,1.25)
\psline[linewidth=.5pt,linestyle=dashed,dash=3pt 2pt](4.7,1.25)(5.3,1.25)
\rput[t](2.75,.65){\small$a_2$}
\psline[linewidth=.5pt,linearc=.2](3,1.55)(3,1.7)(2.5,1.7)(2.5,.75)(3,.75)(3,.95)
\psline[linewidth=.5pt,linestyle=dashed,dash=3pt 2pt](3,.95)(3,1.55)
\rput[t](1.5,.65){\small$a_1$}
\psline[linewidth=.5pt,linearc=.2](1,1.55)(1,1.7)(2,1.7)(2,1.55)
\psline[linewidth=.5pt,linearc=.2](1,.95)(1,.75)(2,.75)(2,.95)
\psline[linewidth=.5pt,linestyle=dashed,dash=3pt 2pt](1,.95)(1,1.55)
\psline[linewidth=.5pt,linestyle=dashed,dash=3pt 2pt](2,.95)(2,1.55)
\endpspicture

%% file: fig2_bk.tex
\pspicture*(6,2.5)
\psline[linewidth=.5pt](1,.25)(5,.25)
\psline(1,2.35)(5,2.35)
\psbezier(1,.25)(0,.25)(0,2.35)(1,2.35)
\psbezier(5,.25)(6,.25)(6,2.35)(5,2.35) %
\pscircle*[linecolor=lightgray](1,1.25){.3}
\pscircle(1,1.25){.3}
\rput[t](1,.9){\tiny$1$}
\pscircle*[linecolor=lightgray](2,1.25){.3}
\pscircle(2,1.25){.3}
\rput[t](2,.9){\tiny$2$}
\pscircle*[linecolor=lightgray](3.5,1.25){.3}
\pscircle(3.5,1.25){.3}
\rput[t](3.5,.9){\tiny$2k+1$}
\pscircle*[linecolor=lightgray](5,1.25){.3}
\pscircle(5,1.25){.3}
\psline[linewidth=.5pt,linestyle=dashed,dash=3pt 2pt](.8,1.45)(1.2,1.05)
\psline[linewidth=.5pt,linestyle=dashed,dash=3pt 2pt](1.2,1.47)(.8,1.05)
%\psline[linewidth=.5pt](.8,1.7)(.65,1.85)
%\psline[linewidth=.5pt](1.2,1.7)(1.35,1.85)
%\psline[linewidth=.5pt](.8,1.3)(.65,1.15)
%\psline[linewidth=.5pt](1.2,1.3)(1.35,1.15)
\psline[linewidth=.5pt,linestyle=dashed,dash=3pt 2pt](1.8,1.45)(2.2,1.05)
\psline[linewidth=.5pt,linestyle=dashed,dash=3pt 2pt](2.2,1.45)(1.8,1.05)
\psline[linewidth=.5pt,linestyle=dashed,dash=3pt 2pt](3.3,1.45)(3.7,1.05)
\psline[linewidth=.5pt,linestyle=dashed,dash=3pt 2pt](3.7,1.45)(3.3,1.05)
%\psline[linewidth=.5pt](1.8,1.7)(1.65,1.85)
%\psline[linewidth=.5pt](2.2,1.7)(2.35,1.85)
%\psline[linewidth=.5pt](1.8,1.3)(1.65,1.15)
%\psline[linewidth=.5pt](2.2,1.3)(2.35,1.15)
%\psline[linewidth=.5pt,linestyle=dashed,dash=3pt 2pt](2.8,1.7)(3.2,1.3)
%\psline[linewidth=.5pt,linestyle=dashed,dash=3pt 2pt](3.2,1.7)(2.8,1.3)
%\psline[linewidth=.5pt](2.8,1.7)(2.65,1.85)
%\psline[linewidth=.5pt](3.2,1.7)(3.35,1.85)
%\psline[linewidth=.5pt](2.8,1.3)(2.65,1.15)
%\psline[linewidth=.5pt](3.2,1.3)(3.35,1.15)
\psline[linewidth=.5pt,linearc=.2](1.2,1.45)(1.35,1.6)(1.65,1.6)(1.8,1.45)
\psline[linewidth=.5pt,linearc=.2](1.2,1.05)(1.35,.9)(1.65,.9)(1.8,1.05)
\psline[linewidth=.5pt,linearc=.2](2.2,1.45)(2.35,1.6)(2.5,1.6)
\psline[linewidth=.5pt,linearc=.2](3,1.6)(3.15,1.6)(3.3,1.45)
\pscircle*(2.55,1.25){.02}
\pscircle*(2.75,1.25){.02}
\pscircle*(2.95,1.25){.02}
\psline[linewidth=.5pt,linearc=.2](2.2,1.05)(2.35,.9)(2.5,.9)
\psline[linewidth=.5pt,linearc=.2](3,.9)(3.15,.9)(3.3,1.05)
%\psline[linewidth=.5pt](2.35,1.85)(2.65,1.85)
%\psline[linewidth=.5pt](2.35,1.15)(2.65,1.15)
\rput[b](1.75,1.9){\small$b_k$}
\psline[linewidth=.5pt,linearc=.2](.8,1.45)(.65,1.6)(.65,1.85)(3.85,1.85)(3.85,1.6)(3.7,1.45)
%\psline[linewidth=.5pt,linearc=.2](.8,1.05)(.65,.9)(.65,.65)(3.85,.65)(3.85,.9)(3.7,1.05)
\psline[linewidth=.5pt,linearc=.2](.8,1.05)(.65,.8)(.65,.55)(4,.55)(4,.9)(3.7,1.05)
\pscircle*(4.05,1.25){.02}
\pscircle*(4.25,1.25){.02}
\pscircle*(4.45,1.25){.02}
%
%\rput[b](3.4,2.35){\small$b_2$}
%\psframe[linewidth=.5pt,framearc=.9](2.4,2.3)(4.6,.6)%(2.55,2.1)(4.45,.9)
%
\endpspicture

%% file: fig3_dl.tex
\pspicture*(6,2.5)
\psline[linewidth=.5pt](1,.25)(5,.25)
\psline(1,2.35)(5,2.35)
\psbezier(1,.25)(0,.25)(0,2.35)(1,2.35)
\psbezier(5,.25)(6,.25)(6,2.35)(5,2.35) %
\pscircle*[linecolor=lightgray](.9,1.25){.3}
\pscircle(.9,1.25){.3}
\rput[t](.9,.9){\tiny$1$}
\pscircle*(1.4,1.25){.02}
\pscircle*(1.6,1.25){.02}
\pscircle*(1.8,1.25){.02}
\pscircle*[linecolor=lightgray](2.3,1.25){.3}
\pscircle(2.3,1.25){.3}
\rput[t](2.3,.9){\tiny$g-l$}
\pscircle*[linecolor=lightgray](3.5,1.25){.3}
\pscircle(3.5,1.25){.3}
%\rput[tl](3.5,.9){\tiny$g-l+1$}
%
\pscircle*(4.05,1.25){.02}
\pscircle*(4.25,1.25){.02}
\pscircle*(4.45,1.25){.02}
\pscircle*[linecolor=lightgray](5,1.25){.3}
\pscircle(5,1.25){.3}
\rput[br](3.25,1.9){\small$d_l$}
\psframe[linewidth=.5pt,framearc=.9](3,.5)(5.5,2)%(2.55,2.1)(4.45,.9)

\endpspicture

%% file: present.tex
\section{\label{pre}The presentation for $\M(F)$}

In \cite{Br} Brown describes a method to produce a presentation of
a group acting on a simply-connected CW-complex. In
\cite{B} Benvenuti uses a special case of Brown's theorem to
obtain a presentation for the orientable mapping class group from
its action on the ordered complex of curves. In this section we
apply the method of \cite{B} to the case of a non-orientable
surface.

\medskip

The following theorem is fundamental for this section.

\begin{theorem}[Ivanov \cite{I}]\label{ivanov}
Let $F=F_g^n$ denote a non-orientable surface of genus $g$ with $n$ boundary
components and $\C(F)$ the complex of curves of $F$. Then
$\C(F)$ is $(g-3)$-connected if $n\in\{0,1\}$, and $(g+n-5)$-connected if
$n\ge 2$.
\foorp\end{theorem}

In particular, except for the surfaces $F_g^n$ where $$(g,n)\in\{(1,n)|n\le 4\}\cup\{(2,n)|n\le 3\}\cup\{(3,n)|n\le 2\}$$
that we call {\it sporadic},
the complex of curves of $F_g^n$ is simply connected.

\medskip

Now we define, following \cite{B}, the {\it ordered complex of
curves} of $F$ denoted by $\Xo$. The $r$-simplices of $\Xo$ are
equivalence classes of {\it ordered} generic $(r+1)$-families of
disjoint curves: $(a_1,\dots,a_r)$ and $(b_1,\dots,b_r)$ represent
the same $(r-1)$-simplex in $\Xo$ if and only if $a_i\simeq b_i^{\pm
1}$ for all $i\in\{1,\dots,r\}$. We denote by
$\left<a_1,\dots,a_r\right>$ the simplex of $\Xo$ represented by the
family $(a_1,\dots,a_r)$. Note that the vertices of $\Xo$ coincide
with those of $\C(F)$ and in general to each $r$-simplex of $\C(F)$
correspond $(r+1)!$ different simplices of $\Xo$ with the same set
of vertices.

\medskip

The following proposition is proved in \cite{B}. The same proof applies to the case
of a non-orientable surface.

\begin{prop}\label{ord}
If $\C(F)$ is simply connected, then $\Xo$ is also simply connected.
\foorp\end{prop}

The mapping class group $\M(F)$ acts on $\Xo$ by $h\langle
a_1,\dots,a_r\rangle=\langle h\circ a_1,\dots,h\circ a_r\rangle.$
Two simplices $\langle a_1,\dots,a_r\rangle$ and $\langle
b_1,\dots,b_r\rangle$ of $\Xo$ are $\M(F)$-equivalent if and only if
the conditions of Proposition \ref{simplices} are satisfied with
$\sigma(i)=i$, $i\in\{1,\dots,r\}$.

\medskip

Let $A=(a_1,\dots,a_r)$ be a generic $r$-family of disjoint
curves. Denote by $\textrm{Stab($\left<A\right>$)}$ the stabilizer
in $\M(F)$ of the simplex of $\Xo$ represented by $A$. The group
$\textrm{Stab($\left<A\right>$)}$ consists of those $h\in\M(F)$
which satisfy $h\circ a_i\simeq a^{\pm 1}_i$ for $i\in\{1,\dots,r\}$.
It is clearly a subgroup of $\rm{Stab}($[A]$)$ and by Proposition
\ref{sequence}, we have the
following exact sequence:
\begin{equation}\label{ordsequence}
1\to{\rm{Stab}}^+([A])\to{\rm{Stab}}(\left<A\right>)\stackrel{\Phi_A}{\to}(\Z_2)^r.
\end{equation}

Here $(\Z_2)^r$ is identified with the subgroup of ${\rm{Cub}}_r$
consisting of those $\phi\in GL(\R^r)$ such that $\phi(e_i)=\pm
e_i$ for all $1\le i\le r$.

\medskip

Denote by $X$ the orbit space $\Xo/\M(F)$ and by $p\colon\Xo\to X$
the canonical projection. The space $X$ inherits from $\Xo$ the
structure of a CW-complex; the $r$-cells of $X$ correspond to the
$\M(F)$-orbits of $r$-simplices of $\Xo$.

By Remark \ref{finite}, $X$ is a finite CW-complex. We denote by
$X^r$ the $r$-skeleton of $X$. Since the edges of $\Xo$ are oriented
and the action of $\M(F)$ preserves the orientation, the edges of
$X$ are also oriented. If $e$ is an edge in either $\Xo$ or $X$ then
we denote by $i(e)$ and $t(e)$ respectively the initial and terminal
vertex of $e$. An edge $e\in X^1$ for which $i(e)=t(e)=v$ is called
a {\it loop} based at $v$.

The advantage of the ordered complex of curves over the ordinary complex of curves is that
$\M(F)$ acts on $(\Xo)^1$ without inversion, which simplifies the statement of Theorem \ref{presentation}
below.

\medskip

\begin{figure}
\begin{center}
\begin{tabular}{cc}
\input{fig4_tri}  & \input{fig5_tri}
\end{tabular}
\caption{\label{tri} A triangle in $X$ and its representative in $\Xo$.}
\end{center}
\end{figure}
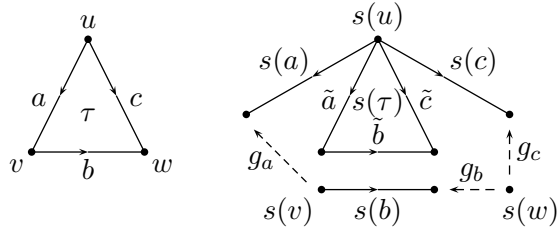

In order to describe a presentation for $\M(F)$ we need to make a number of choices:

(a) We choose a maximal tree $\mathcal{T}$ in $X^1$.

(b) For every $v\in X^0$ we choose a representative $s(v)\in
(\Xo)^0$, and for every $e\in X^1$ a representative $s(e)\in (\Xo)^1$ (i.e. $p(s(v))=v$ and $p(s(e))=e$),
so that $s(i(e))=i(s(e))$ for every $e\in (\Xo)^1$,
and $s(t(e))=t(s(e))$ for every $e\in\mathcal{T}$.
We denote by $S_v$ the stabilizer
$\mathrm{Stab}(s(v))$ and by $S_e$ the stabilizer $\rm{Stab}(s(e))$.

(c) For every $e\in (\Xo)^1$ we choose $g_e\in\M(F)$
such that
$$g_e(s(t(e)))=t(s(e)).$$ If $e\in\mathcal{T}$ then we take $g_e=1$.
Note, then, that the conjugation map $c_e$ given by $g\mapsto g_e^{-1}gg_e$
maps $\rm{Stab}(t(s(e)))$ onto $\rm{Stab}(s(t(e)))$; in particular,
$c_e(S_e)\subseteq S_{t(e)}$.

(d) For every triangle $\tau\in X^2$, with edges
$a$, $b$, $c$ such that $i(c)=i(a)=u$, $t(a)=i(b)=v$,
$t(b)=t(c)=w$, we choose
a representative $s(\tau)$ in $(\Xo)^2$,
such that if $\tilde{a}$, $\tilde{b}$, $\tilde{c}$
are the corresponding edges of $s(\tau)$, then $i(\tilde{c})=i(\tilde{a})=s(u)$ (see Figure \ref{tri}).
We also choose three elements
$$h_{\tau,a}\in S_u,\quad h_{\tau,b}\in S_v,\quad h_{\tau,c}\in S_w,$$
such that $h_{\tau,a}(s(a))=\tilde{a}$,
$h_{\tau,a}g_ah_{\tau,b}(s(b))=\tilde{b}$,
$h_{\tau,a}g_ah_{\tau,b}g_bh_{\tau,c}g^{-1}_c(s(c))=\tilde{c}$. Let
us define $h_\tau=h_{\tau,a}g_ah_{\tau,b}g_bh_{\tau,c}g^{-1}_c.$
Observe, that $h_{\tau}\in S_u$.

The next result is a special case of a general theorem of Brown \cite{Br}
(c.f Theorem 3 of \cite{B}).

\begin{theorem}\label{presentation}
Suppose that $F$ is not sporadic and:\\
(i) for each $v\in X^0$ the group $S_v$ has the
presentation
$S_v=\left< G_v\,|\,R_v\right>,$\\
(ii) for each $e\in X^1$ the stabilizer $S_e$
is generated by $G_e$.

Then $\M(F)$ admits the presentation:
\begin{eqnarray*}
\textrm{generators}&=&\bigcup_{v\in X^0}G_v\cup
\{g_e\,|\,e\in X^1\},\\
\textrm{relations}&=&
\bigcup_{v\in X^0}R_v \cup R^{(1)}\cup R^{(2)}\cup R^{(3)},
\end{eqnarray*}
where:\\
$R^{(1)}=\{g_e=1\,|\,e\in\mathcal{T}\}$.\\
$R^{(2)}=\{g_e^{-1}i_e(g)g_e=c_e(g)\,|\,
g\in G_e, e\in X^1\}$, where $i_e$ is the inclusion $S_e\hookrightarrow S_{i(e)}$ and
$c_e\colon S_e\to S_{t(e)}$ is as in {\rm (c)} above.\\
$R^{(3)}=\{h_{\tau,a}g_ah_{\tau,b}g_bh_{\tau,c}g_c^{-1}=h_\tau
\,|\,\tau\in X^2\}$. \foorp\end{theorem} In Theorem
\ref{presentation}, $i_e(g)$, $c_e(g)$, $h_{\tau,a}$,
$h_{\tau,b}$, $h_{\tau,c}$ and $h_{\tau}$ should be expressed as
words in the generators $\bigcup_{v\in X^0}G_v.$

Suppose that two of the edges of a triangle $\tau\in(X)_2$ belong
to the maximal tree $\mathcal{T}$. Then, using the relations
$R^{(1)}$ and $R^{(3)}$ we can express the generating symbol
corresponding to the third edge as a product of stabilizers of the
representatives for the vertices. The same is true if two of the
symbols for the edges were already expressed as products of
stabilizers. We say that a symbol $g_e$ is {\it determinable} (or
simply that the corresponding edge $e$ is determinable), if using
recursively relations $R^{(1)}$ and $R^{(3)}$, it is possible to
express $g_e$ as a product of elements in $\bigcup_{v\in X^0}G_v.$
Thus, every edge $e\in\mathcal{T}$ is determinable, and if a
triangle in $X^2$ has two determinable edges, then its third edge
is also determinable.

\begin{theorem}\label{determinable}
Suppose that $F^n_g$ is not sporadic.
Then there exists a choice of
the maximal tree $\mathcal{T}$ such that all the edges of $X$
are determinable.
\end{theorem}

\proof

We fix boundary curves $c_1,\dots,c_n$.
For each generic family of disjoint curves $A$ we identify
a generic curve $b$ in $F_A$ with the curve $\rho_A\circ b$
in $F$. For any surface $X$, we denote by $g(X)$ its genus.

\medskip

{\bf Construction of $\mathcal{T}$ for $g\ge 4$.}
Suppose that $g\ge 4$. Let $v_1$ denote the non-separating,
two-sided vertex
$v_1=p([a])$, where $F_{a}$ is
non-orientable.
For each vertex $v$ different from $v_1$, we define an edge $e_v\in X^1$
with initial vertex $v_1$ and terminal vertex $v$ as follows.
We fix a curve $b$, such that $p([b])=v$ and construct $a$ in $F_b$, such that
$p([a])=v_1$. We consider cases.

{\it Case 1:} $b$ is non-separating and $F_b$ is non-orientable.
Since $v\ne v_1$, $b$ must be one-sided and from the comparison of
Euler characteristics we know that $g(F_b)\ge 3$. We define $a$ to
be any two-sided and non-separating curve in $F_b$, such that
$F_{(a,b)}$ is non-orientable.

{\it Case 2:} $b$ is non-separating and $F_b$ is orientable. Now
$F_b$ has genus at least $1$ and hence it contains a
non-separating curve. Let $a$ be any such curve. Note that $F_a$
is non-orientable because we can construct a one-sided curve in
$F_a$ by connecting two boundary points of $F_{(a,b)}$ by an arc.

{\it Case 3:} $b$ is separating, $F_b=N\amalg N'$.
We consider two sub-cases.

{\it Case 3a:} one of the components, say $N$, is orientable. If
$g(N)\ge 1$ then we define $a$ to be any non-separating curve in $N$
(note that $N'$ is non-orientable, and hence so is $F_a$). If
$g(N)=0$, then we define $a$ to be any non-separating, two-sided
curve in $N'$, such that $N'_a$ is non-orientable (such curve
exists, as $g(N')=g\ge 4$).

{\it Case 3b:} both components $N$ and $N'$ are non-orientable.
Assume $g(N)\ge g(N')$. If $g(N)=g(N')$ and $n\ge 1$, then we assume
that $N$ contains the boundary curve $c_1$. If $g(N)\ge 3$ then we
define $a$ to be any non-separating, two-sided curve in $N$, such
that $N_a$ is non-orientable. If $g(N)=2$, then we choose for $a$
any non-separating, two-sided curve in $N$, such that all exterior
boundary curves of $N$ induce the same orientation of $N_a$. If $F$
is closed and $g(N)=g(N')$, then we can not distinguish between $N$
and $N'$. However, whether we choose $a$ in $N$ or $N'$, we obtain
$\M(F)$-equivalent edges $\lr{a,b}$.

In each case we have $p([a])=v_1$ and
we define $e_v=p(\left<a,b\right>)$. By Proposition
\ref{simplices} this definitions do not depend on the choices of
the curves $a$ and $b$. We define the maximal tree
$\mathcal{T}=\{e_v\,|\,v\ne v_1\}$.

\begin{rem}\label{TforClosed}
Suppose that $F$ is closed and consider the curves $a_1$, $a_2$,
$a_3$, $b_k$, $d_l$ in Figures \ref{ai} and \ref{bkdl}. As it was
discussed in Section \ref{orbits}, these curves represent all
vertices of $X$. Clearly $p([a_1])=v_1$ and in the construction of
the maximal tree described above we can take $b$ to be $a_2$ (case
1), $a_3$ (case 2), $b_k$ (case 3a) or $d_l$ (case 3b). Then, in
each case, we can take $a=a_1$. Thus
$$\T=\{\p{a_1,a_2}, \p{a_1,a_3}, \p{a_1,b_k}, \p{a_1,d_l}\,|\,2\le k+1,l\le\frac{g}{2}\}.$$
\end{rem}

\begin{lemma}\label{g>3}
Suppose that $g\ge 4$ and $\mathcal{T}$ is defined as above. Then the following
edges of $X$ are determinable:\\
(i) all the loops based at $v_1$;\\
(ii) all the edges with both ends in non-separating vertices;\\
(iii) all the edges with one end in a non-separating vertex and the other end in a separating vertex;\\
(iv) all the edges with both ends in separating vertices.
\end{lemma}
\proof
Let $e=p(\left<a,b\right>)$ be any edge in $X$ and let $F'$ denote the surface
$F_{(a,b)}$.

(i) Suppose $p([a])=p([b])=v_1$. The surface $F'$
is either connected or it has two connected components, at least one
of which must be non-orientable.

Suppose that $F'$ has a non-orientable connected component of genus
at least 2 or it has two non-orientable components. Then
there exists a one-sided curve $c$ in $F'$ such that
$F_{(a,c)}$ and $F_{(b,c)}$ are non-orientable. By the definition of edges $e_v$ (case 1), we have
that $p(\left<a,c\right>)=p(\left<b,c\right>)=e_{p([c])}$, the triangle $p(\left<a,b,c\right>)$ has two edges
in $\T$, and thus $e$ is determinable.

Suppose now that $F'$ is connected and orientable. Let $a'$,
$a''$, $b'$, $b''$ denote the boundary curves of $F'$ such that
$\rho_{(a,b)}\circ a'=\rho_{(a,b)}\circ a''=a$, $\rho_{(a,b)}\circ
b'=\rho_{(a,b)}\circ b''=b$. Let $c$ be a separating curve in $F'$
such that $\{a',b'\}$ and $\{a'',b''\}$ are in different
components of $F'_c$. Observe that $c$ is non-separating in $F$.
Every one-sided curve in $F$ intersects $a\cup b$ odd number
of times, thus it intersects $c$. Hence $F_c$ is orientable and
$p([c])\ne v_1$. The triangle $\p{a,b,c}$ has edges $e,
e_{p([c])}, e_{p([c])}$ (case 2 in the construction of $\mathcal{T}$), thus $e$ is determinable.

Finally suppose that $F'$ has two components
$N_1$ and $N_2$, such that $N_1$ is non-orientable of genus $1$ and
$N_2$ is orientable. Since $g(N_2)\ge 1$, there is a non-separating
two-sided curve $c$ in $N_2$. Note that $p([c])=v_1$ and the loops
$p(\left<a,c\right>)$, $p(\left<b,c\right>)$ are determinable by
previous arguments, because $F_{(a,c)}$ and $F_{(b,c)}$ are connected. Hence $e$ is also determinable,
by $\p{a,b,c}$.

\medskip

(ii) Suppose that both ends of $e$ are non-separating.
If both of them are one-sided, then $F'$ is connected and has genus at least
$1$ if it is orientable, or at least $2$ if it is
non-orientable. In both cases $F'$ contains a non-separating, two-sided curve $c$.
Now $\p{c,a}=e_p([a])$, $\p{c,b}=e_p([b])$ (case 1 in the construction of $\mathcal{T}$),
hence $e$ is determinable by $\p{c,a,b}$.

Suppose that one end of $e$ is one-sided and the other one is two-sided.
Then $F'$ is connected and the two-sided end is $v_1$. Thus if $b$ is one-sided, then
$e=e_{p([b])}$. If $a$ is one-sided, then we choose any separating curve $c$ in $F'$, such
that $F_c$ is connected. Now $\p{c,a}=e_{p([a])}$ and $\p{c,b}$ is a loop at $v_1$, which is
determinable by (i). Hence $e$ is determinable by $\p{c,a,b}$.

Suppose that both ends of $e$ are two-sided. We can assume that at
least one of the ends is not $v_1$, so $F'$ is orientable. If $F'$
is connected, then we choose a separating generic curve $c$ in
$F'$, such that $F_{(a,c)}$ and $F_{(b,c)}$ are connected. Now
$\p{c,a}$ is either $e_{p([a])}$ (if $F_a$ is orientable) or a
loop at $v_1$ (if $F_a$ is non-orientable) and similarly for
$\p{c,b}$. Hence $e$ is determinable by $\p{c,a,b}$. If $F'$ is
not connected, then $F_a$ and $F_b$ are orientable. Now $F'$ has a
component $N$ with $g(N)\ge1$ and for any non-separating curve $c$
in $N$ we have $\p{c,a}=e_{p([a])}$ and $\p{c,b}=e_{p([a])}$.
Hence $e$ is determinable by $p(\left<c,a,b\right>)$.

\medskip

(iii) Assume, without loss of generality, that $a$ is separating and
$b$ is non-separating. Suppose that both components of $F_a$ have
genus $\ge 1$. Let $a_1$ be a generic curve in $F'$ such that
$\p{a_1,a}\in\T$, and choose any non-separating curve $c$ in the
other component of $F_a$. Notice that $\p{c,a_1}$ is determinable by
(ii), and $\p{c,a}$ is determinable by the triangle $\p{c,a_1,a}$.
Now if $a_1$ and $b$ belong to different components of $F_a$, then
$\p{a_1,b}$ is determinable by (ii), and $e$ is determinable by
$\p{a_1,a,b}$. If $a_1$ and $b$ belong to the same component of
$F_a$, then $e$ is determinable by $\p{c,a,b}$. If one of the
components has genus $0$, then $b$ is contained in the other
component $N$. Now there exists a two-sided generic curve $a_1$ in
$N_b$, such that $N_{a_1}$ is connected and non-orientable. Indeed,
if $N_b$ is orientable, then $g(N_b)\ge1$ and $a_1$ may be any
non-separating curve in $N_b$. If $N_b$ is non-orientable, then
$g(N_b)\ge2$ and we may take $a_1$ to be separating in $N_b$. For
such $a_1$ we have $\p{a_1,a}\in\T$, and $\p{a_1,b}$ is determinable
by (ii). Hence $e$ is determinable by $\p{a_1,a,b}$.

\medskip

To prove (iv) notice that in this case $F'$ must have a non-orientable component. Choose a
one-sided curve $c$ in $F'$ and consider the triangle $p(\left<c,a,b\right>)$.
The assertion follows by (iii).
\foorp

This finishes the proof of Theorem \ref{determinable} for $g\ge 4$.

\medskip

{\bf Construction of $\mathcal{T}$ for $g=3$.} Suppose that $g=3$.
Since $F$ is not sporadic we have $n\ge 3$. Let $v_1$ denote the
non-separating, two-sided vertex $p([a])$, where $F_a$ is non-orientable.
Note that this is the only non-separating, two-sided vertex in $X$.
As we did for $g\ge 4$, for each $v\ne v_1$ we define an edge $e_v$ form $v_1$ to $v$.
We fix $b$ such that $v=p([b])$ and define
$a$ in $F_b$ so that $p([a])=v_1$.

{\it Case 1:} $b$ is one-sided and $F_b$ is non-orientable. Now $F_b$ has genus
$2$. We define $a$ to be any two-sided and non-separating curve in $F_b$, such that
all exterior boundary curves induce the same orientation of $F_b$.

If  $F_b$ is connected and orientable ({\it case 2}) or disconnected
({\it case 3}), then we define $a$ in the same way as we did for
$g\ge 4$. We only remark that in case 2, $b$ is one-sided and hence
$F_b$ has genus $1$; and in case 3a, if $g(N)=0$ then $g(N')=3$,
which suffices to choose two-sided and non-separating $a$ with
$N'_a$ non-orientable.

 As previously we
define $e_v=\p{a,b}$ and $\T=\{e_v\,|\,v\ne v_1\}$.

\begin{lemma}\label{g=3}
Suppose that $g=3$ and $\mathcal{T}$ is defined as above.
Then the following
edges of $X$ are determinable:\\(i) all the loops based at $v_1$;\\
(ii) all the edges with one end in $v_1$;\\
(iii) all the edges with at least one edge in one-sided vertex;\\
(iv) all the edges with both ends in separating vertices.
\end{lemma}

\proof

First observe that every edge in $X$ satisfies one of the conditions
(i)--(iv). Therefore Lemma \ref{g=3} implies Theorem
\ref{determinable} for $g=3$.

 Let $e=p(\left<a,b\right>)$ be any edge in $X$
and $F'=F_{(a,b)}$

(i) If $p([a])=p([b])=v_1$ then $F'$ has two connected components, at least one
of which contains two exterior boundary curves. Let $c$ be a curve in $F'$ bounding a pair of pants
together with two exterior boundary curves. The edge $e$ is determinable by the triangle $p(\left<a,b,c\right>)$
having two edges in $\mathcal{T}$.

\medskip

(ii) Assume $p([a])=v_1$. If $F_b$ is connected and orientable or it
has an orientable component, then $e\in\T$. In the other case
$e\in\T$ if and only if all exterior boundary curves induce the same
orientation of the orientable component of $F'$. Suppose that
$e\notin\T$. Denote by $N$ the connected component of $F_b$ having
genus $2$ and by $N'$ the orientable component of $F'$ (thus
$N'=N_a$). There exists a separating curve $c$ in $N'$, which is
non-separating in $N$ and such that any two exterior boundary curves
induce opposite orientations of $N'$ if and only if they belong to
different components of $N'_c$. The surface $N_c$, which can be
obtained from $N'_c$ by gluing along $a$, is the orientable
component of $N_{(b,c)}$. Note that all exterior boundary curves
induce the same orientation of $N_c$, hence $\p{c,b}=e_{p([b])}$.
The loop $\p{c,a}$ is determinable by (i), thus $e$ is determinable
by $\p{c,a,b}$.

Now assume $p([b])=v_1$ and choose any generic curve $d$ in $F'$.
The edges $\p{b,a}$ and $\p{b,d}$ have initial vertex $v_1$ and we
have already proved that such edges are determinable. Hence
$\p{a,d}$ is determinable by $\p{b,a,d}$, and $e$ by $\p{a,b,d}$.

\medskip

(iii) Suppose that $e$ has both ends in one-sided vertices.
Choose any curve $c$ in $F'$ bounding a
pair of pants together with two exterior boundary curves. Let $d$
be any two-sided non-separating curve in $F_{(a,c)}$. Then
$p([d])=v_1$, and $\p{d,c}$ and $\p{d,a}$ are determinable by
(ii), thus $\p{c,a}$ is determinable by $\p{d,c,a}$. Analogously
$\p{c,b}$ is determinable by a different triangle $\p{d',c,b}$.
Finally $e$ is determinable by $\p{c,a,b}$.

Suppose that $e$ has one vertex in a one-sided vertex $v$ and the
other end in a separating vertex. Assume without
loss of generality  that $a$ is separating and denote by $N$ the
component of $F_a$ which contains $b$, and the other component by
$N'$. If $g(N)=3$ or $g(N)=1$, then
$F'$ contains a
non-separating two-sided curve $c$ and $e$ is determinable by
$\p{c,a,b}$ and (ii).
If $g(N)=2$, then we choose a one-sided
curve $d$ in $N'$ and two-sided, non-separating
curve $c$ in $N$. Now $\p{a,d}$ is determinable by $\p{c,a,d}$ and (ii), and
$\p{b,d}$ is an edge with two one-sided ends, determinable by previous argument.
Finally $e$ is determinable
by $\p{a,b,d}$.

\medskip

(iv) If $e$ has both ends in separating vertices then
$F'$ has a non-orientable connected
component. Choose a one-sided curve $c$ in $F'$
and consider the triangle $p(\left<a,b,c\right>)$. The assertion
follows by (iii). \foorp

\medskip

{\bf Construction of $\mathcal{T}$ for $g=2$.} Suppose that $g=2$
and $n\ge 4$. Let $v_2$ denote the unique one-sided
vertex of $X$. For each separating vertex $v$ we will define an
edge $e_v\in X^1$ from $v_2$ to $v$.
We fix $b$ such that $p([b])=v$ and assume $F_b=N\amalg N'$.
We define $e_v=\p{a,b}$, where $a$ is
a one-sided curve in $F_b$ defined as follows.

{\it Case 1:} one component of $F_b$, say $N$, is orientable.
Then we define $a$ to be any one-sided curve in $N'$.

{\it Case 2:} both components are non-orientable. Assume that $N$
contains the exterior boundary curve $c_1$. We choose $a$ in $N$,
so that all exterior boundary curves of $N$ induce the same orientation
of $N_a$.

Suppose that  $w$ is a two-sided, non-separating vertex of $X$.
Let us choose $b$ such that $p([b])=w$. Now $F_b$ is orientable and
has genus $0$.
We choose a curve $a$ in $F_b$ bounding a
pair of pants together with the exterior boundary curves $c_1$ and $c_2$.
We define $e_w=\p{a,b}$.

We claim that $\T=\{e_v\,|\,v\ne v_2\}$ is a maximal tree in
$X^1$. First notice that $\T'=\{e_v\,|\,\textrm{$v$\ is\
separating}\}$ is a tree, because every edge $e_v\in\T'$ connects
$v$ to $v_2$. Now $\T\backslash\T'=\{e_w\,|\,\textrm{$w$\ is\
two-sided\ and\ non-separating}\}$ and every two-sided and
non-separating vertex $w$ is connected to exactly one vertex of
$\T'$ by $e_w$. It follows that $\T$ indeed is a tree and since it
contains all vertices of $X$ it is a maximal tree.

\begin{lemma}\label{g=2}
Suppose that $g=2$ and $\mathcal{T}$ is defined as above.
Then the following
edges of $X$ are determinable:\\
(i) all the loops based at $v_2$;\\
(ii) all the edges with one end in $v_2$;\\
(iii) all the edges with both ends in two-sided vertices;
\end{lemma}

\proof
Let $e=p(\left<a,b\right>)$ be any edge of $X$ and $F'=F_{(a,b)}$.

(i) Suppose that $p([a])=p([b])=v_2$. Choose any separating generic curve
$c$ in $F'$ such that one component of $F_c$ is orientable.
Then $\p{a,c}=\p{b,c}=e_{p([c])}$ and hence $e$ is determinable by the triangle
$\p{a,b,c}$.

\medskip

(ii) Suppose that $e$ has one end in $v_2$ and the other end in a
separating vertex $v$. Assume without loss of generality, that $a$
is separating. If $F_a$ has an orientable component then for each
one-sided curve $c$ in $F'$ we have $\p{c,a}\in\T$. Now $\p{c,b}$
is determinable by (i), hence $e$ is determinable by $\p{c,a,b}$.
Suppose that both components of $F_a$ are non-orientable and let
$c$ and $d$ be two one-sided curves in different components of $F_a$, such that
$\p{c,a}=e_v$. Since $\p{c,d}$ is determinable by (i), $\p{d,a}$
is determinable by $\p{c,d,a}$. We have $b\cap c=\emptyset$ or
$b\cap d=\emptyset$, hence $e$ is determinable by $\p{c,a,b}$ or
$\p{d,a,b}$.

\medskip

(iii) If both ends of $e$ are separating, then there is a one-sided curve $c$ in
$F'$ and $e$ is determinable  by (ii).
Suppose that $e$ has one separating and one non-separating end.
Assume without loss of generality, that $a$ is non-separating.
Then there is a separating generic curve $c$ in $F'$ such that all
boundary curves of $F$ are contained in one connected component of
$F_c$. In particular, there is a curve $d$ in $F_{(a,c)}$
bounding a pair of pants together with $c_1$ and $c_2$, that is
$\p{a,d}=e_{p([a])}$. The edge $\p{c,d}$ is determinable by the
previous argument, hence $\p{a,c}$ is determinable by $\p{a,c,d}$.
If $c\simeq b^{\pm 1}$ then we can assume $b\cap d=\emptyset$, and
$e$ is determinable by $\p{a,b,d}$. In the other case $e$ is
determinable $\p{a,b,c}$. Finally suppose that both ends of $e$
are non-separating. Since $n\ge 4$, $F'$ contains a  generic curve
$b'$, and $e$ is determinable by $\p{a,b,b'}$.
\foorp

\medskip

%%%%%%%%%%%%%%%%%%%%%%%%%%%%%%%%%%%%%%%%%%%%%%%%%%%%%%%%%%%%%%%%%%%%%%%%%%%%%%%%%%%%%%%%%%%%%%%%%%%%
%%%%%%%%%%%%%%%%%%%%%%%%%%%%%%%%%%%%%%%%%%%%%%%%%%%%%%%%%%%%%%%%%%%%%%%%%%%%%%%%%%%%%%%%%%%%%%%%%%%%
%%%%%%%%%%%%%%%%%%%%%%%%%%%%%%%%%%%%%%%%%%%%%%%%%%%%%%%%%%%%%%%%%%%%%%%%%%%%%%%%%%%%%%%%%%%%%%%%%%%%

{\bf Construction of $\mathcal{T}$ for $g=1$.}

Suppose that $g=1$ and $n\ge 5$.
It follows from Proposition \ref{simplices} that each separating vertex
$p([a])\in X^0$ is uniquely determined by a pair $I,J\subset\{1,\dots,n\}$
such that $I\cap J=\emptyset$, $2\le\#I+\#J\le(n-1)$, and if $N$ is the orientable connected
component of $F_a$ then

\begin{itemize}
\item $c_i$ is a boundary curve of $N$ if and only if $i\in I\cup
J$, \item $c_i$ and $c_j$ induce the same orientation of $N$ if
and only if $\{i,j\}\subseteq I$ or $\{i,j\}\subseteq J$.
\end{itemize}

We denote such vertex by $v_{I,J}$, where we assume $\#I\le\#J$, and if $\#I=\#J$ then
$\mathrm{min}\,I<\mathrm{min}\,J$.
Each one-sided vertex $p([a])$ is uniquely determined by a subset
$I\subset\{1.\dots,n\}$ such that $c_i$ and $c_j$ induce the same
orientation of $F_a$ if and only if $\{i,j\}\subseteq I$ or
$\{i,j\}\subseteq I'$, where $I'=\{1,\dots,n\}\backslash I$. We
denote such vertex by $v_I$, where we assume $\#I\le n/2$, and if
$\#I=n/2$ then $1\in I$ (see Figure \ref{X1}, where we assume that
all boundary curves have positive orientations with respect to the
standard orientation of the plane of the figure).

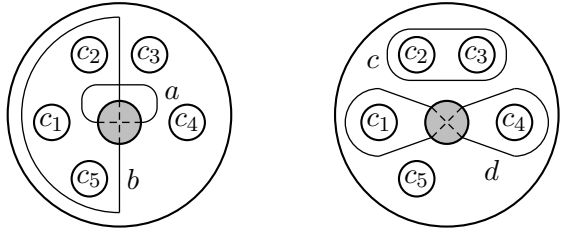
\begin{figure}%[!htbp]
\begin{center}
\begin{tabular}{cc}
\input{fig6_X1} & \input{fig7_X1}
\end{tabular}
\caption{Representatives of different vertices of the complex $X$ of
$F_1^5$: \label{X1}$p([a])=v_\emptyset$, $p([b])=v_{\{3,4\}}$,
$p([c])=v_{\emptyset,\{2,3\}}$, $p([d])=v_{\{1\},\{4\}}$.}
\label{X1}
\end{center}
\end{figure}

If $\#I+\#J\le\#K+\#L$ then $v_{I,J}$ and $v_{K,L}$ are connected by an edge in $X$ if and only if
one of the following conditions is satisfied:
\begin{itemize}
\item $I\subseteq K$, $J\subseteq L$, $\quad \#I+\#J<\#K+\#L$;
\item $I\subseteq L$, $J\subseteq K$, $\quad \#I+\#J<\#K+\#L$;
\item $(I\cup J)\cap(K\cup L)=\emptyset$.
\end{itemize}
Vertices $v_{I}$ and $v_{J,K}$ are connected by an edge if and
only if either $J\subseteq I$, $K\subseteq I'$ or $K\subseteq I$,
$J\subseteq I'$. There are no edges connecting two one-sided vertices
because every two one-sided curves in a surface of genus $1$ intersect.
It follows that $X$ has no loops. Moreover, it
follows from Proposition \ref{simplices} that for each pair
$v,w\in X^0$ there is at most one edge in $X^1$ with initial
vertex $v$ and terminal vertex $w$. If such edge exists, then we
denote it by $\lr{v;w}$. If every two of three vertices $u$, $v$,
$w$ are connected by an edge in $X$, then there are $6$ triangles
in $X^2$ with vertices $u$, $v$, $w$. We denote by $\lr{u;v;w}$
the triangle with edges $\lr{u;v}$, $\lr{u;w}$, $\lr{v;w}$.

We define the maximal tree as
$$\T=\bigcup_{v_{I,J}\in X^0}\left\{\lr{v_{I};v_{I,J}}\right\} \cup
\bigcup_{v_{I}\in X^0\backslash\{v_\emptyset\}}\left\{\lr{v_I;v_{\emptyset,I'}}\right\}.$$

\begin{lemma}\label{g=1}
Suppose that $g=1$ and $\mathcal{T}$ is defined as above.
Then the following
edges of $X$ are determinable:\\
(i) all edges with ends in $v_{I,J}$ and $v_{K,L}$, where $I\subseteq K$, $J\subseteq L$;\\
(ii) all edges with ends in $v_{I,J}$ and $v_{K,L}$, where $(I\cup J)\cap(K\cup L)=\emptyset$;\\
(iii) all edges with ends in $v_{I,J}$ and $v_{K,L}$, where $I\subseteq L$, $J\subseteq K$;\\
(iv) all edges with ends in $v_{I,J}$ and $v_{K}$.
\end{lemma}

\proof
Let $e$ be an edge with ends in $v_{I,J}$ and $v_{K,L}$.

(i) If $I=K$ then $e$ is determinable by a triangle with third
vertex $v_I$. Suppose $I\subsetneq K$, $J=L$. The edge
$\lr{v_{\emptyset,J};v_{\emptyset,K'}}$ is determinable by the
previous argument, hence $\lr{v_K;v_{\emptyset,J}}$ is determinable
by $\lr{v_K;v_{\emptyset,J};v_{\emptyset,K'}}$. If $I=\emptyset$
then $e$ is determinable by the triangle with edges $e$,
$\lr{v_K;v_{\emptyset,J}}$ and $\lr{v_K;v_{K,J}}$. If
$I\ne\emptyset$ then $e$ is determinable by the triangle with edges
$e$, $\lr{v_{I,J};v_{\emptyset,J}}$, $\lr{v_{K,J};v_{\emptyset,J}}$,
whose last two edges are determinable by the previous argument.
Finally, if $I\subsetneq K$ and $J\subsetneq L$ then $e$ is
determinable by the triangle with edges $e$, $\lr{v_{I,J};v_{I,L}}$,
$\lr{v_{K,L};v_{I,L}}$, because the last two edges are determinable
by previous arguments.

\medskip

(ii) If $\#(I\cup J\cup K\cup L)<n$ then $e$ is determinable by a
triangle with third vertex $v_{I\cup K,J\cup L}$, whose remaining
two edges are determinable by (i). If $\#(I\cup J\cup K\cup L)=n$
then we assume $\#(I\cup J)\ge 3$. Then there is a vertex $v_{M,N}$
such that $M\subseteq I$, $N\subseteq J$ and $\#(M\cup N)<\#(I\cup
J)$. Now $\lr{v_{M,N};v_{K,L}}$ is determinable by the previous
argument, and $\lr{v_{M,N};v_{I,J}}$ is determinable by (i). Hence
$e$ is also determinable.

\medskip

(iii) Suppose $J=K$ , $I\subsetneq L$. If $\#J\ge 2$ then the edges
$\lr{v_{\emptyset,L};v_{J,L}}$ and
$\lr{v_{\emptyset,J};v_{\emptyset,L}}$ are determinable by (i) and
(ii), hence any edge connecting $v_{\emptyset,J}$ with $v_{J,L}$ is
determinable. In particular, $e$ is determinable if $I=\emptyset$,
and if $I\ne\emptyset$ then $e$ is determinable by the triangle with
edges $e$, $\lr{v_{J,L};v_{\emptyset,J}}$,
$\lr{v_{I,J};v_{\emptyset,J}}$, whose last edge is determinable by
(i). Suppose $\#J=1$. Then $\#I=1$ and $\#L\ge 2$. Now $e$ is
determinable by a triangle with third vertex $v_{\emptyset,M}$,
where $M=L\backslash I$ if $\#L\ge 3$, and $M=(J\cup L)'$ if $\#L=
2$ ($\#M\ge 2$, since $n\ge 5$). In both cases $e$ is determinable
by (i) and (ii). Finally, if $J\subsetneq K$ and $I\subsetneq L$
then $e$ is determinable by the triangle with edges $e$,
$\lr{v_{I,J};v_{J,L}}$, $\lr{v_{K,L};v_{J,L}}$, because
$\lr{v_{I,J};v_{J,L}}$ is determinable by previous arguments, and
$\lr{v_{K,L};v_{J,L}}$ by (i).

\medskip

(iv) First assume $K=\emptyset$. Then $I=\emptyset$ and if
$v_K=i(e)$ then $e\in\T$. Suppose $v_K=t(e)$. Observe that there is
a vertex $v_{\emptyset,L}$ such that $L\subsetneq J$ or $J\subsetneq
L$. Now $e$ is determinable by
$\lr{v_{\emptyset,J};v_\emptyset;v_{\emptyset,L}}$. Now assume
$K\ne\emptyset$ and $\#J\ge 2$. Any edge connecting $v_{K}$ with
$v_{\emptyset,J}$ is determinable by a triangle with third vertex
$v_{\emptyset,K'}$. In particular, $e$ is determinable if
$I=\emptyset$, and if $I\ne\emptyset$ then $e$ is determinable by
the triangle with edges $e$, $\lr{v_K;v_{\emptyset,J}}$,
$\lr{v_{I,J};v_{\emptyset,J}}$, whose last edge is determinable by
(i). It remains to consider the case $\#I=\#J=1$. It is easy to
check that then there is a triangle with vertices $v_K$, $v_{I,J}$,
$v_{L,M}$, where $I\cup J\subsetneq L\cup M$. The edge connecting
$v_K$ with $v_{L,M}$ is determinable by the previous argument, hence
$e$ is also determinable.

This completes the proof of Lemma \ref{g=1} and Theorem \ref{determinable}
\foorp

We a corollary we obtain the following theorem.

\begin{theorem}\label{reduced}
Suppose that $F=F_g^n$ is not sporadic and
$\mathcal{T}$ is as in Lemma \ref{determinable}. Then it is possible
to express all the generators $g_e$ appearing in Theorem
\ref{presentation} as a product of elements in $\bigcup_{v\in
X^0} G_v$. Hence, the presentation in Theorem
\ref{presentation} reduces to
\begin{equation*}
\M(F)=\langle\bigcup_{v\in X^0}G_v\, | \bigcup_{v\in
X^0}R_v\,\cup\widetilde{R^{(2)}} \cup\widetilde{R^{(3)}}\,\rangle,
\end{equation*}
where $\widetilde{R^{(i)}}$ are the relations obtained
substituting in $R^{(i)}$ the expressions for the generators
$g_e$.\foorp
\end{theorem}

%% file: fig4_tri.tex
\pspicture*(2.5,3)
\rput[tr](.4,.9){\small$v$}
\pscircle*(.5,1){.05}
\rput[tl](2.1,.9){\small$w$}
\pscircle*(2,1){.05}
\rput[b](1.25,2.65){\small$u$}
\pscircle*(1.25,2.5){.05}
\rput(1.25,1.5){\small$\tau$}
\rput[r](.7,1.7){\small$a$}
\psline[linewidth=.5pt]{->}(1.25,2.5)(.875,1.75)
\psline[linewidth=.5pt](.875,1.75)(.5,1)
\rput[l](1.8,1.7){\small$c$}
\psline[linewidth=.5pt]{->}(1.25,2.5)(1.625,1.75)
\psline[linewidth=.5pt](1.625,1.75)(2,1)
\rput[t](1.25,.9){\small$b$}
\psline[linewidth=.5pt]{->}(.5,1)(1.25,1)
\psline[linewidth=.5pt](1.25,1)(2,1)
\endpspicture

%% file: fig5_tri.tex
\pspicture*(4.75,3)
\pscircle*(1.5,1){.05}
\pscircle*(3,1){.05}
\pscircle*(2.25,2.5){.05}
\rput(2.25,1.65){\small$s(\tau)$}
\rput[r](1.7,1.7){\small$\tilde{a}$}
\psline[linewidth=.5pt]{->}(2.25,2.5)(1.875,1.75)
\psline[linewidth=.5pt](1.875,1.75)(1.5,1)
\rput[l](2.8,1.7){\small$\tilde{c}$}
\psline[linewidth=.5pt]{->}(2.25,2.5)(2.625,1.75)
\psline[linewidth=.5pt](2.625,1.75)(3,1)
\rput[b](2.25,1.1){\small$\tilde{b}$}
\psline[linewidth=.5pt]{->}(1.5,1)(2.25,1)
\psline[linewidth=.5pt](2.25,1)(3,1)
\rput[t](2.25,.4){\small$s(b)$}
\pscircle*(1.5,.5){.05}
\pscircle*(3,.5){.05}
\psline[linewidth=.5pt]{->}(1.5,.5)(2.25,.5)
\psline[linewidth=.5pt](2.25,.5)(3,.5)
\rput[br](1.35,2){\small$s(a)$}
\pscircle*(.5,1.5){.05}
\psline[linewidth=.5pt]{->}(2.25,2.5)(1.375,2)
\psline[linewidth=.5pt](1.375,2)(.5,1.5)
\rput[bl](3.2,2){\small$s(c)$}
\pscircle*(4,1.5){.05}
\psline[linewidth=.5pt]{->}(2.25,2.5)(3.125,2)
\psline[linewidth=.5pt](3.125,2)(4,1.5)
\rput[t](.7,1){\small$g_a$}
\psline[linewidth=.5pt,linestyle=dashed,dash=3pt 2pt]{->}(1.3,.6)(.6,1.3)
\rput[b](3.5,.6){\small$g_b$}
\psline[linewidth=.5pt,linestyle=dashed,dash=3pt 2pt]{->}(3.8,.5)(3.2,.5)
\rput[l](4.1,1){\small$g_c$}
\psline[linewidth=.5pt,linestyle=dashed,dash=3pt 2pt]{->}(4,.7)(4,1.3)
%\rput[br](.95,3.55){\small$s(v)$}
%\pscircle*(1,3.5){.05}
%\rput[bl](5.55,3.55){\small$s(u)$}
%\pscircle*(5.5,3.5){.05}
\rput[b](2.25,2.6){\small$s(u)$}
\rput[tr](1.4,.4){\small$s(v)$}
\rput[tl](3.9,.4){\small$s(w)$}
\pscircle*(4,.5){.05}
%
%\rput[b](1.75,3.6){\small$h_{a,t}$}
%\psline[linewidth=.5pt,linestyle=dashed,dash=3pt 2pt]{->}(2.25,3.5)(1.25,3.5)
%\rput[b](4.75,3.6){\small$h_{a,i}$}
%\psline[linewidth=.5pt,linestyle=dashed,dash=3pt 2pt]{->}(4.25,3.5)(5.25,3.5)
%\rput[l](1.1,3){\small$h_{b,i}$}
%\psline[linewidth=.5pt,linestyle=dashed,dash=3pt 2pt]{->}(1,2.75)(1,3.25)
%\rput[l](3.5,3){\small$h_{T,a}$}
%\psline[linewidth=.5pt,linestyle=dashed,dash=3pt 2pt]{->}(3.4,2.75)(3.4,3.25)
%\rput[r](5.4,3){\small$h_{c,i}$}
%\psline[linewidth=.5pt,linestyle=dashed,dash=3pt 2pt]{->}(5.5,2.75)(5.5,3.25)
%
%\rput[tr](2.5,.5){\small$h_{b,t}$}
%\psline[linewidth=.5pt,linestyle=dashed,dash=3pt 2pt]{->}(1.85,.85)(3,.3)
%\rput[tl](4.1,.5){\small$h_{c,t}$}
%\psline[linewidth=.5pt,linestyle=dashed,dash=3pt 2pt]{->}(4.65,.85)(3.5,.3)
%
%\rput[t](2.2,1.45){\small$h_{T,b}$}
%\psline[linewidth=.5pt,linestyle=dashed,dash=3pt 2pt]{->}(2.65,1.6)(1.65,1.6)
%\rput[t](4.35,1.45){\small$h_{T,c}$}
%\psline[linewidth=.5pt,linestyle=dashed,dash=3pt 2pt]{->}(3.85,1.6)(4.85,1.6)
%
\endpspicture

%% file: fig6_X1.tex
\pspicture*(4,3.2)
\pscircle(2,1.6){1.5}
\pscircle*[linecolor=lightgray](2,1.5){.3}
\pscircle(2,1.5){.3}
\pscircle(1.1,1.5){.25}\rput(1.1,1.5){\small$c_1$}
\pscircle(1.6,2.4){.25}\rput(1.6,2.4){\small$c_2$}
\pscircle(2.4,2.4){.25}\rput(2.4,2.4){\small$c_3$}
\pscircle(2.9,1.5){.25}\rput(2.9,1.5){\small$c_4$}
\pscircle(1.6,.75){.25}\rput(1.6,.75){\small$c_5$}
\rput[l](2.1,.75){\small$b$}
\psarc[linewidth=.5pt](2,1.6){1.3}{90}{270}
\psline[linewidth=.5pt](2,2.9)(2,1.8)
\psline[linewidth=.5pt](2,1.2)(2,.3)
\psline[linewidth=.5pt,linestyle=dashed,dash=3pt 2pt](2,1.2)(2,1.8)
\rput[l](2.6,1.9){\small$a$}
\psline[linewidth=.5pt,linestyle=dashed,dash=3pt 2pt](1.7,1.5)(2.3,1.5)
\psline[linewidth=.5pt,linearc=.2](1.7,1.5)(1.5,1.5)(1.5,2)(2.5,2)(2.5,1.5)(2.3,1.5)
\endpspicture

%% file: fig7_X1.tex
\pspicture*(4,3.2)
\pscircle(2,1.6){1.5}
\pscircle*[linecolor=lightgray](2,1.5){.3}
\pscircle(2,1.5){.3}
\pscircle(1.1,1.5){.25}\rput(1.1,1.5){\small$c_1$}
\pscircle(1.6,2.4){.25}\rput(1.6,2.4){\small$c_2$}
\pscircle(2.4,2.4){.25}\rput(2.4,2.4){\small$c_3$}
\pscircle(2.9,1.5){.25}\rput(2.9,1.5){\small$c_4$}
\pscircle(1.6,.75){.25}\rput(1.6,.75){\small$c_5$}
\rput[tr](2.7,1){\small$d$}
\psarc[linewidth=.5pt](1.1,1.5){.45}{90}{270}
\psline[linewidth=.5pt](1.1,1.95)(1.8,1.7)
\psline[linewidth=.5pt](1.1,1.05)(1.8,1.3)
\psline[linewidth=.5pt,linestyle=dashed,dash=3pt 2pt](1.8,1.7)(2.2,1.3)
\psline[linewidth=.5pt,linestyle=dashed,dash=3pt 2pt](1.8,1.3)(2.2,1.7)
\psline[linewidth=.5pt](2.2,1.7)(2.9,1.95)
\psline[linewidth=.5pt](2.2,1.3)(2.9,1.05)
\psarc[linewidth=.5pt](2.9,1.5){.45}{270}{450}
\rput[r](1.1,2.3){\small$c$}
\psframe[linewidth=.5pt,framearc=.9](1.2,2.75)(2.8,2.05)
\endpspicture

%% file: sporadic.tex
\section{\label{sporadic}The sporadic surfaces}
Suppose that $F$ is not sporadic. To obtain a finite presentation
of the group $\M(F)$ using Theorem \ref{reduced} we need finite
presentations for the groups $\textrm{Stab}(s(v))$ and finite sets
of generators of the groups $\textrm{Stab}(s(e))$. By Proposition
\ref{sequence} we can reduce these problems to analogous problems
for the groups $\M(N)$, where $N$ is a connected component of
$F_{s(v)}$ or $F_{s(e)}$. Note that $N$ has either lower genus
than $F$ or equal genus, but less boundary components. If $N$ is
orientable then a finite presentation of $\M(N)$ is known (see
\cite{G} for the most general case). If $N$ is non-orientable and
not sporadic then we can obtain such presentation from Theorem
\ref{reduced}. Thus applying recursively Theorem \ref{reduced} we
obtain a finite presentation for $\M(F)$, provided that we know a
finite presentation of the mapping class group of each sporadic
subsurface.

The groups $\M(F_1^0)$ and $\M(F_1^1)$ are well known to be
trivial (cf. \cite{E}); $\M(F_1^2)$ is generated by Dehn twists
along the boundary curves and is isomorphic to $\Z^2$;
$\M(F_2^0)=\Z_2\times\Z_2$ (\cite{L1}). Simple presentation for
$\M(F_2^1)$ was found in \cite{S}, and for $\M(F_3^0)$ in
\cite{BC}. In this section we determine a finite presentation of
$\M(F_g^n)$ for the remaining sporadic surfaces, i.e. for
$(g,n)\in\{(1,3),(1,4),(2,2),(2,3),(3,1),(3,2)\}.$

\medskip

We begin by introducing the pure mapping class group of a
punctured surface and Birman's exact sequence, which is our main
tool in this section. Let $S$ be an {\it orientable} surface with
$2r$ distinguished points $\Sigma=\{q_1,\dots,q_{2r}\}$ called
{\it punctures}. The {\it pure mapping class group}
$\PM(S,\Sigma)$ is the group of isotopy classes \textsf{rel}
$\Sigma$ of all those diffeomorphisms of $S$ which fix each $q_i$.
Up to isomorphism, this group does not depend on the choice of
$\Sigma$, only on the number of punctures. We also define
$\PM(S,\emptyset)$ to be the ordinary mapping class group $\M(S)$.
Forgetting that $q_{2r-1}$ and $q_{2r}$ are
distinguished defines a homomorphism
$\rho\colon\PM(S,\Sigma)\to\PM(S,\Sigma')$, where
$\Sigma'=\Sigma\backslash\{q_{2r-1}, q_{2r}\}$.
Let $Q=\{(x_1,x_2)\in(S\backslash\Sigma')^2\ |\ x_1\ne x_2\}$.
We define {\it the pure braid group} $PB_2(S\backslash\Sigma')$ as $\pi_1(Q,(q_{2r-1},q_{2r}))$.
If the Euler characteristic of $S\backslash\Sigma'$ is negative, then
there is a short exact sequence due to Birman (see \cite{Bir}):
\begin{equation*}\label{braid}
1\to PB_2(S\backslash\Sigma')\stackrel{j}{\to}\PM(S,\Sigma)\stackrel{\rho}{\to}\PM(S,\Sigma')\to 1,
\end{equation*}
where the homomorphism $j$ is defined as follows. A loop $\beta\in PB_2(S\backslash\Sigma')$ defines an isotopy
of $0$-dimensional submanifold $(q_{2r-1},q_{2r})\subset S\backslash\Sigma'$, which can be extended to an isotopy
$h_t\in\Dif(S,\Sigma')$, $0\le t\le 1$ such that $h_0=1$ and $h_1(q_i)=q_i$ for $1\le i\le 2r$. We define $j(\beta)$ to be the isotopy class in $\Dif(S,\Sigma)$ of $h_1$.

Suppose that $\tau\colon S\to S$ is an orientation reversing
involution of $S$, without fixed points, and such that
$\tau(q_{2k-1})=q_{2k}$ for $1\le k\le r$. Then $S/\tau$ is a
non-orientable surface with $r$ distinguished points
$\Gamma=\{p_1,\dots,p_r\}$. Consider the subgroup
$\PM(S,\Sigma,\tau)$ of $\PM(S,\Sigma)$ consisting of all isotopy
classes which admit a representative which commutes with $\tau$.
It can be shown that two such representatives are isotopic
\textsf{rel} $\Sigma$ if and only if they are isotopic via an
isotopy which commutes with $\tau$ at each time (cf. \cite{BC}).
Since every diffeomorphism of $S/\tau$ has a unique orientation
preserving lift to $S$ which commutes with $\tau$ (the two lifts
differ by $\tau$ which is orientation reversing),
$\PM(S,\Sigma,\tau)$ can be identified with the group of isotopy
classes \textsf{rel} $\Gamma$ of diffeomorphisms of $S/\tau$ which
fix each $p_i$ and preserve the local orientation of $S/\tau$ at
each $p_i$.

It follows from the definition of $j$, that $j(\beta)\in\PM(S,\Sigma,\tau)$ if and only if
$\beta$ is represented by a loop of the form $t\mapsto(a_t,\tau(a_t))$, where $t\mapsto a_t$ is a loop in
$S\backslash\Sigma'$ based at $q_{2r-1}$.
Thus the pre-image $j^{-1}(\PM(S,\Sigma,\tau))$ can be identified with
$\pi_1(S\backslash\Sigma',q_{2r-1})$
 and we obtain the exact sequence:
\begin{equation}\label{braid1}
1\to\pi_1(S\backslash\Sigma')\stackrel{j}{\to}\PM(S,\Sigma,\tau)\stackrel{\rho}{\to}\PM(S,\Sigma',\tau)\to 1.
\end{equation}

\medskip

Suppose now that $F$ is a non-orientable surface of genus $g$ with
$r$ punctures $\Gamma=\{p_1,\dots,p_r\}$. Let
$\mathcal{PM}(F,\Gamma)$ denote the {\it pure mapping class group}
of $F$. It is defined as the group of the isotopy classes
\textsf{rel} $\Gamma$ of all diffeomorphisms of $F$ which fix each
$p_i$. Consider the subgroup $\mathcal{PM}^+(F,\Gamma)$ of
$\mathcal{PM}(F,\Gamma)$, consisting of the isotopy classes of those
diffeomorphisms which preserve the local orientation of $F$ at each
$p_i$. If $S$ is the orientable double cover of $F$ and $F=S/\tau$,
then it follows from above considerations that $\PM^{+}(F,\Gamma)$
can be identified with $\PM(S,\Sigma,\tau)$. Note that
$\pi_1(S\backslash\Sigma')$ can be identified with the subgroup
$\pi^{+}_1(F\backslash\Gamma',p_r)$ of
$\pi_1(F\backslash\Gamma',p_r)$ consisting of the two-sided loops.
With such identifications the sequence (\ref{braid1}) becomes:
\begin{equation}\label{braid3}
1\to\pi^{+}_1(F\backslash\Gamma',p_r)\stackrel{j}{\to}\PM^{+}(F,\Gamma)\stackrel{\rho}{\to}\PM^{+}(F,\Gamma')\to 1,
\end{equation}
where we assume that the Euler characteristic of $F\backslash\Gamma'$ is negative (that is $g+r>3$).

In this paper we use the same symbol to denote a loop and its homotopy class in the
fundamental group. In order for $j$ to be a homomorphism, the product $\alpha\beta$ of two
loops should mean first travel along $\beta$ and then along $\alpha$.

If $\alpha$ is a simple loop in $F$ based at $p_r$, then $j(\alpha)$ is the isotopy class of a diffeomorphism
obtained by sliding $p_r$ once along $\alpha$.

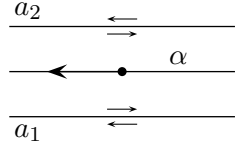
\begin{figure}
\begin{center}
\input{fig28_sl}
\caption{\label{sl}$j(\alpha)=t_{a_1}t_{a_2}$}
\end{center}
\end{figure}

The next two lemmas are proved in \cite{I1}, (6.1).

\begin{lemma}\label{j1}
Let $\alpha\in\pi^{+}_1(F\backslash\Gamma',p_k)$ be a two-sided
simple loop and let $a_1$, $a_2$ denote boundary curves of a
tubular neighborhood of $\alpha$. Then $j(\alpha)=t_{a_1}t_{a_2}$,
where $t_{a_1}$ and $t_{a_2}$ are Dehn twists about $a_1$ and
$a_2$ in the directions indicated by arrows in Figure \ref{sl}.
\foorp
\end{lemma}

The pure mapping class group $\PM(F,\Gamma)$ acts on $\pi^+_1(F\backslash\Gamma')$ in the obvious
way. We denote this action by $h(\alpha)$ for
$h\in\PM(F,\Gamma)$ and $\alpha\in\pi^+_1(F\backslash\Gamma')$.

\begin{lemma}\label{j2}
The homomorphism $j$ is $\PM(F,\Gamma)$-equiveriant. That is  $j(h(\alpha))=hj(\alpha)h^{-1}$ for
$h\in\PM(F,\Gamma)$ and $\alpha\in\pi^+_1(F\backslash\Gamma')$.
\foorp
\end{lemma}

\medskip

Suppose that $\widetilde{F}=F_g^n$ is a non-orientable surface of
negative Euler characteristic (i.e. $g+n>2$) and let
$c_1,\dots,c_n\colon S^1\to\bdr\widetilde{F}$ denote the boundary
curves. Let $F=F_g^0$ be the closed surface with punctures
$\Gamma=\{p_1,\dots,p_n\}$ obtained by gluing a disc with a
puncture $p_i$ to $\bdr\widetilde{F}$ along  $c_i$ for  $1\le i\le
n$. We identify $\widetilde{F}$ with a subsurface of $F$ and
denote by $i_\ast\colon\M(\widetilde{F})\to\PM^+(F,\Gamma)$ the
homomorphism induced by the inclusion $i\colon\widetilde{F}\to F$.
It can be proved, using the same methods as in the proof of
Proposition \ref{kernel}, that $\ker i_\ast$ is a free abelian
group of rank $n$ generated by Dehn twists about the boundary
curves $c_i$. Thus we have the exact sequence
\begin{equation}\label{Mgn}
1\to\Z^n\to\M(F_g^n)\stackrel{i_\ast}{\to}\mathcal{PM}^+(F_g^0,\Gamma)\to
1.
\end{equation}

\begin{rem}\label{centr}
Note that $\ker i_\ast$ is a central subgroup of $\M(F_g^n)$.
Indeed, for every $i\in\{1,\dots,n\}$ and $h\in\M(F_g^n)$ we have
$h t_{c_i}h^{-1}=t_{h(c_i)}=t_{c_i}$.
\end{rem}

%%%%%%%%%%%%%%%%%%%%%%%%%%%%%%%%%%%%%%%%%%%%%%%%%%%%%%%%%%%%%%%%%%%%%%%%%%%%%%%%%%%%%%%%%%%%%%
%%%%%%%%%%%%%%%%%%%%%%%%%%%%%%%%%%%%%%%%%%%%%%%%%%%%%%%%%%%%%%%%%%%%%%%%%%%%%%%%%%%%%%%%%%%%%%%
%%%%%%%%%%%%%%%%%%%%%%%%%%%%%%%%%%%%%%%%%%%%%%%%%%%%%%%%%%%%%%%%%%%%%%%%%%%%%%%%%%%%%%%%%%%%%%
\medskip

We record without proof the following easy lemma.

\begin{lemma}\label{present}
Consider a short exact sequence of groups
$$1\to K\stackrel{i}{\to}G\stackrel{p}{\to}H\to 1$$
and suppose that $K$ and $H$ admit presentations
$$K=\langle G_K\,|\,R_K\rangle,\quad H=\langle G_H\,|\,R_H\rangle.$$
Then $G$ admits the presentation
\begin{equation}\label{pres}
\langle\, i(G_K)\cup\widetilde{G_H}\ |\ i(R_K)\cup\widetilde{R_H}\cup R\,\rangle,
\end{equation}
where:

$i(G_K)=\{i(k)\,|\,k\in G_K\}$, $\widetilde{G_H}=\{\tilde{h}\,|\,h\in G_H\}$, where $\tilde{h}$ is any element in $G$ such that
 $p(\tilde{h})=h$,

$i(R_K)=\{i(k_1)\cdots i(k_n)\,|\,k_1\cdots k_n\in R_K\}$,

$\widetilde{R_H}=\{\tilde{h}_1\cdots\tilde{h}_nw(h_1\cdots h_n)\,|\,h_1\cdots h_n\in R_H\}$,

$R=\{\tilde{h}i(k)\tilde{h}^{-1}w(k,h)\,|\,h\in G_H, k\in G_k\},$\\
where $w(h_1\cdots h_n)$ and $w(k,h)$ are suitable words in generators $i(G_K)$.
\end{lemma}

We can now obtain finite presentations for the mapping class groups
$\M(F_g^n)$ of the sporadic surfaces in the following way. Starting
from known presentations of the groups $\PM^+(F_1^0,\{p_1,p_2\})$,
$\PM^+(F_2^0,\{p_1\})$ and $\M(F_3^0)$, we obtain presentations for
all $\PM^+(F_g^0,\Gamma)$, by applying recursively Lemma
\ref{present} to the sequence (\ref{braid3}). To do this, we need
finite presentations for the groups
$\pi_1^+(F_g^0\backslash\Gamma')$. These can be obtained from
standard presentations of fundamental groups
$\pi_1(F_g^0\backslash\Gamma')$ by the Reidemeister-Schreier method
(see, for example, \cite{MKS}). Once we have found the presentations
for $\PM^+(F_g^0,\Gamma)$, we obtain presentations for $\M(F_g^n)$,
by applying Lemma \ref{present} to the sequence (\ref{Mgn}).

\subsection{Sporadic surfaces of genus 1.}

Until the end of this paper we use the capital letter $A$ to
denote a Dehn twist about the curve labelled as $a$. In order for
this notation to be unambiguous, we have to specify the direction
of the twist $A$ for each curve $a$. Equivalently we may choose an
orientation of a tubular neighborhood of $a$. Then $A$ denotes the
right Dehn twist with respect to the chosen orientation.

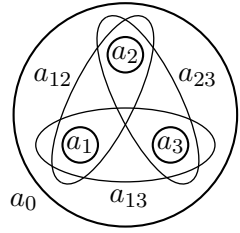
\begin{figure}
\begin{center}
\input{fig8_lantern}
\caption{\label{lantern}The curves of the lantern relation.}
\end{center}
\end{figure}

\medskip

Consider a 2-sphere $S$ with four holes embedded in $F$. Let
$a_0$, $a_1$, $a_2$, $a_3$ denote disjoint boundary curves of $S$,
and $a_{12}$, $a_{13}$, $a_{23}$ separating generic curves such
that $a_{ij}$ separates $a_i$ and $a_j$ from the other two
boundary curves of $S$ (Figure \ref{lantern}). If $A_i$ and
$A_{jk}$ are right Dehn twists with respect to the standard
orientation of the plane of Figure \ref{lantern}, then we have the
well known lantern relation:
\begin{equation}\label{L}
A_0A_1A_2A_3=A_{12}A_{13}A_{23}.
\end{equation}
The lantern relation was discovered by  Dehn \cite{Deh} and rediscovered by Johnson \cite{Joh}.
Note that since $A_{ij}$ commutes with $A_k$, we have:
\begin{equation}\label{L1}
A_{12}A_{13}A_{23}=A_{13}A_{23}A_{12}=A_{23}A_{12}A_{13}.
\end{equation}

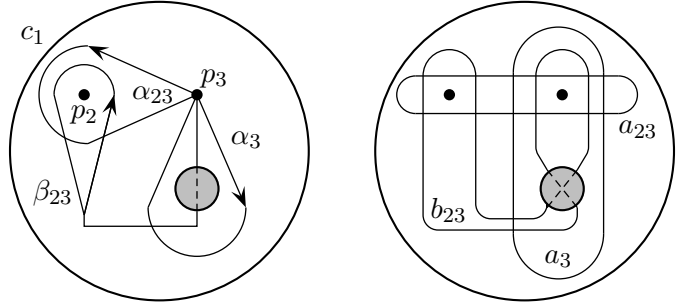
\begin{figure}%[!htbp]
\begin{center}
\begin{tabular}{cc} \input{fig9_Pip13} & \input{fig10_M13}%\\
%a) & b)
\end{tabular}
\caption{\label{F13} Generators of $\pi^+_1(F\backslash\{p_1,p_2\},p_3)$ and
generic curves in $F_1^3$.}
\end{center}
\end{figure}
\begin{figure}
\begin{center}
\begin{tabular}{cc}
\input{fig11_Pip14a} & \input{fig12_Pip14b}
\end{tabular}
\caption{\label{F14} Generators of $\pi^+_1(F\backslash\{p_1,p_2,p_3\},p_4)$.}
\end{center}
\end{figure}
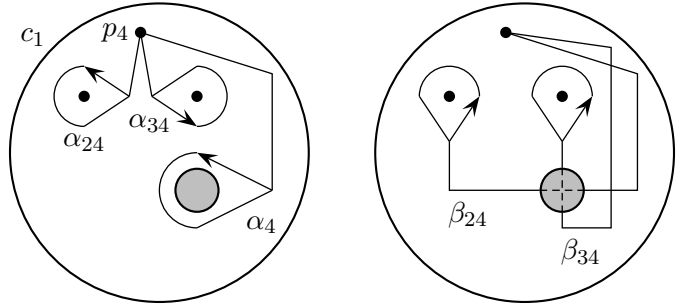

\medskip

Let us fix four points $p_1,\dots,p_4$ in the projective plane
$F=F_1^0$ represented in Figures \ref{F13} and \ref{F14}, where the curve $c_1$
bounds in $F$ a disc containing $p_1$. Let $n\in\{3,4\}$ and
consider the embedding $i\colon\widetilde{F}\to F$, where
$\widetilde{F}=F_1^n$, and the induced homomorphism
$i_\ast\colon\M(\widetilde{F})\to\PM^+(F,\{p_1,\dots,p_n\})$ (if
$n=3$ then we forget that $p_4$ is distinguished). We identify
$\widetilde{F}$ with $i(\widetilde{F})$, and a curve $a$ in
$\widetilde{F}$ with $i\circ a$ in $F$.

Consider the loops $\alpha_i$, $\alpha_{jk}$, $\beta_{jk}$
represented in Figures \ref{F13} and \ref{F14}, where we assume,
that each of them represents a two-sided simple loop in
$\pi^+_1(F\backslash\{p_1,p_2\},p_3)$ or
$\pi^+_1(F\backslash\{p_1,p_2,p_3\},p_4)$. The boundary of a
tubular neighborhood of such loop consist of two two-sided simple
closed curves, one of which is trivial (i.e. it either separates a
M\"obius strip or a disc containing one puncture). We use the
symbol $a_i$ or $a_{jk}$ or $b_{jk}$ to denote the non-trivial
boundary component of the tubular neighborhood of the
corresponding loop (see Figure \ref{F13}). Then by Lemma \ref{j1},
we have $j(\alpha_i)=A_i$, $j(\alpha_{jk})=A_{jk}$,
$j(\beta_{jk})=B_{jk}$. Note that $a_i$, $a_{jk}$, $b_{jk}$ may be
chosen to be generic curves in $\widetilde{F}$.

\begin{theorem}\label{13}
The group $\PM^+(F,\{p_1,p_2,p_3\})$ is free, generated by
$A_3$, $A_{23}$, $B_{23}$. The group $\M(F_1^3)$ is generated by
$A_3$, $A_{23}$, $B_{23}$, $C_1$, $C_2$, $C_3$ and isomorphic to
$\Z^3\times\PM^+(F,\{p_1,p_2,p_3\}).$
\end{theorem}

\proof It can be deduced from Theorem 4.1 of \cite{K} that the group
$\PM^+(F,\{p_1,p_2\})$ is trivial. Thus
$$j\colon\pi^+_1(F\backslash\{p_1,p_2\},p_3)\to\PM^+(F,\{p_1,p_2,p_3\})$$
is an isomorphism.
The fundamental group $\pi_1(F\backslash\{p_1,p_2\},p_3)$ is free on generators
$\alpha_{23}$ and $x$, where $x$ is a one-sided loop, such that $x^2=\alpha_3^{-1}$,
$x\alpha_{23}x^{-1}=\beta_{23}$. Now $\{1,x\}$ is a Schreier system of representatives
of right cosets of $\pi^+_1(F\backslash\{p_1,p_2\},p_3)$ and by the
Reidemeister-Schreier method we obtain that the last group is freely generated
by the loops $\alpha_3$, $\alpha_{23}$, $\beta_{23}$.
Hence the first part of Theorem \ref{13}. The second part follows from
the sequence (\ref{Mgn}). Indeed, the sequence splits as
$\pi^+_1(F\backslash\{p_1,p_2\},p_3)$ is free, and the kernel of $i_\ast$ is central by
Remark \ref{centr}.\foorp

\begin{theorem}\label{14}
The group $\PM^+(F,\{p_1,p_2,p_3,p_4\})$ admits a presentation with
generators
 $\{A_3, A_4, A_{23}, A_{24}, A_{34}, B_{23},
B_{24}, B_{34}, D\}$ and relations:

\noindent
$\qquad(1)\ A_{23}A_4=A_4A_{23},\ A_{24}A_3=A_3A_{24},$

\noindent
$\qquad(2)\ A_3^{-1}A_4A_{34}B_{34}=B_{34}A_3^{-1}A_4A_{34},$

\noindent
$\qquad(3)\ A_4A_{34}A_{24}B_{23}=B_{23}A_4A_{34}A_{24},$

\noindent
$\qquad(4)\ A_{34}A_3^{-1}A_{23}B_{24}=B_{24}A_{34}A_3^{-1}A_{23},$

\noindent
$\qquad(5)\ A_{34}A_{24}A_{23}=A_{24}A_{23}A_{34}=A_{23}A_{34}A_{24},$

\noindent
$\qquad(6)\ B_{34}A_{23}B_{24}=A_{23}B_{24}B_{34}=B_{24}B_{34}A_{23},$

\noindent
$\qquad(7)\ A_{4}A_{34}A_3^{-1}=A_{34}A_3^{-1}A_{4}=A_3^{-1}A_{4}A_{34},$

\noindent
$\qquad(8)\ A^{-1}_{34}B_{24}B_{23}=B_{24}B_{23}A^{-1}_{34}=B_{23}A^{-1}_{34}B_{24},$

\noindent
$\qquad(9)\ A_{24}B_{23}D^{-1}=B_{23}D^{-1}A_{24}=D^{-1}A_{24}B_{23},$

\noindent
$\qquad(10)\ D=A^{-1}_{34}A^{-1}_4B_{34}A_4A_{34}.$

\noindent
The group $\M(F_1^4)$ is isomorphic to
$\Z^4\times\PM^+(F,\{p_1,p_2,p_3,p_4\}).$
\end{theorem}

\proof Let us denote, for simplicity,
$$\pi=\pi^+_1(F\backslash\{p_1,p_2,p_3\},p_4),\quad
G=\PM^+(F,\{p_1,p_2,p_3,p_4\}).$$
The fundamental group $\pi_1(F\backslash\{p_1,p_2,p_3\},p_4)$ is free on generators
$\alpha_{24}$, $\alpha_{34}$ and $x$, where $x$ is a one-sided loop, such that $x^2=\alpha_4$,
$x\alpha_{24}x^{-1}=\beta_{24}$, $x\alpha_{34}x^{-1}=\beta_{34}$. Now $\{1,x\}$ is a Schreier system of representatives
of right cosets of $\pi$ and by the
Reidemeister-Schreier method we obtain that $\pi$ is freely generated
by the loops in Figure \ref{F14}.
By Lemma \ref{pres} applied to sequence (\ref{braid3}) and Theorem
\ref{13}, $G$ admits a presentation with generators
$A_3$, $A_{23}$, $B_{23}$, $A_4=j(\alpha_4)$, $A_{k4}=j(\alpha_{k4})$,
$B_{k4}=j(\beta_{k4})$, $k=2,3$ and relations
$hgh^{-1}\in j(\pi)$ for each $h\in\{A_3, A_{23}, B_{23}\}$,
$g\in\{A_4, A_{k4}, B_{k4}\,|\,k=2,3\}$. We will show that all these relations are consequences of (1-10).
We have:

\noindent
$(1)\Rightarrow A_{23}A_4A_{23}^{-1}, A_3A_{24}A_3^{-1}\in j(\pi)$;
$(2)\Rightarrow A_3B_{34}A_3^{-1}\in j(\pi)$;\\ $(10)\Rightarrow D\in j(\pi)$.
From (5) follows $A_{23}A_{34}A_{23}^{-1}=A_{24}^{-1}A_{34}A_{24}\in j(\pi)$,\\
$A_{23}A_{34}A_{24}A^{-1}_{23}=A_{34}A_{24}\Rightarrow
A_{23}A_{24}A^{-1}_{23}\in j(\pi)$. Analogously we have
$(6-9)\Rightarrow\{A_{23}B_{24}A_{23}^{-1}, A_{23}B_{34}A_{23}^{-1},
A_3A_{34}A_3^{-1}, A_3A_4A_3^{-1}, B_{23}A_{34}B_{23}^{-1},\\
B_{23}B_{24}B_{23}^{-1}, B_{23}DB_{23}^{-1},
B_{23}A_{24}B_{23}^{-1}\}\subset j(\pi).$ From (3) follows\\
$B_{23}A_4A_{34}B_{23}^{-1}\in j(\pi)$; from this and (8) we have
$B_{23}A_4B_{23}^{-1}\in j(\pi)$ and from (10) follows
$B_{23}B_{34}B_{23}^{-1}\in j(\pi)$. Finally we have $(4)\Rightarrow
A_3A_{34}^{-1}B_{24}A_{34}A_3^{-1}=A_{23}B_{24}A_{23}^{-1}$, and by
(6,7) we have $A_3B_{24}A_3^{-1}\in j(\pi)$.

Now we show that relations (1-10) are satisfied in
$\M(\widetilde{F})$, and hance also in $G$. By relation (10), the
generator $D$ is a Dehn twist about the curve
$A^{-1}_{34}A^{-1}_4(b_{34})$ bounding a pair of pants together
with $c_3$ and $c_4$. The relations (1) are obvious. By
considering appropriate embeddings of a 2-sphere with four holes
in $\widetilde{F}$, it is easy to recognize (5-9) as relations of
type (\ref{L1}), i.e. consequences of the lantern relation. In
particular, we have lantern relation
$A_{12}C_3C_4=A_4A_{34}A_3^{-1}$, where $A_{12}$ is Dehn twist
about a curve bounding a pair of pants together with $c_1$ and
$c_2$. Since $B_{34}$ commutes with $A_{12}$, $C_3$ and $C_4$, the
relation (2) holds. By Lemma \ref{j1}, we have
$j(\alpha_4\alpha_{34}\alpha_{24})=A_{14}\in G$, where $a_{14}$
bounds a pair of pants in $\widetilde{F}$ together with $c_1$ and
$c_4$. Thus in $\M(\widetilde{F})$ we have
$A_4A_{34}A_{24}=A_{14}C$, where $C$ is a product of twists
$C_1,\dots,C_4$. Since $B_{23}$ commutes with $A_{14}$ and $C$,
(3) holds. Consider a monomorphism
$j'\colon\pi^+_1(F\backslash\{p_1,p_2,p_4\},p_3)\to G$, defined
like $j$. There exists exactly one loop $\alpha'_{34}$ such that
$j'(\alpha'_{34})=A_{34}\in G$, and we have
$j'(\alpha'_{34}\alpha_3^{-1}\alpha_{23})=A_{13}\in G$, where
$a_{13}$ bounds a pair of pants in $\widetilde{F}$ together with
$c_1$ and $c_3$. Since $B_{24}$ commutes with $A_{13}$, (4) holds.

We have shown that (1-10) are relations in $G$, and all relations
from Lemma \ref{pres} are consequences of (1-10). Hence $G$ admits
presentation with relations (1-10). Since these relations hold
also in $\M(\widetilde{F})$, the sequence (\ref{Mgn}) splits, and
since the kernel of $i_\ast$ is central, we obtain
$\M(\widetilde{F})=\Z^4\times G$.\foorp

\subsection{\label{s2}Sporadic surfaces of genus 2.}

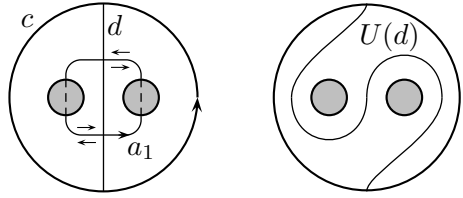
\begin{figure}%[!htbp]
\begin{center}
\input{fig13_U}
\caption{\label{U}The diffeomorphism $U$.}
\end{center}
\end{figure}

\medskip

Consider the Klein bottle $K$ with one hole represented in Figure
\ref{U}. Let $U$ be a diffeomorphism of $K$ interchanging the
shaded discs in Figure \ref{U} and such that $U^2$ is the Dehn
twist about the boundary curve $c$, right with respect to the
standard orientation of the plane of the figure. Up to isotopy,
$U$ acts on the arc $d$ as it is shown in Figure \ref{U} (see \cite{Sz}
for precise definition). We fix Dehn twist $A_1$ about the curve
$a_1$, in the direction indicated by arrows in Figure \ref{U}. The
composition $UA_1$ is the Y-homeomorphism (or cross-cap slide)
introduced by Lickorish \cite{L1}. The next theorem follows
immediately from Theorem A.7 of \cite{S}.

\begin{theorem}\label{21}
The mapping class group $\M(K)$ is generated by $A_1$ and $U$ and admits the presentation
$\langle A_1,U\,|\,UA_1U^{-1}=A_1^{-1}\rangle.$
\foorp
\end{theorem}

\begin{figure}
\begin{center}
\begin{tabular}{ccc}
\input{fig14_lanK1}&
\input{fig15_lanK2}&
\input{fig16_lanK3}
\end{tabular}
\caption{\label{lanK} The surfaces $\widetilde{F}=F_2^2$ and
$\widetilde{F}_{a_1}$.}
\end{center}
\end{figure}
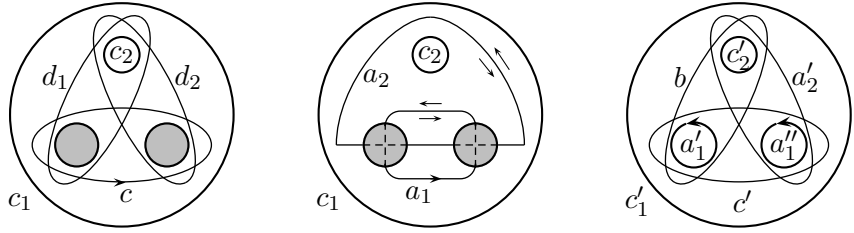

Let $\widetilde{F}=F_2^2$ be the surface obtained by gluing a pair
of pants to $K$, and let $c_1$ and $c_2$ denote the boundary
curves of $\widetilde{F}$ (Figure \ref{lanK}). We extend $U$ by
the identity outside $K$ to a diffeomorphism of $\widetilde{F}$.
Let $C$, $C_1$, $C_2$ $D_1$, $D_2$ be Dehn twists about the curves
represented in Figure \ref{lanK}, right with respect to the
standard orientation of the plane of the figure. We also define
Dehn twist $A_1$, $A_2$ in the indicated directions. Note that
$U^2=C$ and $UD_2U^{-1}=D_1$.

The right hand side of Figure \ref{lanK} represents the four-holed sphere
$\widetilde{F}_{a_1}$ obtained by cutting $\widetilde{F}$ along $a_1$, where
$\rho_{a_1}\circ a'_1=\rho_{a_1}\circ a''_1=a_1$, $\rho_{a_1}(c'_i)=c_i$ for
$i=1,2$, $\rho_{a_1}(c')=c$, $\rho_{a_1}(a'_2)=a_2$,
$\rho_{a_1}(b)=U(a_2)$. If $C'_i$, $C'$,
 $A'_1$, $A''_1$, $A'_2$, $B$ are right Dehn twists with respect to  the standard orientation of the plane of Figure
\ref{lanK}, then
$\rho_\ast(C'_i)=C_i$, $\rho_\ast(C')=C$, $\rho_\ast(A'_1A''_1)=1$, $\rho_\ast(A'_2)=A_2$, and $\rho_\ast(B)=UA_2U^{-1}$.

\begin{lemma}\label{K2rel}
In $\M(\widetilde{F})$ we have $(A_2U)^2=(D_2U)^2=C_1C_2.$
\end{lemma}
\proof
We have the lantern relation $C'_1C'_2A'_1A''_1=A'_2BC'.$
By applying $\rho_\ast$ to both sides we obtain $C_1C_2=A_2(UA_2U^{-1})U^2=(A_2U)^2.$
By another lantern relation we have $C_1C_2=D_2D_1C=D_2(UD_2U^{-1})U^2=(D_2U)^2.$\foorp

\begin{figure}%[!htbp]
\begin{center}
\begin{tabular}{cc}
 \input{fig17_Pi22} & \input{fig18_Pip22}
%a) & b)
\end{tabular}
\caption{\label{Pi22}Generators of $\pi_1(F\backslash\{p_1\},p_2)$ and
$\pi^+_1(F\backslash\{p_1\},p_2)$.}
\end{center}
\end{figure}
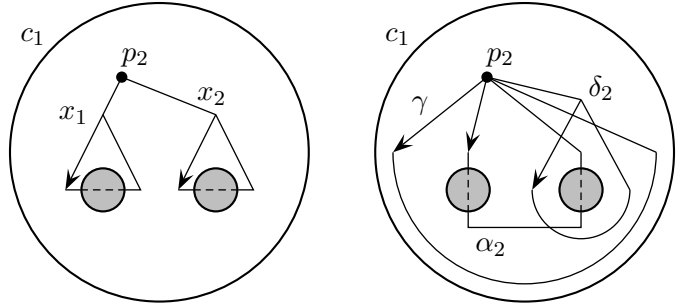

\medskip

Let $F=F_2^0$ be the Klein bottle obtained by gluing a disc with a
puncture $p_i$ to $\bdr\widetilde{F}$ along $c_i$ for $i=1,2$. We
identify $U$, $A_1$, $A_2$, $D_2$, with $i_\ast(U)$,
$i_\ast(A_1)$, $i_\ast(A_2)$, $i_\ast(D_2)$ respectively, where
$i_\ast\colon\M(\widetilde{F})\to\PM^+(F,\{p_1,p_2\})$ is the
homomorphism induced by the inclusion of $\widetilde{F}$ in $F$.

\begin{theorem}\label{P22}
The group $\PM^+(F,\{p_1,p_2\})$ admits a presentation with generators
$\{A_1, A_2, D_2, U\}$
and relations:
$\quad A_1A_2=A_2A_1,\\ UA_1U^{-1}=A^{-1}_1,
\ A_2UD_2=D_2^{-1}A_2U,\ (A_2U)^2=(D_2U)^2=1.$
\end{theorem}
\proof
Consider the exact sequence (\ref{braid3}):
$$1\to\pi^+(F\backslash\{p_1\},p_2)\stackrel{j}{\to}\PM^+(F,\{p_1,p_2\})\to\PM^+(F,\{p_1\})\to 1.$$
By Theorem \ref{21} and sequence (\ref{Mgn}), $\PM^+(F,\{p_1\})$ has presentation
$$\langle A_1,U\,|\,UA_1U^{-1}=A_1^{-1}, U^2=1\rangle.$$
The fundamental group $\pi_1(F\backslash\{p_1\},p_2)$ is free on
generators $x_1$, $x_2$ in Figure \ref{Pi22}. Now $\{1,x_2\}$ is a
Schreier system of representatives of cosets of
$\pi^+_1(F\backslash\{p_1\},p_2)$ and by the Reidemeister-Schreier
method we obtain that the last group is freely generated by
$\delta_2=x^2_2$, $\alpha_2=x_2x_1$ and $x_1x_2^{-1}$. It follows
that $\pi^+_1(F\backslash\{p_1\},p_2)$ is free on generators
$\delta_2$, $\alpha_2$, $\gamma$, where
$\gamma=x_2^2(x_1x_2^{-1})(x_2x_1)$. Observe that
$j(\gamma)=U^{-2}$, $j(\alpha_2)=A_2A_1^{-1}$, $j(\delta_2)=D_2.$
By Lemma \ref{pres}, $\PM^+(F,\{p_1,p_2\})$ admits presentation
with generators $U$, $A_1$, $j(\gamma)$, $j(\alpha_2)$,
$j(\delta_2)$ and relations $UA_1U^{-1}=A_1^{-1}$,
$U^2=(j(\gamma))^{-1}$, and (by Lemma \ref{j2}):
$Uj(\gamma)U^{-1}=j(\gamma)$,
$Uj(\alpha_2)U^{-1}=j(\alpha^{-1}_2\gamma)$,
$Uj(\delta_2)U^{-1}=j(\delta^{-1}_2\gamma)$,
$A_1j(\gamma)A^{-1}_1=j(\gamma)$,
$A_1j(\alpha_2)A^{-1}_1=j(\alpha_2)$,
$A_1j(\delta_2)A^{-1}_1=j(\gamma\alpha^{-1}_2\delta_2\alpha_2)$.
Substituting $j(\gamma)=U^{-2}$, $j(\alpha_2)=A_2A_1^{-1}$,
$j(\delta_2)=D_2$ we obtain a presentation which can easily be
shown to be equivalent to that in Theorem \ref{P22}.\foorp

\begin{theorem}\label{22}
The group $\M(F_2^2)$ admits a presentation with generators $\{C_1, A_1, A_2, D_2, U\}$
and relations: $\quad C_1A_i=A_iC_1,\ \textrm{for $i=1,2$},\\
C_1D_2=D_2C_1,\ C_1U=UC_1,\
A_1A_2=A_2A_1,\ UA_1U^{-1}=A^{-1}_1,\\
A_2UD_2=D_2^{-1}A_2U,\ (A_2U)^2=(D_2U)^2$.
\foorp
\end{theorem}

\proof From sequence (\ref{Mgn}), Theorem \ref{P22} and Lemma \ref{K2rel} we obtain a presentation
for $\M(F_2^2)$ with generators $\{C_1, C_2, A_1, A_2, D_2, U\}$ and relations listed in Theorem \ref{22}
and
\begin{equation}\label{c2rel}
C_1C_2=C_2C_1,\ C_2D_2=D_2C_2,\ C_2U=UC_2,\ C_2A_i=A_iC_2,
\end{equation}
for $i=1,2$ and
\begin{equation}\label{c2def}
(A_2U)^2=C_1C_2.
\end{equation}
We claim that the relations (\ref{c2rel}) are consequences of the relation (\ref{c2def})
and relations from Theorem \ref{22}. Clearly it suffices to check that relations
$$D_2(A_2U)^2=(A_2U)^2D_2,\ U(A_2U)^2=(A_2U)^2U,\ A_i(A_2U)^2=(A_2U)^2A_i,$$
follow from those in Theorem \ref{22}. Observe that
$A_1(A_2U)^2=(A_2U)^2A_1$ follows from $A_1A_2=A_2A_1$ and
$UA_1U^{-1}=A_1^{-1}$. From $A_2UD_2=D_2^{-1}A_2U$ we have
$D_2^{-1}(A_2U)^2D_2=(A_2U)^2$ and
$U(A_2U)^2U^{-1}=\\U(D_2U)^2U^{-1}=D_2^{-1}(D_2U)^2D_2=D_2^{-1}(A_2U)^2D_2=(A_2U)^2.$
Finally we have $A_2^{-1}(A_2U)^2A_2=U(A_2U)^2U^{-1}=(A_2U)^2.$ It
follows that relations (\ref{c2rel}) are redundant, and hence they
can be removed from the presentation. Then the generator $C_2$ can
also be removed together with the relation (\ref{c2def}). \foorp

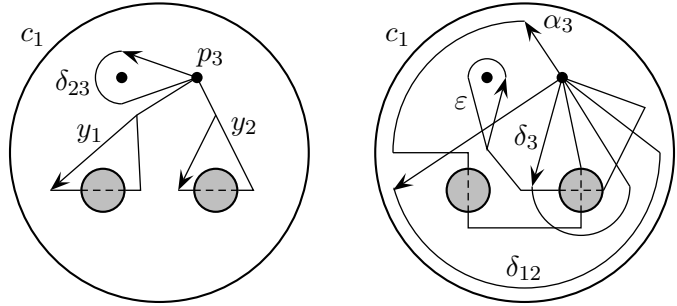
\begin{figure}%[!htbp]
\begin{center}
\begin{tabular}{cc}
\input{fig19_Pi23}&\input{fig20_Pip23}
%a) & b)
\end{tabular}
\caption{\label{Pi23}Generators of $\pi_1(F\backslash\{p_1,p_2\},p_3)$ and
$\pi^+_1(F\backslash\{p_1,p_2\},p_3)$.}
\end{center}
\end{figure}

\medskip

We fix a point $p_3\in F\backslash K$, different from $p_2$ and
$p_1$, and such that $p_3$ and $p_2$  are in different components
of $F\backslash (a_1\cup a_2)$. We identify $U$, $A_2$, $A_1$ and
$D_2$ with elements of $\PM^+(F,\{p_1,p_2,p_3\})$. Let $A_3$ and
$D_3$ be such Dehn twists that $j(\alpha_3)=A_3A_2^{-1}$ and
$j(\delta_3)=D_3$, where $\alpha_3$, $\delta_3$ are the loops in
Figure \ref{Pi23}, and
$j\colon\pi^+(F\backslash\{p_1,p_2\},p_3)\to\PM^+(F,\{p_1,p_2,p_3\})$
is the monomorphism from sequence (\ref{braid3}).
\begin{theorem}\label{P23}
The group $\PM^+(F,\{p_1,p_2,p_3\})$ admits a presentation with generators
$\{A_1, A_2, A_3, D_2, D_3, U\}$ and relations:
\begin{eqnarray*}
&&(1)\ A_iA_j=A_jA_i,\ \textrm{for $i,j\in\{1,2,3\}$};\quad (2)\ UA_1U^{-1}=A^{-1}_1;\\
&&(3)\ A_2UD_2=D_2^{-1}A_2U;\quad (4)\ (A_2U)^2=(D_2U)^2=(UD_2)^2;\\
&&(5)\ (UD_3)^2=(D_3U)^2;\quad (6)\  D_3UD_2U^{-1}=UD_2U^{-1}D_3;\\
&&(7)\ A_3UD_2D_3=UD_2D_3A_3^{-1};\quad (8)\ (UA_3)^2=(UD_2D_3)^{-2};\\
&&(9)\ A_2(A_3UD_2)^2=(A_3UD_2)^2A_2;\\
&&(10)\ A_2A_1^{-1}D_3A_1A_2^{-1}=A_3UD_2D_3^{-1}(A_3UD_2)^{-1};\\%U^{-1}A^{-1}_3D_3^{-1}A_3U;\\
&&(11)\ A_1(A_3UD_2)^2A_1^{-1}=(UD_2)^{-1}(A_3UD_2)^2UD_2.
\end{eqnarray*}
\end{theorem}

\proof
Let us denote, for simplicity,
$$\pi=\pi^+_1(F\backslash\{p_1,p_2\},p_3),\quad
G=\PM^+(F,\{p_1,p_2,p_3\}).$$ The fundamental group
$\pi_1(F\backslash\{p_1,p_2\},p_2)$ is free on generators
$\delta_{23}$, $y_1$, $y_2$ in Figure \ref{Pi23}. Now $\{1,y_2\}$
is a Schreier system of representatives of cosets of $\pi$ and by
the Reidemeister-Schreier method we obtain that the last group is
freely generated by $\delta_{23}$, $\delta_3=y^2_2$,
$\varepsilon=y_2\delta_{23}y^{-1}_2$, $y_2y_1$ and $y_1y_2^{-1}$.
It follows that $\pi$ is free on generators $\delta_{23}$,
$\delta_3$, $\varepsilon$, $\alpha_3$, $\delta_{12}$, where
$\delta_{12}=\delta_3(y_1y_2^{-1})(y_2y_1)$,
$\alpha_3=y_2y_1\delta_{23}$. By Lemmas \ref{j1} and \ref{K2rel}
we have
\begin{equation}\label{sub1}
j(\delta_{23})=(UA_3)^2,\quad j(\delta_{12}\delta_{23})=(UD_2)^{-2}.
\end{equation}
First we show that relations (1-11) are satisfied in $G$:
(1) and (6) are obvious; (4) and (5) follow from Lemma
\ref{K2rel}; (2), (3), (7) are relations of type $ht_ah^{-1}=t^{\pm 1}_{h(a)}$ and hence they can be checked by looking at the
effect of $h$ on the curve $a$;
(10) follows from $A_2A_1^{-1}(\delta_3)=A_3UD_2(\delta^{-1}_3)$;
(8) is equivalent to
$UD_2D_3D^{-1}_2U^{-1}=(UA_3)^{-2}D^{-1}_3(UD_2)^{-2},$
which follows from $UD_2(\delta_3)=\delta^{-1}_{23}\delta^{-1}_3\delta_{12}\delta_{23}$.
It can be checked that $\varepsilon\delta_3=A_3((\delta_{12}\delta_{23})^{-1}\delta_3)$ and hence
$j(\varepsilon)=A_3(UD_2)^2D_3A^{-1}_3D^{-1}_3$; from this and (7) we obtain
\begin{equation}\label{sub2}
j(\varepsilon)=(A_3UD_2)^2.
\end{equation}
Now (9) and (11) follow from (\ref{sub2}) and the equalities $A_2(\varepsilon)=\varepsilon$ and $A_1(\varepsilon)=(UD_2)^{-1}(\varepsilon)$.

By Theorem \ref{P22} and sequence (\ref{braid3}), $G$ admits presentation with generators
$\{A_1, A_2, D_2, U, j(\alpha_3), j(\delta_{3}), j(\delta_{12}), j(\delta_{23}), j(\varepsilon)\}$
and relations (2), (3), $A_1A_2=A_2A_1$, $(A_2U)^2=(D_2U)^2=j(\delta^{-1}_{23}\delta^{-1}_{12})$ and:

\medskip
\noindent
$\mathrm{(i)}\ Uj(\alpha_3)U^{-1}=j(\delta_{23}\alpha^{-1}_3\delta_{12}\delta_{23});\
\mathrm{(ii)}\ Uj(\delta_3)U^{-1}=j(\delta^{-1}_3\delta_{12});\\
\mathrm{(iii)}\ Uj(\delta_{23})U^{-1}=j(\delta_{23});\
\mathrm{(iv)}\ Uj(\delta_{12})U^{-1}=j(\delta_{12});\quad
\mathrm{(v)}\ Uj(\varepsilon)U^{-1}=\\
j(\delta_3^{-1}\delta_{12}\delta_{23}\alpha^{-1}_3\varepsilon\alpha_3\delta^{-1}_{23}\delta^{-1}_{12}\delta_3);\
\mathrm{(vi)}\ D_2j(\alpha_3)D^{-1}_2=j(\delta^{-1}_{23}\delta^{-1}_3\varepsilon\delta_3\alpha_3);\\
\mathrm{(vii)}\ D_2j(\delta_3)D^{-1}_2=j(\delta^{-1}_{23}\delta_3\delta_{23});\
\mathrm{(viii)}\ D_2j(\delta_{23})D^{-1}_2=j(D_2(\delta^{-1}_3)\delta_3\delta_{23});\\
\mathrm{(ix)}\ D_2j(\delta_{12})D^{-1}_2=j(\delta_{12}\delta_{23}D_2(\delta^{-1}_{23}));\
\mathrm{(x)}\ D_2j(\varepsilon)D^{-1}_2=j(D_2(\alpha_3)\alpha^{-1}_3\delta_{23});\\
\mathrm{(xi)}\ A_2j(\alpha_3)A^{-1}_2=j(\alpha_3);\
\mathrm{(xii)}\ A_2j(\delta_3)A^{-1}_2=j(\delta_{12}\delta_{23}\alpha^{-1}_3\varepsilon\delta_3\alpha_3);\\
\mathrm{(xiii)}\ A_2j(\delta_{23})A^{-1}_2=j(\alpha^{-1}_3\delta_{23}\alpha_3);\qquad
\mathrm{(xiv)}A_2j(\delta_{12})A^{-1}_2=\\
j(\delta_{12}\delta_{23})A_2j(\delta_{23}^{-1})A^{-1}_2;\
\mathrm{(xv)}\ A_2j(\varepsilon)A^{-1}_2=j(\varepsilon);\
\mathrm{(xvi)}\ A_1j(\alpha_3)A_1^{-1}=j(\alpha_3);\
\mathrm{(xvii)}\ A_1j(\delta_3)A_1^{-1}=j(\delta_{12}\delta_{23}\alpha^{-1}_3\delta_3\alpha_3\delta_{23}^{-1});\
\mathrm{(xviii)}\ A_1j(\delta_{23})A_1^{-1}=j(\delta_{23});\
\mathrm{(xix)}\ A_1j(\delta_{12})A_1^{-1}=j(\delta_{12});\
\mathrm{(xx)}\ A_1j(\varepsilon)A_1^{-1}=j((UD_2)^{-1}(\varepsilon)).$

\medskip

It remains to check, that the relations (i-xx) above are consequences of (1-11) in Theorem \ref{P23} and
(\ref{sub1}), (\ref{sub2}), $j(\delta_3)=D_3$. We have:

\noindent $\mathrm{(i)}\Leftrightarrow
UA_3A_2^{-1}U^{-1}=(UA_3)^2A_3^{-1}A_2(A_2U)^{-2}\Leftrightarrow
(A_2U)^2=(UA_2)^2\Leftarrow (4);\quad
\mathrm{(ii)}\Leftrightarrow UD_3U^{-1}=D^{-1}_3(UD_2)^{-2}(UA_3)^{-2}
\stackrel{(8)}{=}\\
D_3^{-1}D_2^{-1}U^{-1}D_3UD_2D_3
\stackrel{(6)}{=}
D^{-1}_3U^{-1}D_3UD_3\Leftrightarrow (5);$
$\mathrm{(iii)}\Leftrightarrow (UA_3)^2=(A_3U)^2\Leftarrow (7,8);$
$\mathrm{(iv)}\Leftrightarrow U(UA_3)^2(UD_2)^2U^{-1}=(UA_3)^2(UD_2)^2\Leftarrow (4,7,8);$
$\mathrm{(v)}
\stackrel{(9)}{\Leftrightarrow}
U(A_3UD_2)^2U^{-1}=D^{-1}_3(UD_2)^{-1}(A_3UD_2)^2(UD_2)D_3\stackrel{(7)}{=}
A^{-1}_3D^{-1}_3A_3(UD_2)^2D_3
\stackrel{(4,6,7)}{\iff} D_3(A_3U)^2D_2=(A_3U)^2D_2D_3\stackrel{(8)}{\Leftrightarrow}\\
D_3(UD_2D_3)^{-2}D_2=(UD_2D_3)^{-2}D_2D_3
\Leftarrow(4,5,6)$;\\
$\mathrm{(vi)}\Leftrightarrow
D_2A_3A^{-1}_2D^{-1}_2=
(A_3U)^{-2}D^{-1}_3(A_3UD_2)^2D_3A_3A^{-1}_2
\stackrel{(3,7)}{\iff}\\
D_2=(A_3U)^{-2}D^{-1}_3A_3(UD_2)^2D_3UD^{-1}_2U^{-1}A_3^{-1}
\stackrel{(6,4)}{=}\\
(A_3U)^{-2}D^{-1}_3(A_3U)^2D_2D_3\stackrel{(8)}{=}
(UD_2D_3)^2D_3^{-1}(UD_2D_3)^{-2}D_2D_3\Leftarrow(4,5,6)$;
$\mathrm{(vii)}\Leftrightarrow
D_2D_3D^{-1}_2=
(UD_2D_3)^2D_3(UD_2D_3)^{-2}\Leftarrow(4,5,6);$\\
$\mathrm{(viii)}\Leftrightarrow
D_2(UD_2D_3)^{-2}D_2^{-1}=D_2D^{-1}_3D_2^{-1}D_3(UD_2D_3)^{-2}\Leftarrow (4,6);$\\
$\mathrm{(ix)}\Leftrightarrow
D_2(\delta_{12}\delta_{23})=\delta_{12}\delta_{23}\Leftrightarrow (UD_2)^2=(D_2U)^2
\Leftarrow (4);$\\
$\mathrm{(x)}\Leftrightarrow
D_2(A_3UD_2)^2D^{-1}_2=D_2A_3A^{-1}_2D^{-1}_2A_2A^{-1}_3(A_3U)^2\Leftarrow (3)$;\\
$(1)\Rightarrow \mathrm{(xi)};$
$\mathrm{(xii)}\Leftrightarrow
A_2D_3A^{-1}_2=(UD_2)^{-2}A_2A^{-1}_3(A_3UD_2)^2D_3A_3A^{-1}_2
\stackrel{(7)}{\Leftrightarrow}
A_2(UD_2)^2=(UD_2)^2A_2\Leftarrow (4)$;\
$\mathrm{(xiii)}\Leftrightarrow
A_2(UA_3)^2A_2^{-1}=\\
A_2A_3^{-1}(UA_3)^2A_3A_2^{-1}\stackrel{(1)}{\Leftrightarrow}
(UA_3)^2=(A_3U)^2\Leftarrow (7,8);$\\
$\mathrm{(xiv)}\Leftrightarrow A_2(UD_2)^2=(UD_2)^2A_2\Leftarrow(4);$
$(9)\Rightarrow \mathrm{(xv)}$;\\
$\mathrm{(xvii)}\Leftrightarrow A_1D_3A_1^{-1}=(UD_2)^{-2}A_2A^{-1}_3D_3A_3A^{-1}_2(UA_3)^{-2}
\stackrel{(1,2,4)}{\iff}\\
A_2A_1^{-1}D_3A_1A_2^{-1}=A_3(UD_2)^2D_3(UA_3)^2A_3^{-1}
\stackrel{(8)}{\Leftrightarrow}(10)$;\\
$(1,2,4)\Rightarrow
\mathrm{(xvi,xviii,xix)}$; $\mathrm{(xx)}\Leftrightarrow(11)$.
\foorp

\medskip

Let $\widetilde{F}=F_2^3$ be a subsurface of $F$ such that
boundary curve $c_i\colon S^1\to\bdr\widetilde{F}$ bounds in $F$ a
disc with puncture $p_i$ for $i=1,2,3$. We identify $\{A_1, A_2,
A_3, D_2, D_3, U\}$ with elements of $\M(\widetilde{F})$.

\begin{theorem}\label{23}
The group $\M(F_2^3)$ admits a presentation with generators
$\{A_1, A_2, A_3, D_2, D_3, U, C_1, C_2, C_3\}$ and relations
$(1-7)$, $(9-11)$ from Theorem \ref{P23} and
$\ (8')\ (UA_3)^2(UD_2D_3)^2=(C_1C_2C_3)^2,$
$\ C_iC_j=C_jC_i,\ C_iA_j=A_jC_i,\ C_iD_k=C_iD_k,\ C_iU=UC_i,$
for $i,j\in\{1,2,3\}$, $k\in\{2,3\}$.
\end{theorem}

\proof Let $H$ denote the subgroup of $\M(\widetilde{F})$
generated by the twists $\{C_1, C_2, C_3\}$. It is easy to see
that relations (1-7) and (10) are satisfied in
$\M(\widetilde{F})$.
In the proof of Theorem \ref{P23} we showed that
$j(\varepsilon)=(A_3UD_2)^2$ in $\PM^+(F,\{p_1,p_2,p_3\})$. On the other hand,
by Lemma \ref{j1}, $j(\varepsilon)$ is equal to a Dehn twist
$E$ about a generic curve $e$. Thus in $\M(\widetilde{F})$ we have $E(A_3UD_2)^{-2}\in H$. It can be checked that
in $\M(\widetilde{F})$ we have $A_2EA^{-1}_2=E$ and $A_1EA^{-1}_1=(UD_2)^{-1}E(UD_2)$, and hence (9) and (11) hold,
since $H$ is central.

Let $d_{23}$ and $l$ denote boundary curves of tubular
neighborhoods of the loops $\delta_{23}$ and
$\delta_3\delta_{23}$, such that in $\PM^+(F,\{p_1,p_2,p_3\})$ we
have $D_{23}=j(\delta_{23})$, $LD^{-1}_2=j(\delta_3\delta_{23})$.
The curves $d_{23}$ and $c_1$ bound in $\widetilde{F}$ a Klein
bottle with two holes, while $l$, $c_2$, $c_3$ bound a 4-holed
sphere, together with a curve bounding a M\"obius strip. Thus we
have lantern relation $LC_2C_3=D_{23}D_2D_3$ and relation
$(UA_3)^2=(UL)^2=C_1D_{23}$ from Lemma \ref{K2rel}. Now
$(UA_3)^2=(UL)^2=(UD_{23}D_2D_3(C_2C_3)^{-1})^2=
D^2_{23}(C_2C_3)^{-2}(UD_2D_3)^2=
(UA_3)^4(C_1C_2C_3)^{-2}(UD_2D_3)^2\Leftrightarrow (8')$.

Theorem \ref{23} follows from Theorem \ref{P23} and sequence (\ref{Mgn}).
\foorp

\subsection{Sporadic surfaces of genus 3.}

\begin{figure}%[!htbp]
\begin{center}
\input{fig21_torus}
\caption{\label{torus}The torus $T_3$}
\end{center}
\end{figure}
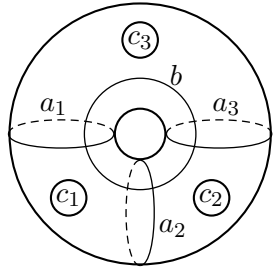

Consider a torus with three holes $T_3$ represented in Figure \ref{torus}, and let $T_2$ be the torus with
two holes obtained by gluing a disc to the boundary of
$T_3$, along the curve $c_2$. We fix in $T_3$ and $T_2$ the orientation induced by the standard orientation of the plane of
Figure \ref{torus}, and let $C_i$, $A_i$, $B$, $i=1,2,3$ denote Dehn twists along the curves in the figure, right with respect to
that orientation. The next theorem follows from the main result of \cite{G}.

\begin{theorem}\label{Tpres}
The group $\M(T_3)$ admits presentation with generators
$\{C_i, A_i, B\,|\, i=1,2,3\}$ and relations:
\begin{equation}\label{C}
C_iC_j=C_jC_i,\ C_iA_j=A_jC_i,\ C_iB=BC_i,
\end{equation}
\begin{equation}\label{B}
A_iA_j=A_jA_i,\
A_iBA_i=BA_iB,
\end{equation}
for $i,j=1,2,3$, and
\begin{equation}\label{S}
(A_1A_2A_3B)^3=C_1C_2C_3.
\end{equation}
A presentation for $\M(T_2)$ may be obtained by adding to the
above presentation relations $C_2=1$ and $A_2=A_3$.\foorp
\end{theorem}

\begin{rem}\label{star}
The relation (\ref{S}) is called ``star'' in \cite{G}. In $\M(T_2)$ it takes form $(A_1A_2^2B)^3=C_1C_3$, and it follows
from relations (\ref{B}) that $(A_1A_2^2B)^3=(A_1^2A_2B)^3$.
\end{rem}

\begin{figure}%[!htbp]
\begin{center}
\begin{tabular}{cc}
\input{fig22_F32} & \input{fig23_F32a}
\end{tabular}
\caption{\label{F32}The surface $\widetilde{F}=F_3^2$}
\end{center}
\end{figure}
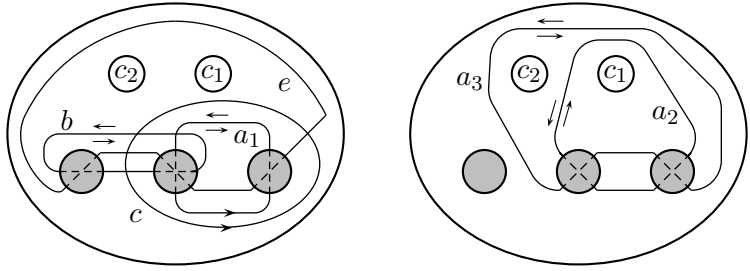

Let $\widetilde{F}=F_3^2$ be the surface obtained by gluing a
M\"obius strip $M$ to the boundary of $T_3$ along $c_3$. We
identify $\widetilde{F}$ with the surface represented in Figure \ref{F32}, where $M$ is a regular
neighborhood of the one-sided curve $e$. Consider an embedding
$\phi\colon K\to\widetilde{F}$, where $K$ is the holed Klein
bottle in Figure \ref{U}, such that $\phi\circ c=c$ and $\phi\circ
a_1=a_1$. We define $U=\phi_\ast(U)$, where $U\colon K\to K$ is
defined in Subsection \ref{s2}. We identify $A_1$, $A_2$, $A_3$,
and $B$ with elements of $\M(\widetilde{F})$ (the directions of these twists are indicated by arrows in
Figure \ref{F32}).

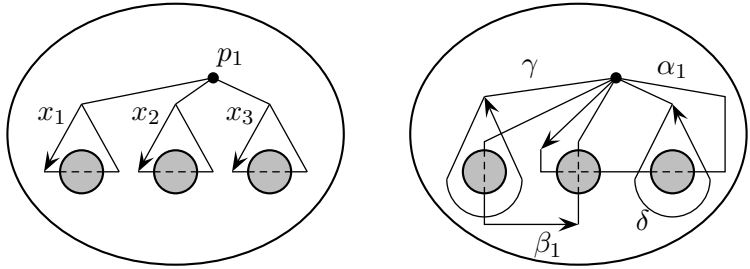
\begin{figure}%[!htbp]
\begin{center}
\begin{tabular}{cc}
\input{fig24_Pi31} & \input{fig25_Pip31}
\end{tabular}
\caption{\label{Pi31}Generators of $\pi_1(F,p_1)$ and
$\pi^+_1(F,p_1)$.}
\end{center}
\end{figure}

Let $F=F_3^0$ be the closed surface obtained by gluing two discs to $\bdr{\widetilde{F}}$.
We fix a point $p_1\in F$ inside the disc bounded by $c_1$, and $p_2\in
F$ inside the disc bounded by $c_2$.

\begin{theorem}\label{P31}
The group $\PM^+(F,\{p_1\})$ admits a presentation with generators
$\{A_1, A_2, B, U\}$ and relations:
\begin{eqnarray*}
&&(1)\ A_1A_2=A_2A_1;\quad (2)\ A_1BA_1=BA_1B,\ A_2BA_2=BA_2B;\\
&&(3)\ UA_1U^{-1}=A_1^{-1};\quad (4)\  UBU^{-1}=A_2^{-1}B^{-1}A_2;\\
&&(5)\ (UA_2)^2=1;\quad (6)\  (A_1A_2^2B)^3=1.
\end{eqnarray*}
\end{theorem}
\proof Let us denote $G=\PM^+(F,\{p_1\})$. Notice that relations (1-6) are
satisfied in $G$: (1) is obvious; $(2,3,4)$ are relations of type
$ht_ah^{-1}=t^{\pm 1}_{h(a)}$; (5) follows from Lemma
\ref{K2rel}; (6) is a star relation (cf. Remark \ref{star}).

Consider the exact sequence (\ref{braid3}):
$$1\to\pi^+(F,p_1)\stackrel{j}{\to}G\to\M(F)\to 1.$$
The fundamental group $\pi_1(F,p_1)$ is generated by the loops
$x_1$, $x_2$, $x_3$ in Figure \ref{Pi31} satisfying one defining
relation $x^2_3x^2_2x^2_1=1$. Now $\{1,x_3\}$ is a Schreier system
of representatives of cosets of $\pi^+_1(F,p_1)$ and by the
Reidemeister-Schreier method we obtain that the last group is
generated by $u_1=x_1x_3^{-1}$, $u_2=x_2x_3^{-1}$, $u_3=x_3x_1$,
$u_4=x_3x_2$ and $u_5=x_3^2$ satisfying two defining relations:
$u_5u_2u_4u_1u_3=1$, $u_5u_4u_2u_3u_1=1$. After Tietze
transformations (c.f. \cite{MKS}) we obtain
$$\pi^+_1(F,p_1)=\langle\alpha_1, \beta_1, \delta, \gamma\,|\,
\beta_1^{-1}\delta^{-1}\gamma^{-1}\alpha^{-1}\delta\alpha_1\beta_1\gamma=1\rangle,$$
where $\alpha_1=u_4$, $\delta=u_5$, $\beta_1=u_2u_3$,
$\gamma=u_1u_3$ are the loops in Figure \ref{Pi31}.
It follows from Theorem 2 of \cite{BC} that $\M(F)$ admits a
presentation with generators $\{A_1, B, U\}$ and relations
$A_1BA_1=BA_1B$, $UA_1U^{-1}=A_1^{-1}$,
$UBU^{-1}=A_1^{-1}B^{-1}A_1$, $U^2=1$, $(A_1^3B)^3=1$. The last
relation is a special form of the star relation (\ref{S}) and it
can be checked that in $G$ we have
$(A_1^3B)^3=j(\beta^{-1}_1\alpha_1\beta_1\alpha^{-1}_1)$. We also
have $UBU^{-1}A_1^{-1}BA_1=j(\beta^{-1}_1\alpha^{-1}_1).$
By Lemma \ref{pres}, $G$ admits presentation with
generators $\{A_1, B, U, j(\alpha_1), j(\beta_1), j(\gamma),
j(\delta)\}$ and relations:
\medskip

\noindent
$\mathrm{(i)}\ A_1BA_1=BA_1B;\ \mathrm{(ii)}\ UA_1U^{-1}=A_1^{-1};\
\mathrm{(iii)}\ UBU^{-1}A_1^{-1}BA_1=j(\beta^{-1}_1\alpha^{-1}_1);
\mathrm{(iv)}\ U^2=j(\gamma);\ \mathrm{(v)}\ (A_1^3B)^3=j(\beta^{-1}_1\alpha_1\beta_1\alpha^{-1}_1);\\
\mathrm{(vi)}\ j(\beta_1^{-1}\delta^{-1}\gamma^{-1}\alpha^{-1}\delta\alpha_1\beta_1\gamma)=1;
\mathrm{(vii)}\ A_1j(\alpha_1)A^{-1}_1=j(\alpha_1);\\ \mathrm{(viii)}\ A_1j(\beta_1)A^{-1}_1=j(\alpha^{-1}_1\beta_1);
\mathrm{(ix)}\ A_1j(\gamma)A^{-1}_1=j(\gamma);\ \mathrm{(x)}\ A_1j(\delta)A^{-1}_1=j(\gamma^{-1}\alpha^{-1}_1\delta\alpha_1);\
\mathrm{(xi)}\ Bj(\alpha_1)B^{-1}=j(\alpha_1\beta_1);\ \mathrm{(xii)}\ Bj(\beta_1)B^{-1}=j(\beta_1);\\
\mathrm{(xiii)}\ Bj(\gamma)B^{-1}=j(\beta^{-1}_1\gamma\delta\beta_1);\
\mathrm{(xiv)}\ Bj(\delta)B^{-1}=j(\delta);\\
\mathrm{(xv)}\ Uj(\alpha_1)U^{-1}=j(\alpha^{-1}_1\gamma^{-1});\
\mathrm{(xvi)}\ Uj(\beta_1)U^{-1}=j(\gamma\delta\alpha_1\beta_1);\\
\mathrm{(xvii)}\ Uj(\gamma)U^{-1}=j(\gamma);\
\mathrm{(xviii)}\ Uj(\delta)U^{-1}=j(\delta^{-1}\gamma^{-1}).
$

We have:
\begin{equation}\label{Rj1}
j(\gamma)=U^2,\quad j(\alpha_1)=A_2A_1^{-1},\quad j(\beta_1)=A_1A_2^{-1}BA_2A^{-1}_1B^{-1}.
\end{equation}
It can be checked that $U^{-1}B(\alpha_1)=\delta\beta_1$, and hence
\begin{equation}\label{Rj2}
j(\delta)=U^{-1}BA_2A_1^{-1}B^{-1}UBA_1A^{-1}_2B^{-1}A_2A^{-1}_1.
\end{equation} %
Let $H$ denote the subgroup of $G$ generated by $\{A_1, A_2, B\}$.
Consider the homomorphism $i_\ast\colon\M(T_2)\to G$ induced by
the inclusion of $T_2$ in $F$. It can be proved, using the same
methods as in the proof of Lemma \ref{kernel}, that $\ker
i_\ast$ is generated by $\{C_1, C_3\}$. Now it follows from
Theorem \ref{Tpres} that $i_\ast(\M(T_2))=H$ and every relation in
$H$ is a consequence of $(1,2,6)$.

We will show that relations (i)-(xviii) after replacing
$j(\alpha_1)$, $j(\beta_1)$, $j(\gamma)$ and $j(\delta)$ by
expressions (\ref{Rj1}, \ref{Rj2}), are consequences of (1-6).
Relations (i,ii) are the same as $(2,3)$; (iv,xi,xvii) are
trivial; (v,vii,viii,xii) are relations in $H$, hence they follow
from $(1,2,6)$. We have\\
$UBU^{-1}A_1^{-1}BA_1\stackrel{(4)}{=}A_2^{-1}B^{-1}A_2A_1^{-1}BA_1
\stackrel{(2)}{=}BA^{-1}_2A_1B^{-1}\Leftrightarrow\mathrm{(iii)};\
(3)\Rightarrow\mathrm{(ix)};\ (1,3,5)\Rightarrow\mathrm{(xv)}$;
(x, xiii, xiv) can easily be reduced to relations in $H$, by using
(1-4).

Let $X=UBA_2^{-1}A_1B^{-1}A_1^{-1}A_2BA_1^{-1}A_2B^{-1}U$, and note that
to prove (1-6) $\Rightarrow$ (xvi, xviii), it suffices to show (1-6) $\Rightarrow X\in H$.
By $(2,3,4)$ we have $UA_1U^{-1}\in H$, $UBU^{-1}\in H$, $BA_2B^{-1}=UB^{-1}U^{-1}$,
thus $X\in H\Leftrightarrow UA_2^{-1}B^{-1}A^{-1}_1A_2BA_1^{-1}A_2B^{-1}U\in H\Leftrightarrow
UA_2^{-1}B^{-1}A_1^{-2}B^{-1}U\in H.$ It can be checked that from (1,2,6) follows $A_2^{-1}B^{-1}A_1^{-2}B^{-1}A_2^{-1}=A_1BA_1^2BA_1$,
hence $X\in H\Leftrightarrow UA_2U\in H\Leftarrow(5)$.
Finally, we have
$j(\beta_1^{-1}\delta^{-1}\gamma^{-1}\alpha^{-1}_1\delta\alpha_1\beta_1\gamma)\stackrel{\mathrm{(xvi)}}{=}
j(\beta_1^{-1}\delta^{-1})U^{-2}j(\alpha_1^{-1})U^{-1}j(\beta_1)U=\\
U^{-1}j(\beta^{-1}_1\alpha^{-1}_1)U^{-1}j(\alpha^{-1}_1)U^{-1}j(\beta_1)U$, thus
$\mathrm{(vi)}\Leftrightarrow (Uj(\alpha_1))^2=1\Leftarrow (1,3,5)$.
\foorp

\begin{theorem}\label{M31}
The group $\M(F_3^1)$ admits a presentation with generators
$\{A_1$, $A_2$, $B$, $U\}$, relations $(1-4)$ from Theorem \ref{P31} and
$(A_2U)^2=(UA_2)^2=(A_1^2A_2B)^3.$
\end{theorem}

\proof Consider the surface $F_3^1$ obtained by gluing a disc to
the boundary of $\tilde{F}$ along $c_2$. Observe that relations
$(1-4)$ from Theorem \ref{P31} are satisfied in $\M(F_3^1)$, and
we have $(A_1^2A_2B)^3=C_1$ (star) and $(A_2U)^2=C_1$ (Lemma
\ref{K2rel}). After replacing the generator $C_1$ in the
presentation of $\M(F_3^1)$ resulting from applying Lemma
\ref{pres} to sequence (\ref{Mgn}), we obtain Theorem
\ref{M31}.\foorp

\begin{figure}%[!htbp]
\begin{center}
\begin{tabular}{cc}
\input{fig26_Pi32} & \input{fig27_Pip32}
\end{tabular}
\caption{\label{Pi32}Generators of $\pi_1(F\backslash\{p_1\},p_2)$ and
$\pi^+_1(F\backslash\{p_1\},p_2)$.}
\end{center}
\end{figure}
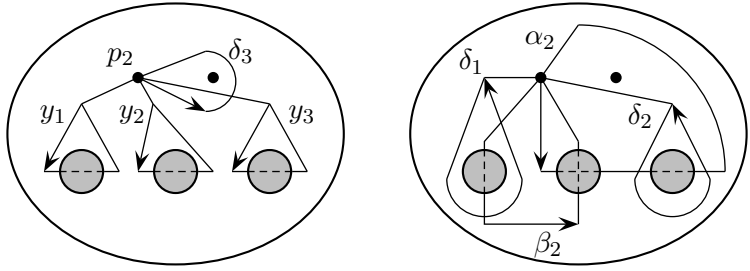

\begin{theorem}\label{P32}
The group $\PM^+(F,\{p_1,p_2\})$ admits a presentation with generators
$\{A_1, A_2, A_3, B, D_1, D_2, D_3, U \}$
and relations:
\begin{eqnarray*}
&&(1)\ A_iA_j=A_jA_i,\ i,j=1,2,3;\
(2)\ A_iBA_i=BA_iB,\ i=1,2,3;\\
&&(3)\ UA_1U^{-1}=A_1^{-1};\ (4)\ UBU^{-1}=A_3^{-1}B^{-1}A_3;\\
&&(5)\ UD_1=D_1U;\ (6)\ UD_3=D_3U;\ (7)\ BD_2=D_2B;\\
&&(8)\ (UA_2)^2=D_1;\ (9)\ (A_1^2A_3B)^3=(UA_3)^2=D_3;\\
&&(10)\ A_2^{-1}UD_2U^{-1}A_2=UB^{-1}D_1^{-1}BU^{-1};\\
&&(11)\ (UD_2)^2D_1D_3=U^2;\ (12)\ (A_1A_2A_3B)^3=1.
\end{eqnarray*}
\end{theorem}

\proof Let us denote $G=\PM^+(F,\{p_1,p_2\})$. The fundamental
group $\pi_1(F\backslash\{p_1\},p_2)$ is free on generators $y_1$,
$y_2$, $y_3$ in Figure \ref{Pi32}. Now $\{1,y_3\}$ is a Schreier
system of representatives of cosets of
$\pi^+_1(F\backslash\{p_1\},p_2)$ and by the Reidemeister-Schreier
method we obtain that the last group is freely generated by
$v_1=y_1y_3^{-1}$, $v_2=y_2y_3^{-1}$, $v_3=y_3y_1$, $v_4=y_3y_2$,
$v_5=y_3^2$. It follows that $\pi^+_1(F\backslash\{p_1\},p_2)$ is
free on generators $\delta_2=v_5$, $\delta_1=v_1v_3$,
$\beta_2=v_2v_3$, $\delta_3=\delta_2v_2v_4\delta_1$,
$\alpha_2=\delta_3v_4$ (see Figure \ref{Pi32}). We introduce Dehn
twists $D_i=j(\delta_i)$, $i=1,2,3$. We also have
$$j(\alpha_2)=A_3A^{-1}_2,\quad
j(\beta_2)=A_3^{-1}A_2BA_2^{-1}A_3B^{-1}.$$

Let us check that relations (1-12) are satisfied in $G$:
$(1,2,12)$ follow from  Theorem \ref{Tpres}; $(3,4,10)$ are
relations of type $ht_ah^{-1}=t^{\pm 1}_{h(a)}$; $(5,6,7)$ are
obvious; $(8,9)$ follow from Lemma \ref{K2rel} and star relation; (11)
follows from the equality
$U(\delta_2)=\delta^{-1}_2\delta^{-1}_3\delta^{-1}_1$ and
relations $(5,6)$.

By Theorem
\ref{P31} and Lemma \ref{pres} for sequence (\ref{braid3}), $G$
admits a presentation with generators $\{A_1, A_2, B, U,
j(\alpha_2), j(\beta_2), j(\delta_i)\,|\,i=1,2,3 \}$ and relations
$(1,2,3)$ and:

\noindent
$\mathrm{(i)}\ UBU^{-1}A_2^{-1}BA_2=j(\beta^{-1}_2\alpha^{-1}_2);\
\mathrm{(ii)}\ (UA_2)^2=j(\delta_1);\ \mathrm{(iii)}\ (A_1A_2^2B)^3=j(\beta^{-1}_2\delta_3^{-1}\alpha_2\beta_2\alpha^{-1}_2);\
\mathrm{(iv)}\ A_1j(\alpha_2)A^{-1}_1=A_2j(\alpha_2)A^{-1}_2=j(\alpha_2);\\
\mathrm{(v)}\ A_1j(\beta_2)A^{-1}_1=j(\alpha^{-1}_2\delta_3\beta_2);\
\mathrm{(vi)}\ A_1j(\delta_1)A^{-1}_1=A_2j(\delta_1)A^{-1}_2=\\
Uj(\delta_1)U^{-1}=j(\delta_1);\
\mathrm{(vii)}\ A_1j(\delta_3)A^{-1}_1=Bj(\delta_3)B^{-1}=Uj(\delta_3)U^{-1}=j(\delta_3);\
\mathrm{(viii)}\ A_1j(\delta_2)A^{-1}_1=j(\delta^{-1}_3\delta^{-1}_1\alpha^{-1}_2\delta_3\delta_2\delta_3^{-1}\alpha_2);\
\mathrm{(ix)}\ A_2j(\beta_2)A^{-1}_2=j(\alpha^{-1}_2\beta_2);\
\mathrm{(x)}\ A_2j(\delta_3)A^{-1}_2=j(\alpha^{-1}_2\delta_3\alpha_2);\
\mathrm{(xi)}\ A_2j(\delta_2)A^{-1}_2=\\
j(\alpha^{-1}_2\delta^{-1}_3\alpha_2\delta_3\delta_2\beta_2\delta^{-1}_1\beta_2^{-1}\alpha^{-1}_2\delta_3\alpha_2);
\mathrm{(xii)}\ Bj(\alpha_2)B^{-1}=j(\alpha_2\beta_2);\\
\mathrm{(xiii)}\ Bj(\beta_2)B^{-1}=j(\beta_2);\
\mathrm{(xiv)}\ Bj(\delta_1)B^{-1}=j(\beta_2^{-1}\delta_1\delta_3\delta_2\beta_2);\\
\mathrm{(xv)}\ Bj(\delta_2)B^{-1}=j(\delta_2);\
\mathrm{(xvi)}\ Uj(\alpha_2)U^{-1}=j(\delta_3\alpha^{-1}_2\delta^{-1}_1);\\
\mathrm{(xvii)}\ Uj(\beta_2)U^{-1}=j(\delta_1\delta_3\delta_2\delta^{-1}_3\alpha_2\beta_2);\
\mathrm{(xviii)}\ Uj(\delta_2)U^{-1}=j(\delta^{-1}_2\delta^{-1}_3\delta^{-1}_1).$

\medskip

We will show that relations
(i-xviii) after substituting $j(\alpha_2)=A_3A_2^{-1}$, $j(\beta_2)=A_3^{-1}A_2BA_2^{-1}A_3B^{-1}$, $j(\delta_i)=D_i$,
are consequences of (1-12).

Let $H$ denote the subgroup of $G$ generated by $\{A_1, A_2, A_3, B\}$.
As in the proof of Theorem \ref{P31}, we have $H=i_\ast(\M(T_3))$, where
$i_\ast$ is the homomorphism induced by the inclusion of $T_3$ in $F$, and every relation in $H$ is a consequence of
$(1,2,12)$, by Theorem \ref{Tpres}.
Note that by the star relation (9), $D_3\in H$.

Relations (i - vii, ix, x, xii, xiii, xv) follow easily from (1-12) or are relations in $H$;
$(8,9) \Rightarrow$ (xvi);
$(5,6,11) \Rightarrow$ (xviii); by $(5,8)$ we have $A_2D_1=D_1A_2$ and (xiv) $\stackrel{\mathrm{(xviii)}}{\Leftrightarrow}
j(\beta_2)BD_1B^{-1}j(\beta^{-1}_2)=Uj(\delta_2^{-1})U^{-1}\Leftrightarrow
A_3^{-1}BA_3D_1A_3^{-1}B^{-1}A_3=A_2^{-1}UD^{-1}_2U^{-1}A_2 \Leftarrow (4,5,10)$;\\
(xvii) $\stackrel{\mathrm{(xviii)}}{\iff}
Uj(\beta_2)U^{-1}=UD^{-1}_2U^{-1}D^{-1}_3j(\alpha_2\beta_2)
\stackrel{(5,6,11)}{\iff}\\
j(\beta_2)=UD_2U^{-1}D_1U^{-1}j(\alpha_2\beta_2)U
\stackrel{\mathrm{(xiv, xviii)}}{\iff}\\
BD_1B^{-1}=
j(\beta_2^{-1})D_1U^{-1}j(\alpha_2\beta_2)U\
\Leftrightarrow
D^{-1}_1A_2A^{-1}_3BA_3A^{-1}_2D_1=\\
U^{-1}j(\alpha_2\beta_2)UB
\stackrel{(4,8)}{\Leftrightarrow}
A_2^{-1}B^{-1}A_2=j(\alpha_2\beta_2)UBU^{-1}\
\Leftarrow UBU^{-1}\in H\Leftarrow (4);\quad$
(viii) $\stackrel{\mathrm{(vii)}}{\Leftrightarrow}
A_1(\delta_3\delta_2\delta^{-1}_3)=
\delta^{-1}_1\alpha^{-1}_2\delta_3\delta_2\delta^{-1}_3\alpha_2\delta^{-1}_3
\stackrel{\mathrm{(xvii)}}{\Leftrightarrow}\\
A_1(\delta^{-1}_1U(\beta_2)\beta^{-1}_2\alpha^{-1}_2)=
\delta^{-1}_1\alpha^{-1}_2\delta^{-1}_1U(\beta_2)\beta^{-1}_2\delta^{-1}_3
\stackrel{\mathrm{(iv,v,vi)}}{\iff}\\
 A_1U(\beta_2)=\alpha^{-1}_2\delta_1^{-1}U(\beta_2)
\stackrel{(8)}{\Leftrightarrow}
A_1Uj(\beta_2)U^{-1}A^{-1}_1=A_3^{-1}U^{-1}A^{-1}_2j(\beta_2)U^{-1}\stackrel{(3)}{\Leftrightarrow}
A_2UA_3UA^{-1}_1j(\beta_2)A_1=j(\beta_2)\Leftarrow UA_3U\in H \Leftarrow
(9)$;\
(xi) $\stackrel{\mathrm{(xiv, xviii)}}{\iff} A_3D_2A^{-1}_3=D^{-1}_3A_3A_2^{-1}D_3D_2B^{-1}UD_2U^{-1}BA_2A_3^{-1}D_3
\stackrel{(\mathrm{vii},7,11)}{\iff}D_3D_2D_3^{-1}=A_2^{-1}B^{-1}D^{-1}_1BA_2
\stackrel{\mathrm{(xvii)}}{\iff} D_1^{-1}Uj(\beta_2)U^{-1}j(\beta_2^{-1}\alpha^{-1}_2)=A^{-1}_2B^{-1}D^{-1}_1BA_2
\stackrel{(8)}{\Leftrightarrow} U^{-1}A^{-1}_3BA_3A^{-1}_2B^{-1}U^{-1}BA_3^{-1}A_2B^{-1}A^{-1}_2=B^{-1}D^{-1}_1B
\Leftarrow (2,4,8)$.
\foorp

\begin{theorem}\label{M32}
The group $\M(F_3^2)$ admits a presentation with generators\\
$\{A_1, A_2, A_3, B, D_1, D_2, D_3, U, C_1, C_2\}$ and relations
$(1-7,9,10)$ from Theorem \ref{P32} and
\begin{eqnarray*}
&(8')\ & (UA_2)^2=D_1C_1,\\
&(11')\ & (UD_2)^2D_1D_3=U^2C_1C^2_2,\\
&(12')\ & (A_1A_2A_3B)^3=C_1C_2=C_2C_1,
\end{eqnarray*}
$C_iA_j=A_jC_i,\ C_iD_k=D_kC_i,\ C_iB=BC_i\ C_iU=UC_i,$
for $i=1,2$, $j,k=1,2,3$.
\end{theorem}

\proof The relations $(1-7,9,10)$ from Theorem \ref{P32} are satisfied in $\M(\widetilde{F})=\M(F_3^2)$;
$(8')$ follows from Lemma \ref{K2rel}; $(12')$ is the star relation;
$(11')$ follows from Lemma \ref{K2rel} and lantern relation
$C_1C_2U^2=((UD_2)^2C_2^{-1})D_1D_3$. Now Theorem \ref{M32} follows from Theorem \ref{P32} and Lemma \ref{pres} for
sequence (\ref{braid3}).\foorp

%% file: fig28_sl.tex
\pspicture*(3.5,2)
\rput[t](.5,.3){\small$a_1$}
\psline[linewidth=.5pt](.25,.4)(3.25,.4)
\rput[b](2.5,1.1){\small$\alpha$}
\psline[linewidth=.5pt](.25,1)(3.25,1)
\rput[b](.5,1.7){\small$a_2$}
\psline[linewidth=.5pt](.25,1.6)(3.25,1.6)
\pscircle*(1.75,1){.06}%\rput[lb](1.8,1.1){\small$p_k$}
\psline[linewidth=.5pt, arrowsize=4pt 3]{<-}(.75,1)(1.8,1)
\psline[linewidth=.4pt]{->}(1.55,.5)(1.95,.5)
\psline[linewidth=.4pt]{<-}(1.55,.3)(1.95,.3)
\psline[linewidth=.4pt]{<-}(1.55,1.7)(1.95,1.7)
\psline[linewidth=.4pt]{->}(1.55,1.5)(1.95,1.5)
\endpspicture

%% file: fig8_lantern.tex
\pspicture*(3.75,3.2)
\pscircle(2,1.6){1.5}\rput(.65,.45){\small$a_0$}
\pscircle(2,2.4){.25}\rput(2,2.4){\small$a_2$}
\pscircle(1.4,1.2){.25}\rput(1.4,1.2){\small$a_1$}
\pscircle(2.6,1.2){.25}\rput(2.6,1.2){\small$a_3$}
\psbezier[linewidth=.5pt](1.15,.7)(.65,.95)(1.75,3.15)(2.25,2.9)
\psbezier[linewidth=.5pt](1.15,.7)(1.65,.45)(2.75,2.65)(2.25,2.9)
\rput[r](1.3,2.1){\small$a_{12}$}
\psbezier[linewidth=.5pt](2.85,.7)(2.35,.45)(1.25,2.65)(1.75,2.9)
\psbezier[linewidth=.5pt](2.85,.7)(3.35,.95)(2.25,3.15)(1.75,2.9)
\rput[l](2.7,2.1){\small$a_{23}$}
\psellipse[linewidth=.5pt](2,1.2)(1.2,.5)
\rput[t](2.05,.6){\small$a_{13}$}
\endpspicture

%% file: fig9_Pip13.tex
\pspicture*(4.5,4.5)
\rput[tl](.4,3.9){\small$c_1$}
\pscircle(2.25,2.25){2}
\pscircle*[linecolor=lightgray](2.75,1.75){.3}
\pscircle(2.75,1.75){.3}
\pscircle*(2.75,3){.075}
\pscircle*(1.25,3){.075}
%\pscircle*(2,3.85){.075}
%\rput[l](2.15,3.85){\small$p_4$}
\rput[bl](2.8,3.1){\small$p_3$}
\rput[t](1.25,2.85){\small$p_2$}
\psline[linewidth=.5pt](2.75,3)(2.1,1.5)
\psline[linewidth=.5pt, arrowsize=4pt 3]{->}(2.75,3)(3.4,1.5)
\psarc[linewidth=.5pt](2.75,1.5){.65}{180}{360}
\rput[l](3.2,2.4){\small$\alpha_3$}
\psline[linewidth=.5pt, arrowsize=4pt 3]{->}(2.75,3)(1.3,3.65)
\psline[linewidth=.5pt](2.75,3)(1.3,2.35)
\psarc[linewidth=.5pt](1.3,3){.65}{90}{270}
\rput[r](2.45,3){\small$\alpha_{23}$}
\psarc[linewidth=.5pt](1.25,3){.4}{0}{180}
\psline[linewidth=.5pt](1.25,1.4)(1.65,3)
\psline[linewidth=.5pt, arrowsize=4pt 3]{->}(1.25,1.4)(1.65,3)
\psline[linewidth=.5pt](1.25,1.4)%(1.05,2.2)
%\psline[linewidth=.5pt](1.05,2.2)
(.85,3)
\rput[r](1.1,1.7){\small$\beta_{23}$}
\psline[linewidth=.5pt](1.25,1.4)(1.25,1.25)(2.75,1.25)(2.75,1.45)
\psline[linewidth=.5pt,linestyle=dashed,dash=3pt 2pt](2.75,1.45)(2.75,2.05)
\psline[linewidth=.5pt](2.75,2.05)(2.75,3)
%
%\psbezier[linewidth=.5pt](2.75,3)(1.5,.5)(4,.5)(2.75,3)
\endpspicture

%% file: fig10_M13.tex
\pspicture*(4.5,4.5)
\pscircle(2.25,2.25){2}
\pscircle*[linecolor=lightgray](2.75,1.75){.3}
\pscircle(2.75,1.75){.3}
\pscircle*(2.75,3){.075}
\pscircle*(1.25,3){.075}
%\pscircle*(2,3.85){.075}
%
\psline[linewidth=.5pt,linearc=.2](2.55,2)(2.4,2.25)(2.4,3.25)
\psline[linewidth=.5pt,linearc=.2](2.95,2)(3.1,2.25)(3.1,3.25)
%\psline[linewidth=.5pt](2.4,2.25)(2.4,3.25)
%\psline[linewidth=.5pt](3.1,2.25)(3.1,3.25)
\psarc[linewidth=.5pt](2.75,3.25){.35}{0}{180}
\psline[linewidth=.5pt,linestyle=dashed,dash=3pt 2pt](2.55,2)(2.95,1.5)
\psline[linewidth=.5pt,linestyle=dashed,dash=3pt 2pt](2.95,2)(2.55,1.5)
%
%\psline[linewidth=.5pt](2.95,1.5)(2.95,1.2)
\psarc[linewidth=.5pt](1.25,3.25){.35}{0}{180}
%\psline[linewidth=.5pt](.9,3.25)(.9,1.2)
\psline[linewidth=.5pt,linearc=.2](1.6,3.25)(1.6,1.35)(2.55,1.35)(2.55,1.5)
\psline[linewidth=.5pt,linearc=.2](.9,3.25)(.9,1.2)(2.95,1.2)(2.95,1.5)
\psline[linewidth=.5pt](.8,3.25)(3.5,3.25)
\psline[linewidth=.5pt](.8,2.75)(3.5,2.75)
\psarc[linewidth=.5pt](.8,3){.25}{90}{270}
\psarc[linewidth=.5pt](3.5,3){.25}{270}{90}
\psline[linewidth=.5pt](3.3,1.15)(3.3,3.3)
\psline[linewidth=.5pt](2.1,1.15)(2.1,3.3)
\psarc[linewidth=.5pt](2.7,3.3){.6}{0}{180}
\psarc[linewidth=.5pt](2.7,1.15){.6}{180}{360}
\rput[b](2.7,.7){\small$a_{3}$}
\rput[bl](1,1.3){\small$b_{23}$}
\rput[tl](3.5,2.65){\small$a_{23}$}
%
%
%
%
%\psbezier[linewidth=.5pt](5,2.75)(5,2.25)(3.9,2.25)(3.9,1.5)
\endpspicture

%% file: fig11_Pip14a.tex
\pspicture*(4.5,4.5)
%
%\rput[tl](.4,3.9){\small$q_1$}
\rput[tl](.4,3.9){\small$c_1$}
\pscircle(2.25,2.25){2}
\pscircle*[linecolor=lightgray](2.75,1.75){.3}
\pscircle(2.75,1.75){.3}
\pscircle*(2.75,3){.075}
\pscircle*(1.25,3){.075}
\pscircle*(2,3.85){.075}
\rput[r](1.85,3.85){\small$p_4$}
%\rput[bl](2.8,3.1){\small$q_3$}
%\rput[t](1.25,2.85){\small$q_2$}
%
%\psarc[linewidth=.5pt](2.75,3){.6}{180}{360}
%\psarc[linewidth=.5pt](2.75,3){.4}{270}{450}
%\psline[linewidth=.5pt](2,3.85)(2.75,3.4)
%\psline[linewidth=.5pt, arrowsize=4pt 3]{->}(2,3.85)(2.1,3)
%\psline[linewidth=.5pt](2.1,3)(2.75,2.6)
%
\psarc[linewidth=.5pt](2.75,3){.4}{270}{450}
\psline[linewidth=.5pt](2.15,3)(2.75,3.4)
\psline[linewidth=.5pt, arrowsize=4pt 3]{->}(2.15,3)(2.75,2.6)
\psline[linewidth=.5pt](2,3.85)(2.15,3)
\rput[tr](2.4,2.75){\small$\alpha_{34}$}
\psarc[linewidth=.5pt](1.25,3){.4}{90}{270}
\psline[linewidth=.5pt](2,3.85)(1.85,3)
\psline[linewidth=.5pt, arrowsize=4pt 3]{->}(1.85,3)(1.25,3.4)
\psline[linewidth=.5pt](1.85,3)(1.25,2.6)
\rput[t](1.25,2.5){\small$\alpha_{24}$}
%
%\psarc[linewidth=.5pt](2.75,3){.2}{0}{180}
%\psline[linewidth=.5pt, arrowsize=4pt 3]{->}(2.75,2.4)(2.55,3)
%\psline[linewidth=.5pt](2.75,2.4)(2.95,3)
%
\psline[linewidth=.5pt](2,3.85)(3.75,3.3)(3.75,1.75)
\psarc[linewidth=.5pt](2.75,1.75){.5}{90}{270}
\psline[linewidth=.5pt, arrowsize=4pt 3]{->}(3.75,1.75)(2.75,2.25)
\psline[linewidth=.5pt](3.75,1.75)(2.75,1.25)
\rput[tl](3.4,1.4){\small$\alpha_{4}$}
%
%\psline[linewidth=.5pt]
%\psline[linewidth=.5pt]
%\psbezier[linewidth=.5pt](2.75,3)(1.5,.5)(4,.5)(2.75,3)
\endpspicture

%% file: fig12_Pip14b.tex
\pspicture*(4.5,4.5)
%
%\rput[tl](.4,3.9){\small$q_1$}
\pscircle(2.25,2.25){2}
\pscircle*[linecolor=lightgray](2.75,1.75){.3}
\pscircle(2.75,1.75){.3}
\pscircle*(2.75,3){.075}
\pscircle*(1.25,3){.075}
\pscircle*(2,3.85){.075}
\psarc[linewidth=.5pt](2.75,3){.4}{0}{180}
\psline[linewidth=.5pt](2.75,2.4)(2.35,3)
\psline[linewidth=.5pt, arrowsize=4pt 3]{->}(2.75,2.4)(3.15,3)
\psline[linewidth=.5pt](2.75,2.4)(2.75,2.05)
\psline[linewidth=.5pt](2,3.85)(3.75,3.3)(3.75,1.75)(3.05,1.75)
\psline[linewidth=.5pt](2.75,1.45)(2.75,1.25)(3.4,1.25)(3.4,3.65)(2,3.85)
\rput[t](3,1.1){\small$\beta_{34}$}
\psarc[linewidth=.5pt](1.25,3){.4}{0}{180}
\psline[linewidth=.5pt](1.25,2.4)(.85,3)
\psline[linewidth=.5pt, arrowsize=4pt 3]{->}(1.25,2.4)(1.65,3)
\psline[linewidth=.5pt](1.25,2.4)(1.25,1.75)(2.45,1.75)
\rput[t](1.5,1.6){\small$\beta_{24}$}
\psline[linewidth=.5pt,linestyle=dashed,dash=3pt 2pt](2.75,2.05)(2.75,1.45)
\psline[linewidth=.5pt,linestyle=dashed,dash=3pt 2pt](2.45,1.75)(3.05,1.75)
\endpspicture

%% file: fig13_U.tex
\pspicture*(7,3)
\rput[tl](.4,2.6){\small$c$}
%\pscircle(1.5,1.5){1.25}
\psarc[arrowsize=2pt 2.5]{->}(1.5,1.5){1.25}{0}{360}
\pscircle*[linecolor=lightgray](1,1.5){.25}
\pscircle(1,1.5){.25}
\psline[linewidth=.5pt,linestyle=dashed,dash=3pt 2pt](1,1.75)(1,1.25)
\pscircle*[linecolor=lightgray](2,1.5){.25}
\pscircle(2,1.5){.25}
\psline[linewidth=.5pt,linestyle=dashed,dash=3pt 2pt](2,1.75)(2,1.25)
\psline[linewidth=.5pt](1.5,2.75)(1.5,.25)
\psline[linewidth=.5pt,linearc=.2](1,1.75)(1,2)(2,2)(2,1.75)
\psline[linewidth=.5pt,linearc=.2](1,1.25)(1,1)(2,1)(2,1.25)
\psline[linewidth=.5pt, arrowsize=2pt 2.5]{->}(1.7,1)(1.85,1)
\psline[linewidth=.4pt]{<-}(1.15,.9)(1.4,.9)
\psline[linewidth=.4pt]{->}(1.15,1.1)(1.4,1.1)
\psline[linewidth=.4pt]{->}(1.6,1.9)(1.85,1.9)
\psline[linewidth=.4pt]{<-}(1.6,2.1)(1.85,2.1)
\rput[tl](1.55,2.6){\small$d$}
\rput[t](2,.9){\small$a_1$}
%
%
%\psline{->}(3,1.5)(3.5,1.5)
%\rput[b](3.25,1.7){\small$U$}
%
%
\pscircle(5,1.5){1.25}
\pscircle*[linecolor=lightgray](4.5,1.5){.25}
\pscircle(4.5,1.5){.25}
%\psline(4.35,1.65)(4.65,1.35)
%\psline(4.65,1.65)(4.35,1.35)
%
\pscircle*[linecolor=lightgray](5.5,1.5){.25}
\pscircle(5.5,1.5){.25}
\psbezier[linewidth=.5pt](4,1.5)(4,.75)(5,.75)(5,1.5)
\psbezier[linewidth=.5pt](5,1.5)(5,2.25)(6,2.25)(6,1.5)
\psbezier[linewidth=.5pt](4,1.5)(4,2)(5,2.5)(5,2.75)
\psbezier[linewidth=.5pt](6,1.5)(6,1)(5,.5)(5,.25)
\rput[tl](4.9,2.5){\small$U(d)$}
%
%
%\psline{->}(6.5,1.5)(7,1.5)
%\rput[b](6.75,1.7){\small$A_1$}
%
%
\endpspicture

%% file: fig14_lanK1.tex
\pspicture*(3.75,3.2)
\pscircle(2,1.6){1.5}\rput(.65,.45){\small$c_1$}
\pscircle(2,2.4){.25}\rput(2,2.4){\small$c_2$}
\pscircle*[linecolor=lightgray](1.4,1.2){.3}
\pscircle(1.4,1.2){.3}
\pscircle*[linecolor=lightgray](2.6,1.2){.3}
\pscircle(2.6,1.2){.3}
\psbezier[linewidth=.5pt](1.15,.7)(.65,.95)(1.75,3.15)(2.25,2.9)
\psbezier[linewidth=.5pt](1.15,.7)(1.65,.45)(2.75,2.65)(2.25,2.9)
\rput[r](1.3,2.1){\small$d_1$}
\psbezier[linewidth=.5pt](2.85,.7)(2.35,.45)(1.25,2.65)(1.75,2.9)
\psbezier[linewidth=.5pt](2.85,.7)(3.35,.95)(2.25,3.15)(1.75,2.9)
\rput[l](2.7,2.1){\small$d_2$}
\psellipse[linewidth=.5pt](2,1.2)(1.2,.5)
\psline[linewidth=.5pt, arrowsize=2pt 2.5]{->}(2,.7)(2.05,.7)
\rput[t](2.05,.6){\small$c$}
\endpspicture

%% file: fig15_lanK2.tex
\pspicture*(3.75,3.2)
\pscircle(2,1.6){1.5}\rput(.65,.45){\small$c_1$}
\pscircle(2,2.4){.25}\rput(2,2.4){\small$c_2$}
\pscircle*[linecolor=lightgray](1.4,1.2){.3}
\pscircle(1.4,1.2){.3}
\pscircle*[linecolor=lightgray](2.6,1.2){.3}
\pscircle(2.6,1.2){.3}
\psline[linewidth=.5pt, arrowsize=2pt 2.5]{->}(2,.75)(2.15,.75)
\psline[linewidth=.5pt,linearc=.2](1.4,1.5)(1.4,1.65)(2.6,1.65)(2.6,1.5)
\psline[linewidth=.5pt,linearc=.2](1.4,.9)(1.4,.75)(2.6,.75)(2.6,.9)
\psline[linewidth=.5pt,linestyle=dashed,dash=3pt 2pt](1.4,.9)(1.4,1.5)
\psline[linewidth=.5pt,linestyle=dashed,dash=3pt 2pt](2.6,.9)(2.6,1.5)
\rput[t](1.85,.65){\small$a_1$}
\psline[linewidth=.4pt]{<-}(1.85,1.75)(2.15,1.75)
\psline[linewidth=.4pt]{->}(1.85,1.55)(2.15,1.55)
\rput[lt](1.1,2.2){\small$a_2$}
\psbezier[linewidth=.5pt](.75,1.2)(.75,2)(1.5,2.9)(2,2.9)
\psbezier[linewidth=.5pt](3.25,1.2)(3.25,2)(2.5,2.9)(2,2.9)
\psline[linewidth=.5pt](.75,1.2)(1.1,1.2)
\psline[linewidth=.5pt,linestyle=dashed,dash=3pt 2pt](1.1,1.2)(1.7,1.2)
\psline[linewidth=.5pt](1.7,1.2)(2.3,1.2)
\psline[linewidth=.5pt,linestyle=dashed,dash=3pt 2pt](2.3,1.2)(2.9,1.2)
\psline[linewidth=.5pt](2.9,1.2)(3.25,1.2)
\psline[linewidth=.4pt]{<-}(2.85,2.1)(2.65,2.35)
\psline[linewidth=.4pt]{->}(3.05,2.2)(2.85,2.5)
\endpspicture

%% file: fig16_lanK3.tex
\pspicture*(3.75,3.2)
\pscircle(2,1.6){1.5}\rput(.65,.45){\small$c'_1$}
\pscircle(2,2.4){.25}\rput(2,2.4){\small$c'_2$}
%\pscircle(1.4,1.2){.3}\rput(1.4,1.2){\small$a'_1$}
%\pscircle(2.6,1.2){.3}\rput(2.6,1.2){\small$a''_1$}
\psarc[arrowsize=2pt 2.5]{->}(1.4,1.2){.3}{110}{470}\rput(1.4,1.2){\small$a'_1$}
\psarc[arrowsize=2pt 2.5]{->}(2.6,1.2){.3}{110}{470}\rput(2.6,1.2){\small$a''_1$}
\psbezier[linewidth=.5pt](1.15,.7)(.65,.95)(1.75,3.15)(2.25,2.9)
\psbezier[linewidth=.5pt](1.15,.7)(1.65,.45)(2.75,2.65)(2.25,2.9)
\rput[r](1.3,2.1){\small$b$}
\psbezier[linewidth=.5pt](2.85,.7)(2.35,.45)(1.25,2.65)(1.75,2.9)
\psbezier[linewidth=.5pt](2.85,.7)(3.35,.95)(2.25,3.15)(1.75,2.9)
\rput[l](2.7,2.1){\small$a'_2$}
\psellipse[linewidth=.5pt](2,1.2)(1.2,.5)
\rput[t](2.05,.6){\small$c'$}
%\psline[linewidth=.5pt, arrowsize=2pt 2.5]{<-}(1.95,.7)(2,.7)
%
\endpspicture

%% file: fig17_Pi22.tex
\pspicture*(4.5,4.5)
\rput[tl](.4,3.9){\small$c_1$}
\pscircle(2.25,2.25){2}
\pscircle*[linecolor=lightgray](3,1.75){.3}
\pscircle(3,1.75){.3}
\pscircle*[linecolor=lightgray](1.5,1.75){.3}
\pscircle(1.5,1.75){.3}
%
%\pscircle*(2.75,3.25){.075}\rput[bl](2.75,3.4){\small$p_3$}
\pscircle*(1.75,3.25){.075}\rput[bl](1.75,3.4){\small$p_2$}
\rput[l](2.75,3){\small$x_2$}
\psline[linewidth=.5pt, arrowsize=4pt 3]{->}(3,2.75)(2.5,1.75)
\psline[linewidth=.5pt](2.5,1.75)(2.7,1.75)
\psline[linewidth=.5pt](1.75,3.25)(3,2.75)(3.5,1.75)(3.3,1.75)
\psline[linewidth=.5pt,linestyle=dashed,dash=3pt 2pt](2.7,1.75)(3.3,1.75)
\rput[r](1.3,2.75){\small$x_1$}
\psline[linewidth=.5pt, arrowsize=4pt 3]{->}(1.5,2.75)(1,1.75)
\psline[linewidth=.5pt](1,1.75)(1.2,1.75)
\psline[linewidth=.5pt](1.75,3.25)(1.5,2.75)(2,1.75)(1.8,1.75)
\psline[linewidth=.5pt,linestyle=dashed,dash=3pt 2pt](1.2,1.75)(1.8,1.75)
\endpspicture

%% file: fig18_Pip22.tex
\pspicture*(4.5,4.5)
\rput[tl](.4,3.9){\small$c_1$}
\pscircle(2.25,2.25){2}
\pscircle*[linecolor=lightgray](3,1.75){.3}
\pscircle(3,1.75){.3}
\pscircle*[linecolor=lightgray](1.5,1.75){.3}
\pscircle(1.5,1.75){.3}
%
%\pscircle*(2.75,3.25){.075}\rput[bl](2.75,3.4){\small$p_3$}
\pscircle*(1.75,3.25){.075}\rput[bl](1.75,3.4){\small$p_2$}
%
%\psarc[linewidth=.5pt](2.25,2.25){1.75}{180}{450}
\rput[bl](.75,2.75){\small$\gamma$}
\psarc[linewidth=.5pt](2.25,2.25){1.75}{180}{360}
\psline[linewidth=.5pt, arrowsize=4pt 3]{->}(1.75,3.25)(.5,2.25)
\psline[linewidth=.5pt](1.75,3.25)(4,2.25)
%\psline[linewidth=.5pt](1.75,3.25)(2.25,4)
%
\rput[bl](3.1,2.95){\small$\delta_2$}
\psarc[linewidth=.5pt](3,1.75){.65}{180}{360}
\psline[linewidth=.5pt, arrowsize=4pt 3]{->}(3,2.95)(2.35,1.75)
\psline[linewidth=.5pt](3,2.95)(3.65,1.75)
\psline[linewidth=.5pt](1.75,3.25)(3,2.95)
%\pscircle*(2,3.85){.075}
\rput[tl](1.6,1.1){\small$\alpha_2$}
\psline[linewidth=.5pt](3,2.05)(3,2.25)
\psline[linewidth=.5pt,linestyle=dashed,dash=3pt 2pt](3,2.05)(3,1.45)
\psline[linewidth=.5pt](3,1.45)(3,1.25)(1.5,1.25)(1.5,1.45)
\psline[linewidth=.5pt](1.5,2.05)(1.5,2.25)
\psline[linewidth=.5pt,linestyle=dashed,dash=3pt 2pt](1.5,2.05)(1.5,1.45)
\psline[linewidth=.5pt, arrowsize=4pt 3]{->}(1.75,3.25)(1.5,2.25)
\psline[linewidth=.5pt](1.75,3.25)(3,2.25)
%\psline[linewidth=.5pt](1.5,1.45)(1.5,1.25)

%\rput[l](2.15,3.85){\small$p_4$}
%\rput[bl](2.8,3.1){\small$p_3$}
%\rput[t](1.25,2.85){\small$p_2$}
%
%\psline[linewidth=.5pt, arrowsize=4pt 3]{->}(2.75,3)(2.1,1.5)
%\psarc[linewidth=.5pt](2.75,1.5){.65}{180}{360}
%\psline[linewidth=.5pt,linestyle=dashed,dash=3pt 2pt](2.75,1.45)(2.75,2.05)
%\psbezier[linewidth=.5pt](2.75,3)(1.5,.5)(4,.5)(2.75,3)
\endpspicture

%% file: fig19_Pi23.tex
\pspicture*(4.5,4.5)
\rput[tl](.4,3.9){\small$c_1$}
\pscircle(2.25,2.25){2}
\pscircle*[linecolor=lightgray](3,1.75){.3}
\pscircle(3,1.75){.3}
\pscircle*[linecolor=lightgray](1.5,1.75){.3}
\pscircle(1.5,1.75){.3}
\pscircle*(2.75,3.25){.075}\rput[bl](2.75,3.4){\small$p_3$}
\pscircle*(1.75,3.25){.075}%\rput[bl](1.75,3.4){\small$p_2$}
\rput[tr](1.35,3.35){\small$\delta_{23}$}
\psarc[linewidth=.5pt](1.75,3.25){.35}{90}{270}
\psline[linewidth=.5pt, arrowsize=4pt 3]{->}(1.75,2.9)(2.75,3.25)(1.75,3.6)
\rput[l](3.2,2.6){\small$y_2$}
\psline[linewidth=.5pt, arrowsize=4pt 3]{->}(3,2.75)(2.5,1.75)
\psline[linewidth=.5pt](2.5,1.75)(2.7,1.75)
\psline[linewidth=.5pt](2.75,3.25)(3,2.75)(3.5,1.75)(3.3,1.75)
\psline[linewidth=.5pt,linestyle=dashed,dash=3pt 2pt](2.7,1.75)(3.3,1.75)
\rput[r](1.5,2.5){\small$y_1$}
\psline[linewidth=.5pt, arrowsize=4pt 3]{->}(1.95,2.75)(.8,1.75)
\psline[linewidth=.5pt](.8,1.75)(1.2,1.75)
\psline[linewidth=.5pt](2.75,3.25)(1.95,2.75)(2,1.75)(1.8,1.75)
\psline[linewidth=.5pt,linestyle=dashed,dash=3pt 2pt](1.2,1.75)(1.8,1.75)
\endpspicture

%% file: fig20_Pip23.tex
\pspicture*(4.5,4.5)
\rput[tl](.4,3.9){\small$c_1$}
\pscircle(2.25,2.25){2}
\pscircle*[linecolor=lightgray](3,1.75){.3}
\pscircle(3,1.75){.3}
\pscircle*[linecolor=lightgray](1.5,1.75){.3}
\pscircle(1.5,1.75){.3}
\pscircle*(2.75,3.25){.075}%\rput[bl](2.75,3.4){\small$p_3$}
\pscircle*(1.75,3.25){.075}%\rput[bl](1.75,3.4){\small$p_2$}
\rput[br](2.45,2.3){\small$\delta_3$}
\psarc[linewidth=.5pt](3,1.8){.65}{180}{360}
\psline[linewidth=.5pt, arrowsize=4pt 3]{->}(2.75,3.25)(2.35,1.8)
%\psline[linewidth=.5pt](2.43,2.5)(2.36,1.75)
\psline[linewidth=.5pt](2.75,3.25)(3.65,1.8)
%\psline[linewidth=.5pt](1.75,3.25)(3,2.95)
%
%\rput[tl](1.6,1.1){\small$\alpha_2$}
%\psline[linewidth=.5pt](2.75,3.25)
%\psline[linewidth=.5pt,linestyle=dashed,dash=3pt 2pt](3,2.05)(3,1.45)
%
\rput[tr](1.5,3,25){\small$\varepsilon$}
\psarc[linewidth=.5pt](1.75,3.25){.25}{0}{180}
\psline[linewidth=.5pt, arrowsize=4pt 3]{->}(1.5,3.25)(1.75,2.3)(2.0,3.25)
\psline[linewidth=.5pt](1.75,2.3)(2.2,1.75)(2.7,1.75)
\psline[linewidth=.5pt](3.3,1.75)(3.85,2.85)(2.75,3.25)
%
%\rput[tr](1.35,3){\small$\delta_{23}$}
%\psarc[linewidth=.5pt](1.75,3.25){.5}{90}{270}
%\psline[linewidth=.5pt, arrowsize=4pt 3]{->}(1.75,2.75)(2.75,3.25)(1.75,3.75)
%
\rput[bl](2.5,3.85){\small$\alpha_{3}$}
\psarc[linewidth=.5pt](2.25,2.25){1.75}{90}{180}
\psline[linewidth=.5pt,linestyle=dashed,dash=3pt 2pt](3,2.05)(3,1.45)
\psline[linewidth=.5pt,linestyle=dashed,dash=3pt 2pt](1.5,2.05)(1.5,1.45)
\psline[linewidth=.5pt](1.5,1.45)(1.5,1.25)(3,1.25)(3,1.45)
\psline[linewidth=.5pt](2.75,3.25)(3,2.05)
\psline[linewidth=.5pt](.5,2.25)(1.5,2.25)(1.5,2.05)
\psline[linewidth=.5pt, arrowsize=4pt 3]{->}(2.75,3.25)(2.25,4)
\psline[linewidth=.5pt,linestyle=dashed,dash=3pt 2pt](3.3,1.75)(2.7,1.75)
\rput[b](2.25,.55){\small$\delta_{12}$}
\psline[linewidth=.5pt, arrowsize=4pt 3]{->}(2.75,3.25)(.5,1.75)
\psline[linewidth=.5pt](2.75,3.25)(4.05,2.25)
\psarc[linewidth=.5pt](2.25,2.25){1.8}{195}{360}

\endpspicture

%% file: fig21_torus.tex
\pspicture*(4,4)
\pscircle(2,3.25){.25}\rput(2,3.25){\small$c_3$}
\pscircle(1.05,1.15){.25}\rput(1.05,1.15){\small$c_1$}
\pscircle(2.95,1.15){.25}\rput(2.95,1.15){\small$c_2$}
%\rput[tl](.4,2.6){\small$c$}
\pscircle(2,2){1.75}
\pscircle(2,2){.35}
\rput[lb](2.4,2.65){\small$b$}
\pscircle[linewidth=.5pt](2,2){.75}
\psbezier[linewidth=.5pt](.25,2)(.25,1.75)(1.65,1.75)(1.65,2)
\psbezier[linewidth=.5pt, linestyle=dashed,dash=3pt 2pt](.25,2)(.25,2.25)(1.65,2.25)(1.65,2)
\rput[b](.85,2.25){\small$a_1$}
\psbezier[linewidth=.5pt](2.35,2)(2.35,1.75)(3.75,1.75)(3.75,2)
\psbezier[linewidth=.5pt, linestyle=dashed,dash=3pt 2pt](2.35,2)(2.35,2.25)(3.75,2.25)(3.75,2)
\rput[b](3.15,2.25){\small$a_3$}
\psbezier[linewidth=.5pt](2,1.65)(2.25,1.65)(2.25,.25)(2,.25)
\psbezier[linewidth=.5pt, linestyle=dashed,dash=3pt 2pt](2,1.65)(1.75,1.65)(1.75,.25)(2,.25)
\rput[l](2.25,.75){\small$a_2$}
%\rput[b](4.5,2.5){$y'$}
%\psline{->}(4.25,2.25)(4.75,2.25)
%
%
%\psbezier(7,4)(7,3.5)(5.5,3)(5.5,2.25)
\endpspicture

%% file: fig22_F32.tex
\pspicture*(5,4)
\psellipse(2.5,2)(2.25,1.75)
\pscircle*[linecolor=lightgray](2.5,1.5){.3}
\pscircle(2.5,1.5){.3}
\pscircle*[linecolor=lightgray](1.25,1.5){.3}
\pscircle(1.25,1.5){.3}
\pscircle*[linecolor=lightgray](3.75,1.5){.3}
\pscircle(3.75,1.5){.3}
\pscircle(1.85,2.8){.25}\rput(1.85,2.8){\small$c_2$}
\pscircle(3,2.8){.25}\rput(3,2.8){\small$c_1$}
\rput[tr](2.05,1){\small$c$}
\psellipse[linewidth=.5pt](3.125,1.6)(1.3,.85)
\psline[linewidth=.5pt, arrowsize=2pt 2.5]{->}(3.2,.75)(3.25,.75)
\rput[t](3.45,2.05){\small$a_1$}
\psline[linewidth=.4pt]{->}(2.9,2.05)(3.2,2.05)
\psline[linewidth=.4pt]{<-}(2.9,2.25)(3.2,2.25)
\psline[linewidth=.5pt,linearc=.2](2.5,1.8)(2.5,2.15)(3.75,2.15)(3.75,1.8)
\psline[linewidth=.5pt,linestyle=dashed,dash=3pt 2pt](2.5,1.8)(2.5,1.2)
\psline[linewidth=.5pt,linestyle=dashed,dash=3pt 2pt](3.75,1.8)(3.75,1.2)
\psline[linewidth=.5pt,linearc=.2](2.5,1.2)(2.5,.95)(3.75,.95)(3.75,1.2)
\psline[linewidth=.5pt, arrowsize=2pt 2.5]{->}(3.25,.95)(3.3,.95)
\rput[b](1.05,2.05){\small$b$}
\psline[linewidth=.4pt]{->}(1.4,1.9)(1.7,1.9)
\psline[linewidth=.4pt]{<-}(1.4,2.1)(1.7,2.1)
\psline[linewidth=.5pt,linearc=.2](2.8,1.5)(2.9,1.5)(2.9,2)(.75,2)(.75,1.5)(.95,1.5)
\psline[linewidth=.5pt,linestyle=dashed,dash=3pt 2pt](.95,1.5)(1.55,1.5)
\psline[linewidth=.5pt](1.55,1.5)(2.2,1.5)
\psline[linewidth=.5pt,linestyle=dashed,dash=3pt 2pt](2.2,1.5)(2.8,1.5)
\psline[linewidth=.5pt,linestyle=dashed,dash=3pt 2pt](1,1.25)(1.5,1.75)
\psline[linewidth=.5pt](1.5,1.75)(2.25,1.75)
\psline[linewidth=.5pt,linestyle=dashed,dash=3pt 2pt](2.25,1.75)(2.75,1.25)
\psline[linewidth=.5pt](2.75,1.25)(3.5,1.25)
\psline[linewidth=.5pt,linestyle=dashed,dash=3pt 2pt](3.5,1.25)(4,1.75)
%
%\psline[linewidth=.5pt](1,1.25)(.85,1.1)
\psline[linewidth=.5pt](4,1.75)(4.5,2.25)
\psbezier[linewidth=.5pt](.5,2.25)(.25,1.75)(.85,1.1)(1,1.25)
\psbezier[linewidth=.5pt](.5,2.25)(1,3.25)(2,3.5)(2.5,3.5)
\psbezier[linewidth=.5pt](4.5,2.25)(4,3.25)(3,3.5)(2.5,3.5)
\rput(3.95,2.65){\small$e$}
\endpspicture

%% file: fig23_F32a.tex
\pspicture*(5,4)
\psellipse(2.5,2)(2.25,1.75)
\pscircle*[linecolor=lightgray](2.5,1.5){.3}
\pscircle(2.5,1.5){.3}
\pscircle*[linecolor=lightgray](1.25,1.5){.3}
\pscircle(1.25,1.5){.3}
\pscircle*[linecolor=lightgray](3.75,1.5){.3}
\pscircle(3.75,1.5){.3}
\pscircle(1.85,2.8){.25}\rput(1.85,2.8){\small$c_2$}
\pscircle(3,2.8){.25}\rput(3,2.8){\small$c_1$}
%
%\rput[tr](2.05,1){\small$c$}
%\psellipse[linewidth=.5pt](3.125,1.6)(1.3,.85)
%\psline[linewidth=.5pt, arrowsize=2pt 2.5]{->}(3.2,.65)(3.25,.65)
%
%\rput[t](3.45,2.05){\small$a_1$}
%\psline[linewidth=.4pt]{->}(2.9,2.05)(3.2,2.05)
%\psline[linewidth=.4pt]{<-}(2.9,2.25)(3.2,2.25)
%\psline[linewidth=.5pt,linearc=.2](2.5,1.8)(2.5,2.15)(3.75,2.15)(3.75,1.8)
%\psline[linewidth=.5pt,linestyle=dashed,dash=3pt 2pt](2.5,1.8)(2.5,1.2)
%\psline[linewidth=.5pt,linestyle=dashed,dash=3pt 2pt](3.75,1.8)(3.75,1.2)
%\psline[linewidth=.5pt,linearc=.2](2.5,1.2)(2.5,.95)(3.75,.95)(3.75,1.2)
%\psline[linewidth=.5pt, arrowsize=2pt 2.5]{->}(3.25,.95)(3.3,.95)
%
%\rput[b](1.05,2.05){\small$b$}
%\psline[linewidth=.4pt]{->}(1.4,1.9)(1.7,1.9)
%\psline[linewidth=.4pt]{<-}(1.4,2.1)(1.7,2.1)
%\psline[linewidth=.5pt,linearc=.2](2.8,1.5)(2.9,1.5)(2.9,2)(.75,2)(.75,1.5)(.95,1.5)
%\psline[linewidth=.5pt,linestyle=dashed,dash=3pt 2pt](.95,1.5)(1.55,1.5)
%\psline[linewidth=.5pt](1.55,1.5)(2.2,1.5)
%\psline[linewidth=.5pt,linestyle=dashed,dash=3pt 2pt](2.2,1.5)(2.8,1.5)
%
%
%
\psline[linewidth=.5pt,linestyle=dashed,dash=3pt 2pt](2.25,1.25)(2.75,1.75)
\psline[linewidth=.5pt](2.75,1.75)(3.5,1.75)
\psline[linewidth=.5pt,linestyle=dashed,dash=3pt 2pt](2.25,1.75)(2.75,1.25)
\psline[linewidth=.5pt](2.75,1.25)(3.5,1.25)
\psline[linewidth=.5pt,linestyle=dashed,dash=3pt 2pt](3.5,1.25)(4,1.75)
\psline[linewidth=.5pt,linestyle=dashed,dash=3pt 2pt](3.5,1.75)(4,1.25)
%
%\psline[linewidth=.5pt](1,1.25)(.85,1.1)
\psline[linewidth=.5pt,linearc=.35](4,1.75)(4.15,1.9)(3.3,3.25)(3,3.25)(2.5,3.25)(2.15,1.85)(2.25,1.75)
\psline[linewidth=.4pt]{->}(2.3,2.1)(2.4,2.45)
\psline[linewidth=.4pt]{<-}(2.1,2.1)(2.2,2.45)
\rput[r](3.85,2.25){\small$a_2$}
\psline[linewidth=.5pt,linearc=.4](2.25,1.25)(2,1.25)(1.3,2.5)(1.35,3.4)(3.5,3.4)(4.4,2.5)(4.4,1.25)(4,1.25)
\psline[linewidth=.4pt]{->}(1.95,3.3)(2.3,3.3)
\psline[linewidth=.4pt]{<-}(1.95,3.5)(2.3,3.5)
\rput[r](1.25,2.7){\small$a_3$}
%\psbezier[linewidth=.5pt](.5,2.25)(.25,1.75)(.85,1.1)(1,1.25)
%\psbezier[linewidth=.5pt](.5,2.25)(1,3.25)(2,3.5)(2.5,3.5)
%\psbezier[linewidth=.5pt](4.5,2.25)(4,3.25)(3,3.5)(2.5,3.5)
%\rput(3.95,2.65){\small$e$}
%
%
\endpspicture

%% file: fig24_Pi31.tex
\pspicture*(5,4)
\psellipse(2.5,2)(2.25,1.75)
\pscircle*[linecolor=lightgray](2.5,1.5){.3}
\pscircle(2.5,1.5){.3}
\pscircle*[linecolor=lightgray](1.25,1.5){.3}
\pscircle(1.25,1.5){.3}
\pscircle*[linecolor=lightgray](3.75,1.5){.3}
\pscircle(3.75,1.5){.3}
%
%
%\pscircle*(2,2.75){.075}%\rput[br](1.95,2.9){\small$p_2$}
\pscircle*(3,2.75){.075}\rput[bl](3.05,2.9){\small$p_1$}
%
%\rput[bl](3.85,2.25){\small$c$}
%\psellipse[linewidth=.5pt](3.125,1.5)(1.3,.85)
%
%\rput[t](3.125,.7){\small$\alpha_1$}
%\psline[linewidth=.5pt](3,2.75)(3.125,2.1)
%\psline[linewidth=.5pt](2.5,1.8)(2.5,2.1)(3.75,2.1)(3.75,1.8)
%\psline[linewidth=.5pt,linestyle=dashed,dash=3pt 2pt](2.5,1.8)(2.5,1.2)
%\psline[linewidth=.5pt,linestyle=dashed,dash=3pt 2pt](3.75,1.8)(3.75,1.2)
%\psline[linewidth=.5pt, arrowsize=4pt 3]{->}(2.5,1.2)(2.5,.8)
%\psline[linewidth=.5pt](2.5,.8)(3.75,.8)(3.75,1.2)
%
%\rput[br](2.2,2.2){\small$\beta_1$}
%\psline[linewidth=.5pt](3,2.75)(1.9,1.95)
%\psline[linewidth=.5pt, arrowsize=4pt 3]{->}(2.8,1.5)(3,1.5)(3,1.95)(.55,1.95)
%\psline[linewidth=.5pt](.55,1.95)(.55,1.5)(.95,1.5)
%\psline[linewidth=.5pt,linestyle=dashed,dash=3pt 2pt](.95,1.5)(1.55,1.5)
%\psline[linewidth=.5pt](1.55,1.5)(2.2,1.5)
%\psline[linewidth=.5pt,linestyle=dashed,dash=3pt 2pt](2.2,1.5)(2.8,1.5)
%\psline[linewidth=.5pt]
\rput[r](1.05,2.25){\small$x_1$}
\psline[linewidth=.5pt, arrowsize=4pt 3]{<-}(.75,1.5)(1.25,2.4)
\psline[linewidth=.5pt](1.75,1.5)(1.25,2.4)
\psline[linewidth=.5pt](3,2.75)(1.25,2.4)
\psline[linewidth=.5pt](1.55,1.5)(1.75,1.5)
\psline[linewidth=.5pt,linestyle=dashed,dash=3pt 2pt](.95,1.5)(1.55,1.5)
\psline[linewidth=.5pt](.95,1.5)(.75,1.5)
\rput[r](2.3,2.25){\small$x_2$}
\psline[linewidth=.5pt, arrowsize=4pt 3]{<-}(2,1.5)(2.5,2.4)
\psline[linewidth=.5pt](4.25,1.5)(3.75,2.4)
\psline[linewidth=.5pt](3,2.75)(3.75,2.4)
\psline[linewidth=.5pt](2.8,1.5)(3,1.5)
\psline[linewidth=.5pt,linestyle=dashed,dash=3pt 2pt](2.2,1.5)(2.8,1.5)
\psline[linewidth=.5pt](2.2,1.5)(2,1.5)
\rput[r](3.55,2.25){\small$x_3$}
\psline[linewidth=.5pt, arrowsize=4pt 3]{<-}(3.25,1.5)(3.75,2.4)
\psline[linewidth=.5pt](3,1.5)(2.5,2.4)
\psline[linewidth=.5pt](3,2.75)(2.5,2.4)
\psline[linewidth=.5pt](4.05,1.5)(4.25,1.5)
\psline[linewidth=.5pt,linestyle=dashed,dash=3pt 2pt](3.45,1.5)(4.05,1.5)
\psline[linewidth=.5pt](3.45,1.5)(3.25,1.5)
%\psline[linewidth=.5pt, arrowsize=4pt 3]{->}(2.75,3)(2.1,1.5)
%\psarc[linewidth=.5pt](2.75,1.5){.65}{180}{360}
%\psline[linewidth=.5pt,linestyle=dashed,dash=3pt 2pt](2.75,1.45)(2.75,2.05)
%\psbezier[linewidth=.5pt](2.75,3)(1.5,.5)(4,.5)(2.75,3)
\endpspicture

%% file: fig25_Pip31.tex
\pspicture*(5,4)
\psellipse(2.5,2)(2.25,1.75)
\pscircle*[linecolor=lightgray](2.5,1.5){.3}
\pscircle(2.5,1.5){.3}
\pscircle*[linecolor=lightgray](1.25,1.5){.3}
\pscircle(1.25,1.5){.3}
\pscircle*[linecolor=lightgray](3.75,1.5){.3}
\pscircle(3.75,1.5){.3}
%
%
%\pscircle*(2,2.75){.075}%\rput[br](1.95,2.9){\small$p_2$}
\pscircle*(3,2.75){.075}%\rput[bl](3.05,2.9){\small$p_1$}
%
%\rput[bl](3.85,2.25){\small$c$}
%\psellipse[linewidth=.5pt](3.125,1.5)(1.3,.85)
%
%\rput[t](3.125,.7){\small$\alpha_1$}
%\psline[linewidth=.5pt](3,2.75)(3.125,2.1)
%\psline[linewidth=.5pt](2.5,1.8)(2.5,2.1)(3.75,2.1)(3.75,1.8)
%\psline[linewidth=.5pt,linestyle=dashed,dash=3pt 2pt](2.5,1.8)(2.5,1.2)
%\psline[linewidth=.5pt,linestyle=dashed,dash=3pt 2pt](3.75,1.8)(3.75,1.2)
%\psline[linewidth=.5pt, arrowsize=4pt 3]{->}(2.5,1.2)(2.5,.8)
%\psline[linewidth=.5pt](2.5,.8)(3.75,.8)(3.75,1.2)
%
\rput[b](3.75,2.75){\small$\alpha_1$}
\psline[linewidth=.5pt, arrowsize=4pt 3]{->}(3,2.75)(2,1.8)
\psline[linewidth=.5pt,linestyle=dashed,dash=3pt 2pt](2.2,1.5)(2.8,1.5)
\psline[linewidth=.5pt,linestyle=dashed,dash=3pt 2pt](3.45,1.5)(4.05,1.5)
\psline[linewidth=.5pt](2,1.8)(2,1.5)(2.2,1.5)
\psline[linewidth=.5pt](2.8,1.5)(3.45,1.5)
\psline[linewidth=.5pt](4.05,1.5)(4.45,1.5)(4.45,2.5)(3,2.75)
%\rput[br](2.2,2.2){\small$\beta_1$}
%\psline[linewidth=.5pt](3,2.75)(1.9,1.95)
%\psline[linewidth=.5pt, arrowsize=4pt 3]{->}(2.8,1.5)(3,1.5)(3,1.95)(.55,1.95)
%\psline[linewidth=.5pt](.55,1.95)(.55,1.5)(.95,1.5)
%\psline[linewidth=.5pt,linestyle=dashed,dash=3pt 2pt](.95,1.5)(1.55,1.5)
%\psline[linewidth=.5pt](1.55,1.5)(2.2,1.5)
%\psline[linewidth=.5pt,linestyle=dashed,dash=3pt 2pt](2.2,1.5)(2.8,1.5)
%\psline[linewidth=.5pt]
\rput[t](2.1,.7){\small$\beta_1$}
\psline[linewidth=.5pt,linestyle=dashed,dash=3pt 2pt](2.5,1.8)(2.5,1.2)
\psline[linewidth=.5pt,linestyle=dashed,dash=3pt 2pt](1.25,1.8)(1.25,1.2)
\psline[linewidth=.5pt](3,2.75)(2.5,1.9)
\psline[linewidth=.5pt](2.5,1.8)(2.5,1.9)
%\psline[linewidth=.5pt](3.75,2.1)(3.75,1.8)
\psline[linewidth=.5pt](2.5,1.2)(2.5,.8)
\psline[linewidth=.5pt, arrowsize=4pt 3]{<-}(2.5,.8)(1.25,.8)(1.25,1.2)
\psline[linewidth=.5pt](1.25,1.8)(1.25,1.9)(3,2.75)
\rput[bl](1.75,2.75){\small$\gamma$}
\psarc[linewidth=.5pt](1.25,1.4){.5}{180}{360}
\psline[linewidth=.5pt](.75,1.4)(1.25,2.5)
\psline[linewidth=.5pt, arrowsize=4pt 3]{->}(1.75,1.4)(1.25,2.5)
\psline[linewidth=.5pt](3,2.75)(1.25,2.5)
\rput[tr](3.45,1){\small$\delta$}
\psarc[linewidth=.5pt](3.75,1.4){.5}{180}{360}
\psline[linewidth=.5pt](3.25,1.4)(3.75,2.4)
\psline[linewidth=.5pt, arrowsize=4pt 3]{->}(4.25,1.4)(3.75,2.4)
\psline[linewidth=.5pt](3,2.75)(3.75,2.4)
%\psline[linewidth=.5pt, arrowsize=4pt 3]{->}(2.75,3)(2.1,1.5)
%\psarc[linewidth=.5pt](2.75,1.5){.65}{180}{360}
%\psline[linewidth=.5pt,linestyle=dashed,dash=3pt 2pt](2.75,1.45)(2.75,2.05)
%\psbezier[linewidth=.5pt](2.75,3)(1.5,.5)(4,.5)(2.75,3)
\endpspicture

%% file: fig26_Pi32.tex
\pspicture*(5,4)
\psellipse(2.5,2)(2.25,1.75)
\pscircle*[linecolor=lightgray](2.5,1.5){.3}
\pscircle(2.5,1.5){.3}
\pscircle*[linecolor=lightgray](1.25,1.5){.3}
\pscircle(1.25,1.5){.3}
\pscircle*[linecolor=lightgray](3.75,1.5){.3}
\pscircle(3.75,1.5){.3}
\pscircle*(2,2.75){.075}\rput[br](1.95,2.9){\small$p_2$}
\pscircle*(3,2.75){.075}%\rput[bl](3.05,2.9){\small$p_1$}
%
%\rput[bl](3.85,2.25){\small$c$}
%\psellipse[linewidth=.5pt](3.125,1.5)(1.3,.85)
%
\rput[l](3.2,3.1){\small$\delta_3$}
\psarc[linewidth=.5pt](2.9,2.7){.4}{270}{450}
\psline[linewidth=.5pt, arrowsize=4pt 3]{->}(2.9,3.1)(2,2.75)(2.9,2.3)
\rput[r](1.05,2.25){\small$y_1$}
\psline[linewidth=.5pt, arrowsize=4pt 3]{<-}(.75,1.5)(1.25,2.4)
\psline[linewidth=.5pt](1.75,1.5)(1.25,2.4)
\psline[linewidth=.5pt](2,2.75)(1.25,2.4)
\psline[linewidth=.5pt](1.55,1.5)(1.75,1.5)
\psline[linewidth=.5pt,linestyle=dashed,dash=3pt 2pt](.95,1.5)(1.55,1.5)
\psline[linewidth=.5pt](.95,1.5)(.75,1.5)
\rput[r](2.1,2.25){\small$y_2$}
\psline[linewidth=.5pt, arrowsize=4pt 3]{<-}(2,1.5)(2.2,2.4)
\psline[linewidth=.5pt](2.8,1.5)(3,1.5)
\psline[linewidth=.5pt](3,1.5)(2.2,2.4)
\psline[linewidth=.5pt](2,2.75)(2.2,2.4)
\psline[linewidth=.5pt,linestyle=dashed,dash=3pt 2pt](2.2,1.5)(2.8,1.5)
\psline[linewidth=.5pt](2.2,1.5)(2,1.5)
\rput[l](4,2.25){\small$y_3$}
\psline[linewidth=.5pt, arrowsize=4pt 3]{<-}(3.25,1.5)(3.75,2.4)
\psline[linewidth=.5pt](4.25,1.5)(3.75,2.4)
\psline[linewidth=.5pt](2,2.75)(3.75,2.4)
\psline[linewidth=.5pt](4.05,1.5)(4.25,1.5)
\psline[linewidth=.5pt,linestyle=dashed,dash=3pt 2pt](3.45,1.5)(4.05,1.5)
\psline[linewidth=.5pt](3.45,1.5)(3.25,1.5)
%\psline[linewidth=.5pt, arrowsize=4pt 3]{->}(2.75,3)(2.1,1.5)
%\psarc[linewidth=.5pt](2.75,1.5){.65}{180}{360}
%\psline[linewidth=.5pt,linestyle=dashed,dash=3pt 2pt](2.75,1.45)(2.75,2.05)
%\psbezier[linewidth=.5pt](2.75,3)(1.5,.5)(4,.5)(2.75,3)
\endpspicture

%% file: fig27_Pip32.tex
\pspicture*(5,4)
\psellipse(2.5,2)(2.25,1.75)
\pscircle*[linecolor=lightgray](2.5,1.5){.3}
\pscircle(2.5,1.5){.3}
\pscircle*[linecolor=lightgray](1.25,1.5){.3}
\pscircle(1.25,1.5){.3}
\pscircle*[linecolor=lightgray](3.75,1.5){.3}
\pscircle(3.75,1.5){.3}
\pscircle*(2,2.75){.075}%\rput[br](1.95,2.9){\small$p_2$}
\pscircle*(3,2.75){.075}%\rput[bl](3.05,2.9){\small$p_1$}
%
%\rput[bl](3.85,2.25){\small$c$}
%\psellipse[linewidth=.5pt](3.125,1.5)(1.3,.85)
%
\rput[t](2.1,.7){\small$\beta_2$}
\psline[linewidth=.5pt,linestyle=dashed,dash=3pt 2pt](2.5,1.8)(2.5,1.2)
\psline[linewidth=.5pt,linestyle=dashed,dash=3pt 2pt](1.25,1.8)(1.25,1.2)
\psline[linewidth=.5pt](2,2.75)(2.5,1.9)
\psline[linewidth=.5pt](2.5,1.8)(2.5,1.9)
%\psline[linewidth=.5pt](3.75,2.1)(3.75,1.8)
\psline[linewidth=.5pt](2.5,1.2)(2.5,.8)
\psline[linewidth=.5pt, arrowsize=4pt 3]{<-}(2.5,.8)(1.25,.8)(1.25,1.2)
\psline[linewidth=.5pt](1.25,1.8)(1.25,1.9)(2,2.75)
%\psline[linewidth=.5pt](2.5,.8)(3.75,.8)(3.75,1.2)
%
\rput[r](2.2,3.25){\small$\alpha_2$}
%\psline[linewidth=.5pt](2,2.75)(1.9,1.95)
\psline[linewidth=.5pt,linestyle=dashed,dash=3pt 2pt](2.2,1.5)(2.8,1.5)
\psline[linewidth=.5pt,linestyle=dashed,dash=3pt 2pt](3.45,1.5)(4.05,1.5)
\psline[linewidth=.5pt, arrowsize=4pt 3]{->}(2,2.75)(2,1.5)
\psline[linewidth=.5pt](2,1.5)(2.2,1.5)
\psline[linewidth=.5pt](2.8,1.5)(3.45,1.5)
\psline[linewidth=.5pt](4.05,1.5)(4.45,1.5)
\psarc[linewidth=.5pt](2.5,1.5){1.95}{0}{90}
\psline[linewidth=.5pt](2,2.75)(2.5,3.45)
%\psline[linewidth=.5pt, arrowsize=4pt 3]{->}(2.8,1.5)(3,1.5)(3,1.95)(.55,1.95)
%\psline[linewidth=.5pt](.55,1.95)(.55,1.5)(.95,1.5)
%\psline[linewidth=.5pt](1.55,1.5)(2.2,1.5)
%\psline[linewidth=.5pt]
%
\rput[br](1.25,2.8){\small$\delta_1$}
\psarc[linewidth=.5pt](1.25,1.4){.5}{180}{360}
\psline[linewidth=.5pt](.75,1.4)(1.25,2.75)
\psline[linewidth=.5pt, arrowsize=4pt 3]{->}(1.75,1.4)(1.25,2.75)
\psline[linewidth=.5pt](2,2.75)(1.25,2.75)
\rput[tr](3.5,2.4){\small$\delta_2$}
\psarc[linewidth=.5pt](3.75,1.4){.5}{180}{360}
\psline[linewidth=.5pt](3.25,1.4)(3.75,2.4)
\psline[linewidth=.5pt, arrowsize=4pt 3]{->}(4.25,1.4)(3.75,2.4)
\psline[linewidth=.5pt](2,2.75)(3.75,2.4)
%
%
%\rput[tl](2.8,2.2){\small$\delta_3$}
%\psarc[linewidth=.5pt](2.9,2.7){.4}{270}{450}
%\psline[linewidth=.5pt, arrowsize=4pt 3]{->}(2.9,3.1)(2,2.75)(2.9,2.3)
%
%\psline[linewidth=.5pt, arrowsize=4pt 3]{->}(2.75,3)(2.1,1.5)
%\psarc[linewidth=.5pt](2.75,1.5){.65}{180}{360}
%\psline[linewidth=.5pt,linestyle=dashed,dash=3pt 2pt](2.75,1.45)(2.75,2.05)
%\psbezier[linewidth=.5pt](2.75,3)(1.5,.5)(4,.5)(2.75,3)
\endpspicture